\documentclass[reqno,11pt]{amsart}

\usepackage{amsmath}
\usepackage{amsthm}
\usepackage{mathrsfs}
\usepackage{amsfonts}
\usepackage{parskip} 
\usepackage{hyperref}
\usepackage{epsfig}

\usepackage{amsmath,amssymb,amsfonts,amsthm,amstext,verbatim}
\usepackage{enumerate,color,graphicx,
}
\usepackage{psfrag}  
\usepackage{url}
\usepackage{mdwlist}   
\usepackage{mathrsfs}
\usepackage{hyperref}
\usepackage{epsfig}
\usepackage{etoolbox}

\makeatletter
\patchcmd{\@part}{\Huge}{\normalsize}{}{} 
\makeatother

\numberwithin{equation}{section}

\allowdisplaybreaks

\newtheorem{theorem}{Theorem}[section]

\newtheorem{proposition}[theorem]{Proposition}
\newtheorem{lemma}[theorem]{Lemma}

\newtheorem{definition}[theorem]{Definition}
\newtheorem{corollary}[theorem]{Corollary}

\newcommand{\ep}{\varepsilon}

\newcommand{\R}{\mathbb{R}}
\newcommand{\N}{\mathbb{N}}
\newcommand{\Z}{\mathbb{Z}}

\newcommand{\dist}{\vert\vert}

\newcommand{\Cgjthree}{C_4}

\newcommand{\eps}{\ep}

\newcommand{\psap}{\mathsf{P}}
\newcommand{\psapleft}[1]{\mathsf{P}_{#1}^{\mathsf{left}}}

\newcommand{\psaw}{\mathsf{W}}
\newcommand{\psawtd}{\mathsf{W}^{3,\mathsf{E}}}
\newcommand{\barp}{\hat{p}}
\newcommand{\carp}{\hat{c}}
\newcommand{\thethat}{\hat{\thet}}
\newcommand{\muhat}{\hat{\mu}}
\newcommand{\psaptd}{\mathsf{P}^{3,\mathsf{E}}}
\newcommand{\sawtd}{\mathrm{SAW}^{3,\mathsf{E}}}
\newcommand{\saptd}{\mathrm{SAP}^{3,\mathsf{E}}}

\newcommand{\esaw}{{\mathbb{E}_{\mathsf{W}}}}

\newcommand{\sapell}{\mathrm{SAP}^{\mathrm{left}}}
\newcommand{\sapellward}{\mathrm{SAP}^{\mathfrak{l}}}
\newcommand{\sapr}{\mathrm{SAP}^{\mathrm{right}}}
\newcommand{\saprward}{\mathrm{SAP}^{\mathfrak{r}}}

\newcommand{\fpart}{\mathrm{First}}

\newcommand{\reflect}{\mathcal{R}}

\newcommand{\lattv}{\Z^d}

\newcommand{\northeast}{\mathrm{NE}}
\newcommand{\eastnorth}{\mathrm{EN}}
\newcommand{\eastsouth}{\mathrm{ES}}
\newcommand{\westnorth}{\mathrm{WN}}
\newcommand{\southeast}{\mathrm{SE}}
\newcommand{\base}{\mathrm{SE}}
\newcommand{\rightv}{\mathrm{Right}}

\newcommand{\vas}{\varsigma}
\newcommand{\vaso}{\varsigma^{\rm loc}}
\newcommand{\gamin}{\Gamma^{{\rm in}}}
\newcommand{\gamout}{\Gamma^{{\rm out}}}
\newcommand{\gaemp}{\ga^{{\rm empty}}}
\newcommand{\ellempty}{l^{{\rm empty}}}

\newcommand{\alphamac}{\alpha'}
\newcommand{\deltamac}{\delta'}
\newcommand{\deltanew}{\varphi}

\newcommand{\kmac}{k^*}

\newcommand{\ZZ}{\mathbb{Z}}     

\newcommand\hPhi{\mathrm{High}\fpart}

\theoremstyle{remark}

\newcounter{mycount}

\def\mik{1}
\newcommand\cpsfrag[2]{\ifnum\mik=1\psfrag{#1}{#2}\fi}

\newcommand{\saw}{\mathrm{SAW}}
\newcommand{\rhssaw}{\mathrm{RHSSAW}}
\newcommand{\sawfree}{\mathrm{SAW}^{*}}
\newcommand{\sap}{\mathrm{SAP}}

\newcommand{\dubbub}{\mathrm{RGJ}}

\newcommand{\Cpcp}{\zeta}
\newcommand{\maceta}{\zeta}
\newcommand{\Cpoly}{C_1}
\newcommand{\Ch}{C_6}

\newcommand{\Cgjone}{C_2}
\newcommand{\Cgjtwo}{C_3}

\newcommand{\conindex}{C_8}

\newcommand{\calnew}{c}

\newcommand{\epn}{\eps}
\newcommand{\epnmac}{\eps}
\newcommand{\eroom}{\epsilon_1}
\newcommand{\eclo}{\epsilon_2}
\newcommand{\esm}{\eta}

\newcommand{\thet}{\theta}
\newcommand{\thetc}{\xi}
\newcommand{\hcb}{\highclose_{1/2 + \chi}}

\newcommand{\reg}{\mathsf{R}}
\newcommand{\regul}{\mathsf{K}}
\newcommand{\irr}{\mathsf{S}}
\newcommand{\regnew}[3]{\mathsf{R}_{#1,#2}^{#3}}
\newcommand{\indexc}{L}
\newcommand{\indexcj}{j}

\newcommand{\litnum}{\tfrac{1}{10}}

\newcommand{\closes}{\text{ closes}}

\newcommand\Ga{\Gamma}     
\newcommand\ga{\gamma}     
\newcommand\Ptmp{P}  
\newcommand\Wtmp{\rm W}  
\newcommand\Atmp{\rm A}  

\newcommand{\hidden}[1]{}

\newcommand\mcond{\, \middle| \,}
\newcommand\cond{\, | \,}

\newcommand\strongjoin{\mathsf{GlobalMJ}}
\newcommand\strongtdjoin{\mathsf{GlobalJoin}}
\newcommand\globaljoin{\mathsf{GJ}}
\newcommand\ymax{y_{{\rm max}}}
\newcommand\ymin{y_{{\rm min}}}
\newcommand\xmax{x_{{\rm max}}}
\newcommand\xmin{x_{{\rm min}}}

\newcommand\righttip{\eastsouth}

\newcommand\phicl[3]{\hPhi_{#1,#2}^{#3}}

\newcommand\reverse{\overleftarrow}

\newcommand\pexper{P_{{\rm res}}}
\newcommand\eexper{E_{{\rm res}}}
\newcommand\lmid{\mathrm{l}_{{\rm mid}}}
\newcommand\notick{\mathsf{NoCharm}}

\newcommand{\modify}[1]{#1_{{\rm mod}}}

\newcommand{\highclose}{\mathsf{HCP}}
\newcommand{\highpolynum}{\mathsf{HPN}}

\newcommand{\eone}{\theta^{(1)}}
\newcommand{\etwo}{\theta^{(2)}}
\newcommand{\ethree}{\theta^{(3)}}

\newcommand{\macess}{\mathsf{High}\Theta}

\newcommand{\lset}{K'}
\newcommand{\tilden}{m'}
\newcommand{\proj}{\mathsf{Proj}}

\newcommand{\leftv}{\mathrm{Left}}
\newcommand{\lrpp}{\mathsf{LeftRightPolyPair}^{3,\mathsf{E}}}
\newcommand{\xspan}{x_{{\rm span}}}

\newcommand{\targetmany}{\mathsf{TMHA}}

\title[Self-avoiding polygons and walks]{On self-avoiding polygons and walks: \\ counting, joining and closing}
\date{}

\author[A.~Hammond]{Alan Hammond}
\address{Departments of Mathematics and Statistics, 
  U.C. Berkeley,
  Berkeley, CA, 94720-3840, U.S.A.}
\email{alanmh@stat.berkeley.edu}

\begin{document}
\maketitle

\vspace{-0.9cm}

\begin{abstract}
For $d \geq 2$ and $n \in \N$, let $c_n = c_n(d)$ denote the number of length~$n$ self-avoiding walks beginning at the origin in the integer lattice $\Z^d$, and, for even~$n$, let $p_n = p_n(d)$ denote the number of length~$n$ self-avoiding polygons in $\Z^d$ up to translation. Then the probability   under the uniform law $\psaw_n$ on self-avoiding walks $\Gamma$ of any given odd length $n$ beginning at the origin that $\Gamma$ {\em closes} -- i.e., that $\Gamma$'s endpoint is a neighbour of the origin -- is given by $\psaw_n \big( \Gamma \closes \big) = 2(n+1) p_{n+1}/c_n$. The polygon and walk cardinalities share a common exponential growth: $\lim_n c_n^{1/n} = \lim_{n \in 2\N} p_n^{1/n} = \mu$ (where the common value $\mu \in (0,\infty)$ is called the connective constant). Madras~\cite{Madras95} has shown that $p_n \leq C n^{-1/2} \mu^n$ in dimension $d=2$, while the closing probability was recently shown in~\cite{ontheprob} to satisfy $\psaw_n \big( \Gamma \closes \big) \leq n^{-1/4 + o(1)}$ in any dimension $d \geq 2$. Here we establish that
\begin{itemize}
\item  $\psaw_n \big( \Gamma \closes \big) \leq n^{-1/2 + o(1)}$ for any $d \geq 2$;
\item   $\psaw_n \big( \Gamma \closes \big) \leq n^{-4/7 + o(1)}$ for a subsequence of odd $n$, if $d = 2$;
\item $p_n \leq n^{-3/2 + o(1)} \mu^n$ for a set of even $n$ of full density when $d=2$.
\end{itemize}
We also argue that the closing probability is bounded above by $n^{-(1-1/d) + o(1)}$ on a full density set when $d \geq 3$ for a certain variant of self-avoiding walk.  
\end{abstract}
\maketitle

\vspace{0.1cm}

\tableofcontents

\newpage

\part{Introduction}

Self-avoiding walk was introduced in the 1940s by Flory and Orr \cite{Flory,Orr47} as a model of a long polymer chain in a system of such chains at very low concentration. It is well known among the basic models of discrete statistical mechanics for posing problems that are simple to state but difficult to solve. Two recent surveys are the lecture notes~\cite{BDGS11}  and \cite[Section 3]{Lawler13}. 

\section{The model} 

We will denote by $\N$ the set of non-negative integers.
Let $d \geq 2$. For $u \in \R^d$, let $\dist u \dist$ denote the Euclidean norm of $u$. 
A {\em walk} of length $n  \in \N$ with $n > 0$ is a map $\gamma:\{0,\cdots,n \} \to \lattv$ 
such that $\dist \gamma(i + 1) - \gamma(i) \dist = 1$ for each $i \in \{0,\cdots,n-1\}$.
An injective walk is called {\em self-avoiding}. 
Let $\saw_n$ denote the set of self-avoiding walks~$\gamma$ of length~$n$ that start at $0$, i.e., with $\gamma(0) = 0$. 
We denote by $\psaw_n$ the uniform law on $\saw_n$. 
The walk under the law $\psaw_n$ will be denoted by $\Ga$. 

\medskip

A walk $\gamma \in \saw_n$ is said to close (and to be closing)
if $\dist \gamma(n) \dist = 1$. 
When the {\em missing edge} connecting $\gamma(n)$ and $\gamma(0)$
is added, a polygon results.
\begin{definition}
Let $\ga:\{0,\ldots,n-1\} \to \Z^d$ be a closing self-avoiding walk. For  $1 \leq i \leq n-1$, let $u_i$ denote the unordered nearest neighbour edge in $\Z^d$ with endpoints $\ga(i-1)$ and $\ga(i)$. Let $u_n$ denote $\ga$'s missing edge, with endpoints $\ga(n-1)$ and~$\ga(0)$. 
We call the collection of edges 
$\big\{ u_i: 1 \leq i \leq n \big\}$ the polygon of~$\ga$. A self-avoiding polygon in $\Z^d$ is defined to be any polygon of a closing self-avoiding walk in~$\Z^d$. The polygon's length is its cardinality. 
\end{definition}
We will usually omit the adjective self-avoiding in referring to walks and polygons. Recursive and algebraic structure has been used to analyse polygons in such domains as strips, as~\cite{BMB09} describes.

Note that the polygon of a closing walk has length that exceeds the walk's by one. Polygons have even length and closing walks, odd.


\section{Main results}

The {\em closing probability} is  $\psaw_n \big( \Gamma \closes \big)$. In~\cite{ontheprob}, an upper bound on this quantity of $n^{-1/4 + o(1)}$ was proved in general dimension.

In this paper, we use a variety of approaches to prove several closing probability upper bounds.

First, we revisit the closing probability upper bound method of~\cite{ontheprob}, styling it the {\em snake method via Gaussian pattern fluctuation}. This technique cannot hope to show that the closing probability decays faster than $n^{-1/2}$. We show that this decay can be proved in general dimension.
\begin{theorem}\label{t.closingprob}
Let $d \geq 2$. For any $\eps > 0$ and $n \in 2\N + 1$ sufficiently high,
$$
\psaw_n \big( \Gamma \closes \big) \leq n^{-1/2 + \eps} \, .
$$
\end{theorem}

Define the polygon number $p_n$ to be the number of length~$n$ polygons up to translation, and the walk number $c_n$ to be the number of length~$n$ walks beginning at the origin. As we shall soon review, the limiting exponential growth rates
$\lim_{n \in \N} c_n^{1/n}$ and 
$\lim_{n \in 2\N} p_n^{1/n}$ exist and coincide. Writing $\mu$ for the common value, we
define the real-valued {\em polygon number deficit} and {\em walk number excess} exponents $\thet_n$ and~$\xi_n$ according to the formulas
\begin{equation}\label{e.thetn}
 p_n = n^{-\thet_n} \cdot \mu^n \, \, \, \textrm{for $n \in 2\N$} 
\end{equation}
and
\begin{equation}\label{e.xin}
 c_n = n^{\xi_n} \cdot \mu^n \, \, \, \textrm{for $n \in \N$} \, .
\end{equation}

The closing probability may be written in terms of the polygon and walk numbers.
 There are $2n$ closing walks whose polygon is a given polygon of length~$n$, since there are $n$ choices of missing edge and two of orientation. Thus,
\begin{equation}\label{e.closepc}
\psaw_n \big( \Gamma \closes \big) = \frac{2(n+1) p_{n+1}}{c_n} \, ,
\end{equation}
for any $n \in \N$ (but non-trivially only for odd values of $n$).


 As we will shortly review,  it is a straightforward fact that $\liminf_{n \in 2\N} \thet_n$ and each $\xi_n$ are non-negative in any dimension~$d \geq 2$. Hara and Slade~\cite[Theorem 1.1]{HS92} used the lace expansion to prove that $c_n \sim C \mu^n$ when $d \geq 5$, while Madras and Slade \cite[Theorem 6.1.3]{MS93} have proved that $\thet_n \geq d/2 + 1$ in these dimensions for {\em spread-out} models, in which the vertices of $\Z^d$ are connected by edges below some bounded distance. Thus, the conclusion that the closing probability decays as fast as $n^{-d/2}$ has been reached for such models when $d \geq 5$. This conclusion may be expected to be sharp, but the opposing lower bound is not known to the best of the author's knowledge.
 
Madras~\cite{Madras91} has proved a bound on the moment generating function of the sequence $\big\{ p_n: n \in 2\N \big\}$ which when $d=3$ would assert~$\lim_{n \in 2\N} \thet_n \geq 1$ were this limit known to exist.
More relevantly for us, he has shown in~\cite{Madras95} using a polygon joining technique that $\theta_n \geq 1/2 - o(1)$ for $d=2$. We develop this technique to prove a stronger lower bound valid for typical high $n$.

\begin{definition}
The limit supremum density of 
a set $A$ of even, or odd, integers~is
$$
\limsup_n \frac{\big\vert A \cap [0,n] \big\vert}{\vert 2\N \cap [0,n] \vert} = \limsup_n \, n^{-1} \big\vert A \cap [0,2n] \big\vert \, .
$$
When the corresponding limit infimum density equals the limit supremum density, we naturally call it the limiting density. 
\end{definition}
\begin{theorem}\label{t.polydev}
 Let $d=2$. For any $\delta > 0$,
the  limiting density of the set of $n \in 2\N$ for which $\thet_n \geq 3/2 - \delta$ is equal to one. 
\end{theorem}


Applying~(\ref{e.closepc}) and $\xi_n \geq 0$ to Theorem~\ref{t.polydev}, we reobtain Theorem~\ref{t.closingprob} when $d=2$ on a subsequence of odd $n$ of full density.
Theorem~\ref{t.closingprob} is applicable in all dimensions, however.

As we rework the method of~\cite{ontheprob} to prove Theorem~\ref{t.closingprob}, we take the opportunity to present the technique in a general guise. The snake method is a proof-by-contradiction technique for deriving closing probability upper bounds that involves constructing sequences of laws of self-avoiding walks conditioned on increasingly severe avoidance events. We also exploit it in a new way  to prove closing probability upper bounds below $n^{-1/2}$ in two dimensions.
\begin{theorem}\label{t.thetexist}
Let $d = 2$. 
\begin{enumerate} 
\item
For any $\eps > 0$, the bound 
$$
\psaw_n \big( \Gamma \closes \big) \leq n^{-4/7 + \eps} 
$$
holds on a set of $n \in 2\N + 1$ of limit supremum density at least $1/{1250}$. 
\item
Suppose that the limits $\thet : = \lim_{n \in 2\N} \thet_n$
and $\xi : = \lim_{n \in \N} \xi_n$
exist in $[0,\infty]$. 
Then $\thet + \xi \geq 5/3$. Since 
\begin{equation}\label{e.closepc.formula}
\psaw_n \big( \Gamma \closes \big) = n^{-\theta - \xi + 1 + o(1)} 
\end{equation}
as $n \to \infty$ through odd values of $n$
by (\ref{e.closepc}), the closing probability is seen to be bounded above by $n^{-2/3 + o(1)}$. 
\end{enumerate} 
\end{theorem} 
(When $\thet + \xi = \infty$, (\ref{e.closepc.formula}) should be interpreted as asserting a superpolynomial decay in $n$ for the left-hand side.)

In our view, Theorem~\ref{t.thetexist}(1) is the most conceptually interesting result in this paper. Its proof brings together many of the ideas harnessed here, using the snake method and the polygon joining technique at once. Theorem~\ref{t.thetexist}(2) is only a conditional result, but it serves a valuable expository purpose: its proof is that of the theorem's first part with certain technicalities absent.

However far from rigorous they may be, the validity of the hypotheses of Theorem~\ref{t.thetexist}(2) are uncontroversial. For example, the limiting value $\thet$ is predicted to exist and to satisfy a relation with the Flory exponent~$\nu$ for mean-squared radius of gyration. The latter exponent is specified by the putative formula 
\begin{equation}\label{e.nu}
       \esaw_n \, \dist \Gamma(n) \dist^2   = n^{2\nu + o(1)}  \, ,
\end{equation}
where $\esaw_n$ denotes the expectation associated with~$\psaw_n$ (and where note that $\Gamma(n)$ is the non-origin endpoint of $\Gamma$); in essence, $\dist \Ga(n) \dist$  is supposed to be typically of order $n^{\nu}$. 
The hyperscaling relation that is expected to hold between $\thet$ and $\nu$ is 
\begin{equation}\label{e.hyperscaling}
    \thet   =   1 +  d \nu
\end{equation}
where the dimension $d \geq 2$ is arbitrary. In $d = 2$, $\nu = 3/4$ and thus $\theta = 5/2$ is expected. That $\nu = 3/4$ was predicted by the Coulomb gas formalism \cite{Nie82,Nie84} and then by conformal field theory \cite{Dup89,Dup90}. 
Hara and Slade~\cite{HS92,HS92b} used the lace expansion to show that $\nu = 1/2$ when $d \geq 5$ by demonstrating that, for some constant $D \in (0,\infty)$,  $\esaw_n \, \dist \Gamma_n \dist^2 - Dn$ is $O(n^{-1/4 + o(1)})$. This value of $\nu$ is anticipated in four dimensions as well, since $\esaw_n \, \dist \Gamma_n \dist^2$ is expected to grow as $n \big( \log n \big)^{1/4}$. 
In fact, the continuous-time weakly self-avoiding walk in $d=4$ has been the subject of an extensive recent investigation of Bauerschmidt, Brydges and Slade.  In \cite{BBSlog}, a $\log^{1/4}$ correction to the susceptibility is derived,
relying on 
 rigorous renormalization group analysis developed in a five-paper series~\cite{BSRigRen1,BSRigRen2,BSRigRen3,BSRigRen4,BSRigRen5}.



We mention that $\xi = 11/32$ is expected when $d=2$; in light of the $\thet = 5/2$ prediction and~(\ref{e.closepc}),
 $\psaw_n \big( \Ga \closes\big) = n^{-\psi + o(1)}$ with $\psi = 59/32$ is expected. The $11/32$ value was predicted by Nienhuis in \cite{Nie82} and can also be deduced from calculations concerning SLE$_{8/3}$: see \cite[Prediction~5]{LSW}. We will refer to $\psi$ as the {\em closing} exponent. Incidentally, the possibility that the half-integer value of~$\thet$ in~$d=2$ may indicate that polygons are more tractable than walks has been mooted in Tony Guttmann's survey~\cite{Guttmann12}: certainly the present article builds on the theme of~\cite{ontheprob} to provide further evidence that polygons, and in particular the closing probability, may be some of the more tractable aspects of the theory of self-avoiding walk.

Except for Theorem~\ref{t.closingprob}, we have chosen to keep an intent focus on the two-dimensional case in this article, in the hope that we maintain a focus on concepts rather than technicalities by doing so. It would be interesting however to try higher-dimensional versions of results such as Theorem~\ref{t.thetexist}.
 We do in fact present one further result, in two or more dimensions. We select a variant model tailored to eliminate technical difficulties.

\begin{definition}
The maximal edge local time of a nearest neighbour walk $\gamma:\{ 0,\cdots,n\} \to \Z^d$ is the maximum number of times that $\ga$ traverses an edge of $\Z^d$; more formally, it is the maximum cardinality of a subset $I \subseteq \{ 0,\cdots,n-1 \}$ such that the unordered sets  $\big\{ \ga(i),\ga(i+1) \big\}$ and $\big\{ \ga(j),\ga(j+1) \big\}$ coincide for each pair $(i,j) \in I^2$. Call $\ga$ $k$-edge self-avoiding if its maximal edge local time is at most~$k \in \N$.
Note that even $1$-edge self-avoiding walk satisfies a weaker avoidance constraint than does self-avoiding walk.

When considering (as we will) $3$-edge self-avoiding walks, we say that a walk $\ga$ as above closes if $\dist \ga(n) \dist = 0$.  Two such walks may be identified if they coincide after reparametrization by cyclic shift or reversal. A $3$-edge self-avoiding polygon is an equivalence class under this relation on closing walks. The length of such a polygon is the length of any of its members (and is $n$ if one of these members is $\ga$ as above).  Note that these definitions entail that
not only polygons but also
 closing walks have even length.

 Write $\carp_n$ and $\psawtd_n$ for the cardinality of, and uniform law on, the set $\sawtd_n$ of length~$n$ $3$-edge self-avoiding walks beginning at $0$.  
For $n \in 2\N$, let $\barp_n$ 
denote the number of
$3$-edge self-avoiding polygons of length~$n$ up to translation. 
\end{definition}
With these definitions,~(\ref{e.closepc}) has the counterpart
\begin{equation}\label{e.closepc.td}
\psawtd_n \big( \Gamma \closes \big) = \frac{2 n \barp_n}{\carp_n} 
\end{equation}
for any $n \in 2\N$. 

The connective constant $\lim_{n \in 2\N} \barp_n^{1/n}$ also exists for this model and we denote it by~$\muhat$. We define a real sequence $\big\{ \thethat_n : n \in \N \big\}$ so that
\begin{equation}\label{e.thetn.td}
 \barp_n = n^{-\thethat_n} \cdot \muhat^n \, \, \, \textrm{for $n \in 2\N$} \, .
\end{equation}
 
\begin{theorem}\label{t.threed}
Let $d \geq 2$. 
 For any $\delta > 0$,
the set of $n \in 2\N$ for which $\thethat_n \geq 2 - 1/d  - \delta$
has limiting  density equal to one. 
\end{theorem}
As we will later explain, we also have~$\carp_n \geq \muhat^n$. Also using~(\ref{e.closepc.td}), the next inference is immediate.
\begin{corollary}\label{c.threed}
Let $d \geq  2$. 
 For any $\delta > 0$,
the set of $n \in 2\N$ for which 
$$
\psawtd_n \big( \Gamma \closes \big) \leq n^{- (1 - 1/d) + \delta}
$$
has  limiting  density equal to one. 
\end{corollary}

Theorem~\ref{t.threed}'s proof is a general dimensional analogue of the polygon joining Theorem~\ref{t.polydev}.
Its main interest probably lies in the case when $d$ equals three or four. Known conclusions made by the lace expansion when $d \geq 5$ sometimes apply to spread-out models rather than the nearest neighbour one, but the method presumably offers a more plausible route to sharp conclusions in these dimensions.

\medskip

\noindent{\bf The structure of the paper.}
The paper consists of three further parts. 

The first of these, 
Part~\ref{c.boundspn}, uses the technique of polygon joining in order to prove the polygon number bound Theorems~\ref{t.polydev} and~\ref{t.threed}.

Part~\ref{c.thesnakemethod} presents the probabilistic snake method in a general guise, and uses it to prove the general dimensional closing probability upper bound Theorem~\ref{t.closingprob}.

The combinatorial and probabilistic ideas of these two parts are combined in the final Part~\ref{c.aboveonehalf}
when  Theorem~\ref{t.thetexist} is proved by using
the snake method alongside  the polygon joining technique.

\medskip

\noindent{\bf Acknowledgments.} 
I am very grateful to a referee for a thorough discussion of an earlier version of the article. Indeed, Theorem~\ref{t.polydev}'s present form is possible on the basis of a suggestion made by the referee, and this strengthened form has led to improvements in Theorem~\ref{t.thetexist}(1) and Theorem~\ref{t.threed}. 
I thank Hugo Duminil-Copin and Ioan Manolescu for many stimulating and valuable conversations about the central ideas in the paper. I thank Wenpin Tang for useful comments on a draft version. I would also like to thank Itai Benjamini for suggesting the problem of studying upper bounds on the closing probability.

\part{Bounds on polygon number via polygon joining}\label{c.boundspn}

This part is principally devoted to introducing Madras' polygon joining technique and using it to prove  Theorems~\ref{t.polydev} and~\ref{t.threed}.
The part has four sections. In the first, Section~\ref{s.prelude}, polygon joining and its most basic applications are described; a heuristic derivation is made of the lower bound in the hyperscaling relation~(\ref{e.hyperscaling}). This exposition is made because it provides a useful framework for discussing many of the paper's ideas. Section~\ref{s.three} provides a rigorous treatment of polygon joining that proves Theorem~\ref{t.polydev}. Certain variations to the joining technique will be employed in the final Part~\ref{c.aboveonehalf} to prove Theorem~\ref{t.thetexist}, and these changes are detailed in the next Section~\ref{s.pjprep}. In Section~\ref{s.threed}, the general dimensional Theorem~\ref{t.threed} is proved, by varying the argument for Theorem~\ref{t.polydev}.

\section{An heuristic prelude to polygon number bounds via joining}\label{s.prelude}

\subsection{Some general notation and tools}

We present many of our proofs for the case~$d=2$, partly in the hope that doing so will encourage the reader to visualise examples of the constructs being discussed. For this reason, some of the notation we now introduce is specifically adapted to the two dimensional case. 

\subsubsection{Multi-valued maps}\label{s.mvp}
For a finite set $B$, let $\mathcal{P}(B)$ denote its power set.
Let $A$ be another finite set. A {\em multi-valued map} from $A$ to $B$ is a function $\Psi: A \to \mathcal{P}(B)$. An arrow is a pair $(a,b) \in A \times B$ for which $b \in \Psi(a)$; such an arrow is said to be outgoing from $a$ and incoming to $b$. We consider multi-valued maps in order to find lower bounds on $\vert B \vert$, and for this, we need upper (and lower) bounds on the number of incoming (and outgoing) arrows. This next lemma is an example of such a lower bound.
\begin{lemma}\label{l.mvm}
Let $\Psi: A \to \mathcal{P}(B)$.
Set $m$ to be the minimum over $a \in A$ of the number of arrows outgoing from $a$, and $M$ to be the maximum over $b \in B$ of the number of arrows incoming to~$b$. Then
  $\vert B \vert \geq m M^{-1} \vert 
A \vert$.
\end{lemma}
\noindent{\bf Proof.} The quantities $M \vert B \vert$ and $m \vert A \vert$
are upper and lower bounds on the total number of arrows. \qed

\subsubsection{Denoting walk vertices and subpaths}
For $i,j \in \N$ with $i \leq j$, we write $[i,j]$ for  $\big\{ k \in \N:  i \leq k \leq j \big\}$. For a walk $\ga:[0,n] \to \Z^d$ and $j \in [0,n]$, we write $\ga_j$ in place of $\ga(j)$. For $0 \leq i \leq j \leq n$, $\ga_{[i,j]}$ denotes the subpath $\ga_{[i,j]}: [i,j] \to \Z^d$ given by restricting $\ga$. 

\subsubsection{Notation for certain corners of polygons}

\begin{definition}\label{d.corners}
The Cartesian unit vectors are denoted by $e_1$ and $e_2$ and the coordinates of $u \in \Z^2$ by $x(u)$ and $y(u)$.
For a finite set of vertices $V \subseteq \Z^2$, we define the northeast vertex $\northeast(V)$ in $V$ to be that element of $V$ of maximal $e_2$-coordinate; should there be several such elements, we take $\northeast(V)$ to be the one of maximal $e_1$-coordinate. That is, $\northeast(V)$ is the uppermost element of $V$, and the rightmost among such uppermost elements if there are more than one. Using the four compass directions, we may similarly define eight elements of $V$, 
including the lexicographically minimal and maximal elements of $V$, $\mathsf{WS}(V)$ and $\mathsf{EN}(V)$. We extend the notation to any self-avoiding walk or polygon $\gamma$, writing for example $\northeast(\gamma)$ for $\northeast(V)$, where $V$ is the vertex set of $\gamma$. 
For a polygon or walk $\ga$, set $\ymax(\ga) = y\big(\northeast(\gamma)\big)$, $\ymin(\ga) = y\big(\southeast(\gamma)\big)$, $\xmax(\ga) = x\big(\eastnorth(\gamma)\big)$ and $\xmin(\ga) = x\big(\westnorth(\gamma)\big)$.
The height~$h(\ga)$ of $\ga$ 
 is $\ymax(\gamma) - \ymin(\gamma)$  and its width~$w(\gamma)$ is  $\xmax(\gamma) - \xmin(\gamma)$.
\end{definition}

\subsubsection{Polygons with northeast vertex at the origin}

For $n \in 2\N$, let $\sap_n$ denote the set of length $n$ polygons $\phi$ such that $\northeast(\phi) = 0$. The set $\sap_n$ is in bijection with equivalence classes of length~$n$ polygons where polygons are identified if one is a translate of the other. Thus, $p_n =   \vert \sap_n \vert$. 

We write $\psap_n$ for the uniform law on $\sap_n$. 
A polygon  sampled with law $\psap_n$ will be denoted by $\Ga$, as a walk with law $\psaw_n$ is.

There are $2n$ ways of tracing the vertex set of a polygon $\phi$ of length~$n$: $n$ choices of starting point and two of orientation. We now select one of these ways. Abusing notation, we may write $\phi$ as a map from $[0,n]$ to $\Z^2$, setting $\phi_0 = \northeast(\phi)$, $\phi_1 = \northeast(\phi) - e_1$, and successively defining $\phi_j$ to be the previously unselected vertex for which $\phi_{j-1}$ and $\phi_j$ form the vertices incident to an edge in $\phi$, with the final choice $\phi_n = \northeast(\phi)$ being made. Note that $\phi_{n-1} = \northeast(\phi) - e_2$.

\subsubsection{Cardinality of a finite set $A$} This is denoted by either $\# A$ or $\vert A \vert$.


\subsubsection{Plaquettes}

The shortest non-empty polygons contain four edges. Certain such polygons play an important role in several arguments and we introduce notation for them now.

\begin{figure}
    \begin{center}
      \includegraphics[width=0.75\textwidth]{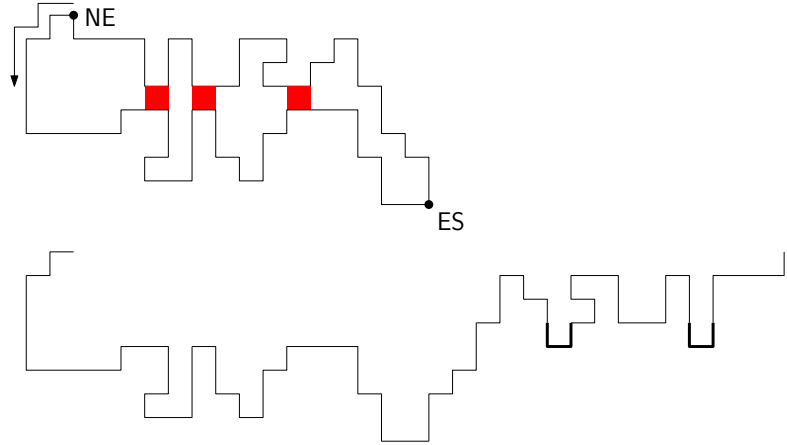}
    \end{center}
    \caption{In the upper sketch, a polygon $\phi$ with the beginning of the trajectory from $\northeast(\phi)$ indicated by an arrow. The polygon has three global join plaquettes, shaded red and labelled in the text $P^1$, $P^2$ and $P^3$ left to right (the order in which they are encountered counterclockwise along $\phi$ from $\northeast(\phi)$). In the lower picture, $\mathscr{S}_\kappa(\phi)$ with $\kappa = \big\{ 1,3 \big\}$ is depicted. The two bold subpaths indicate the reflections of the two diversions made by $\phi^{2,\kappa}$ that traverse $P^1$ and $P^3$ the long way around.}\label{f.reflect}
  \end{figure}

\begin{definition}
A plaquette is a polygon with four edges. Let $\phi$ be a polygon. A plaquette $P$ is called a join plaquette of~$\phi$ if $\phi$ and $P$ intersect at precisely the two horizontal edges of $P$. (The boundaries of the three shaded red squares of the polygon in the upcoming Figure~\ref{f.reflect} are join plaquettes of the polygon.)
Note that when $P$ is a join plaquette of $\phi$, the operation of removing the two horizontal edges in $P$ from $\phi$ and then adding in the two vertical edges in $P$ to~$\phi$ results in two disjoint polygons whose lengths sum to the length of $\phi$. We use symmetric difference notation and denote the output of this operation by $\phi \, \Delta \, P$.

The operation may also be applied in reverse: for two disjoint polygons $\phi^1$ and $\phi^2$, each of which contains one vertical edge of a plaquette~$P$, the outcome $\big( \phi^1 \cup \phi^2 \big) \, \Delta \, P$ of removing $P$'s vertical edges and adding in 
its horizontal ones is a polygon whose length is the sum of the lengths of $\phi^1$ and $\phi^2$.
\end{definition}

\subsection{Polygon superadditivity: polygon joining in its simplest guise}

Here are some very classical facts regarding the growth of walk and polygon numbers. 

\begin{proposition}\label{p.basic}
Consider any dimension $d \geq 2$.
\begin{enumerate}
\item For $n,m \in \N$, $c_{n+m} \leq c_n c_m$.
\item For $n,m \in 2\N$, $p_{n+m} \geq \tfrac{1}{d-1} p_n p_m$.
\item  The limits $\mu_1 : = \lim_{n \in \N} c_n^{1/n}$ and $\mu_2 : = \lim_{n \in \N} p_{2n}^{1/{2n}}$ exist; moreover, $c_n \geq \mu_1^n$ and $p_{2n} \leq (d-1)\mu_2^{2n}$.
\item The limits are equal: $\mu_1 = \mu_2$.
\end{enumerate}
\end{proposition} 
The common exponential growth constant $\mu_1 = \mu_2$ is called the connective constant and denoted by $\mu$. It is easily seen that $\mu \in \big[ d,2d - 1 \big]$. Duminil-Copin and Smirnov~\cite{DS10} have proved  using parafermionic observables that~$\mu$ equals~$\sqrt{2 + \sqrt{2}}$ for the honeycomb lattice. These observables have also been used in~\cite{BBDDG13} and~\cite{Gla13}, the second of which computes  the connective constant for a weighted walk model.

Polygon joining arguments are fundamental to the ideas in this paper. Perhaps the most basic such argument proves Proposition~\ref{p.basic}(2). 
We take the opportunity to present this joining idea in a simple guise.

\medskip

\noindent{\bf Proof of a portion of Proposition~\ref{p.basic}.} 
{\em (1)}. 
A walk $\gamma$ of length $n + m$ with $\gamma_0 = 0$ can be severed at $\gamma_n$ to form two subpaths, $\gamma_{[0,n]}$ and $\gamma_{[n,n+m]}$.
The first of these (called $\gamma^a$) is a length~$n$ walk beginning at $0$. The second becomes a length~$m$ walk (called $\gamma^b$) from~$0$ after translation by $-\ga_n$ (and reindexing of the domain from $[n,n+m]$ to $[0,m]$). The application $\gamma \to \big( \gamma^a, \gamma^b \big)$ is injective. Thus, $c_{n + m} \leq c_n c_m$.

\medskip 

\noindent{{\em (2)}.} 
Polygons cannot be severed, but they can be joined in pairs. We consider only the case of $d=2$ and equal polygon length $m =n \in 2\N$. The reader may consult for~\cite[Theorem 3.2.3]{MS93} for the general case. 
Consider a pair $\phi,\phi' \in \sap_n$. 
Relabel $\phi'$ by translating it so that $\westnorth(\phi')$ is one unit to the right of $\eastnorth(\phi)$. The plaquette $P$ whose upper left vertex is $\eastnorth(\phi)$ has one vertical edge in $\phi$ and one in $\phi'$. 
Note that $\chi : = \big( \phi \, \cup \, \phi' \big) \, \Delta \, P$ is a polygon of length~$2n$ which is either an element of~$\sap_{2n}$ or may be associated to such an element by making a translation. Moreover, the application $\big( \phi,\phi' \big) \to \chi$ is injective, because from~$\chi$ we can detect the location of the plaquette~$P$. The reader may wish to confirm this property, using that  $\phi$ and $\phi'$ have the same length. This injectivity implies that $p_{2n} \geq p_n^2$.

\medskip

\noindent{{\em (3)}.} 
The $\N$-indexed sequences with terms $\log c_n$ and $- \log \big( p_{2n}(d - 1)^{-1} \big)$ are subadditive by the two preceding parts.
Thus Fekete's lemma \cite[Lemma 1.2.1]{St97} implies the result. 

\medskip

\noindent{{\em (4)}.} 
That the expression (\ref{e.closepc}) is at most one implies that $\mu_1 \geq \mu_2$; the opposite inequality 
 follows directly from the following lower bound on the closing probability (specifically, from the resulting subexponential decay). \qed
\begin{lemma}\label{l.closeproblb}
There is a constant $\calnew > 0$ such that, for $n \in 2\N + 1$,
$$
 \psaw_n \big( \Gamma \closes \big) \geq \exp \big\{ - \calnew n^{1/2} \big\} \, .
$$
\end{lemma}
\noindent{\bf Proof.} The Hammersley-Welsh lower bound~\cite{hammersleywelsh} on the number of self-avoiding {\em bridges} in terms of walk number is a classical unfolding argument that is recounted in \cite[Chapter 3]{MS93}. It implies that there exists a constant $c_{HW} > 0$ such that, for all $n \in \N$, 
\begin{equation}\label{e.hw}
  c_n \leq e^{c_{HW} n^{1/2}} \mu^n \, ,
\end{equation} 
where here $\mu$ is specified to be $\mu_1 = \lim_n c_n^{1/n}$.
We may now use the bound $p_n \geq e^{-c n^{1/2}} \mu^n$ in \cite[Theorem 3]{kestenone} and~(\ref{e.closepc}) to conclude. Alternatively (but similarly), the lemma is proved from~(\ref{e.hw}) by further unfolding in~\cite{DKY11}: see Lemma~$5$ and the proof of Proposition~$3$ in that paper. \qed

\subsection{The polygon number deficit  exponent via a hyperscaling relation}

We review the standard heuristic derivation of the relation~(\ref{e.hyperscaling}); see equation~(1.4.14) of~\cite{MS93} for a representation of~(\ref{e.hyperscaling}) written in terms of the exponent $\alpha_{{\rm sing}} = 3 - \thet$, and Section~2.1 of this text for a slightly more detailed presentation of the derivation.

To describe the derivation, we hypothesise the existence of the limit
\begin{equation}\label{e.thetexist}
  \thet : = \lim_{n \in 2\N} \thet_n
\end{equation}
and suppose also that $\nu$ exists, given by ~(\ref{e.nu}).
Recall from~(\ref{e.xin}) the walk counterpart $\thetc_n$ to the polygon number deficit exponent, $\thet_n$.
Note the absence of a minus sign in~(\ref{e.xin}) compared to (\ref{e.thetn}). Note also that
Proposition~\ref{p.basic}(1) and~(2) imply that $\thetc_n \geq 0$ and~$\liminf \thet_n \geq 0$ for $n \in \N$.

Observe that the proof of Proposition~\ref{p.basic}(1) in fact shows that $\tfrac{c_{n+m}}{c_n c_m}$
equals the probability that independent samples of $\psaw_n$ and $\psaw_m$ {\em avoid} each other, i.e., intersect only at the origin. Supposing the existence of  the limit 
\begin{equation}\label{e.thetcexist}
\thetc := \lim_n \thetc_n \, ,
\end{equation}
we see that $\thetc$ would have the interpretation 
\begin{equation}\label{e.avoidexp}
 \Big( \psaw_n \times \psaw_n \Big) \Big( \Gamma^1 \, \textrm{avoids} \, \, \,  \Gamma^2 \Big) = n^{-\thetc + o(1)} \, .
\end{equation}
Assuming~(\ref{e.nu}), the endpoint $\Ga_n$ under $\psaw_n$ presumably adopts a location that is fairly uniform over a ball about the origin of radius of order $n^\nu$. For generic locations~$x$ in this ball, we infer that $\psaw_n \big( \Gamma_n = x \big) = n^{-d \nu + o(1)}$. For $x$ close to the origin, however, this probability is in addition penalized since $\big\{ \Gamma_n = x \big\}$  entails that the initial and final segments of $\Gamma$ come close without touching one another. When $\dist x \dist = 1$, the order of this additional penalty is given by~(\ref{e.avoidexp}). Thus, we expect
$$
  \psaw_n \big( \Gamma \closes \big) = n^{-d \nu - \thetc + o(1)} \, ,
$$ 
as $n \to \infty$ through odd values.
Using (\ref{e.closepc}) and our hypotheses (\ref{e.thetexist}) and (\ref{e.thetcexist}), we also find that
$$
  \psaw_n \big( \Gamma \closes \big) = n^{1 - \thet - \thetc   + o(1)} \, . 
$$
The $\thetc$ terms cancel and we obtain~(\ref{e.hyperscaling}).

\subsection{The hyperscaling relation lower bound argued by polygon joining}
We now present another heuristic derivation of the lower bound in~(\ref{e.hyperscaling}). The new argument takes longer to present than the first and is no more convincing. However, it is conceptually central to this paper, and it provides a useful framework for understanding some of the main ideas in 
 the proofs of our principal results. 

In fact, these concepts are illustrated by rederiving merely the lower bound  
\begin{equation}\label{e.thetonesided}
\thet \geq 1 +  2 \nu \, 
\end{equation}
when $d=2$.
We will derive~(\ref{e.thetonesided}) in three steps, arguing that $\thet \geq \nu$ and $\thet \geq 1 + \nu$ in the first and second steps. In light of the three step derivation, it is of interest to reflect on why we may also expect $\theta \leq 1 + 2\nu$; but we do not do so here. We also mention that essentially the same heuristics point to $\thet \geq 1 + d \nu$ in dimension~$d \geq 2$, with steps one and two reaching the conclusions that $\thet \geq (d-1) \nu$ and $\theta \geq 1 + (d-1)\nu$.

\medskip

\noindent{\bf Step one.} 
This step is an argument of Madras in~\cite{Madras95}.
When $d=2$, we will argue that 
\begin{equation}\label{e.madrasheur}
p_{2n} \geq n^{\nu} p_n^2 \, , 
\end{equation}
where recall that (\ref{e.nu}) specifies $\nu$. That $\theta \geq \nu$ follows directly (provided the two exponents exist).
The derivation of~(\ref{e.madrasheur}) develops the polygon joining argument in the proof of Proposition~\ref{p.basic}(2). The length $n$ polygons $\phi$ and $\phi'$ were joined in only one alignment, after displacement of $\phi'$ to a given location.
Madras argues under~(\ref{e.nu}) that there are at least an order of $n^\nu$ locations to which $\phi'$ may be translated and then attached to $\phi$. A total of $n^{\nu}$ distinct length $2n$ polygons results, and we obtain~(\ref{e.madrasheur}). 

 Where are these new locations for joining? Orient $\phi$ and $\phi'$ so that the height of each is at least its width; thus each height is at least of order $n^\nu$. Translate $\phi'$ vertically so that some vertices in $\phi$ and $\phi'$ share their $y$-coordinate, and then push  $\phi'$ to the right if need be so that this polygon is to the right of $\phi$. Then push $\phi'$ back to the left stopping just before the two polygons overlap. The two polygons contain vertices at distance one, and so it is plausible that in the locale of this vertex pair, we may typically find a plaquette $P$ whose left and right vertical edges are occupied by $\phi$ and $\phi'$. The operation in the proof of Proposition~\ref{p.basic}(2) applied with plaquette $P$ then yields a length $2n$ polygon. 
This construction began with a vertical translation of $\phi'$,
with different choices of this translation resulting in different outcomes for the joined polygon; since there is an order of $n^\nu$ different heights that may be used for this translation, we see that the bound~(\ref{e.madrasheur}) results.

There is in fact a technical difficulty in implementing this argument: in some configurations, no such plaquette~$P$ exists. Madras developed a local joining procedure that overcomes this difficulty in two dimensions. His procedure will be reviewed shortly, in Section~\ref{s.madrasjoin}. 

\medskip

\noindent{\bf Step two.} The next step in this second derivation of~(\ref{e.hyperscaling}) is to argue that
\begin{equation}\label{e.thetnew}
\theta \geq 1 + \nu \, .
\end{equation} 

We do so by arguing heuristically in favour of a strengthening of~(\ref{e.madrasheur}),
\begin{equation}\label{e.joinstronger}
 p_{2n} \geq n^{\nu - o(1)} \sum_{j = - n/2}^{n/2} p_{n -j} p_{n+j}  \, .
\end{equation}

Expressed using the polygon number deficit exponents, we would then have $n^{-\thet_{2n}} \geq n^\nu \sum_{j = - n/2}^{n/2} (n-j)^{-\thet_{n -j}} (n+j)^{-\thet_j}$. Using (\ref{e.thetexist}), the bound~(\ref{e.thetnew}) results.

To argue for (\ref{e.joinstronger}), note that, in deriving $p_{2n} \geq n^\nu p^2_n$, each polygon pair $(\phi,\phi') \in \sap_n \times \sap_n$ resulted in $n^\nu$ distinct length~$2n$ polygons. The length pair $(n,n)$ may be varied to be of the form $(n-j,n+j)$ for any $j \in [0,n/2]$. We are constructing a multi-valued map 
$$
\Psi: \bigcup_{\vert j \vert \leq n/2} \sap_{n-j} \times \sap_{n+j} \to \mathcal{P} \big( \sap_{2n} \big) \, 
$$
to the power set of $\sap_{2n}$ which associates to each polygon pair $(\phi,\phi')$ in the domain an order of $n^\nu$
elements of $\sap_{2n}$.  Were $\Psi$ injective, we would obtain~(\ref{e.joinstronger}). (The term {\em injective} is being misused: we mean that no two arrows of $\Psi$ are incoming to the same element of~$\sap_{2n}$.) The map is not injective but it is plausible that it only narrowly fails to be so: that is, abusing notation in a similar fashion, for typical $\chi \in \textrm{Range}(\phi)$, the cardinality of $\Psi^{-1}(\chi)$ is at most $n^{o(1)}$. A definition is convenient before we argue this.
\begin{definition}
Let $\phi$ be a polygon, and let $P$ be one of $\phi$'s join plaquettes. Let  $\phi^1$ and $\phi^2$ denote the disjoint polygons of which $\phi \, \Delta \, P$ is comprised.  If each of $\phi^1$ and $\phi^2$ has at least one quarter of the length of $\phi$, then we call $P$ a {\em macroscopic} join plaquette.
\end{definition}

To see that $\Psi$ is close to being injective, note that each pre-image of $\chi \in \textrm{Range}(\phi)$ under $\Psi$ corresponds to a macroscopic join plaquette of $\chi$.
Each macroscopic join plaquette entails a probabilistically costly macroscopic {\em four-arm} event, where four walks of length of order~$n$ must approach the plaquette without touching each other.
That $\chi$ belongs to $\textrm{Range}(\phi)$ amounts to saying that $\chi$ is an element of $\sap_{2n}$ having at least one macroscopic join plaquette. 
The four-arm costs 
make it plausible that a typical such polygon has only a few such plaquettes, gathered together in a small neighbourhood.
Thus, $\Psi$ is plausibly close to injective so that~(\ref{e.joinstronger}), and~(\ref{e.thetnew}), results.
The reader who would prefer less vagueness is directed to the upcoming Section~\ref{s.three} where a rigorous implementation of step two is offered.

\medskip

\noindent{\bf Step three.} A further term of $\nu$ is now sought to move from~(\ref{e.thetnew}) to~(\ref{e.thetonesided}).
Let $\sap_{2n}^{\rm mj}$ denote the subset of $\sap_{2n}$ whose elements have a macroscopic join plaquette.
Note that $\Psi$ is a map into the power set of $\sap_{2n}^{\rm mj}$, so that we have derived
\begin{equation}
 p^{{\rm mj}}_{2n} \geq n^{\nu - o(1)} \sum_{j = - n/2}^{n/2} p_{n -j} p_{n+j}  \, , 
\end{equation}
where $p^{{\rm mj}}_{2n} = \big\vert  \sap_{2n}^{\rm mj} \big\vert$. Since $\psap_{2n} \big( \Gamma \in \sap_{2n}^{\rm mj} \big) = p^{{\rm mj}}_{2n}/p_{2n}$, the following, alongside~(\ref{e.thetnew}), will imply (\ref{e.thetonesided}):
\begin{equation}\label{e.mjpbound}
\psap_{2n} \big( \Gamma \in \sap_{2n}^{\rm mj} \big) \leq C n^{-\nu} \, .
\end{equation}
To argue for this, consider a typical polygon $\phi \in \sap_{2n}$. It occupies horizontal coordinates including $0$ in an interval of length of order $n^{\nu}$. Suppose that this interval contains $\big[ -5 n^\nu/4 ,  5n^\nu/4 \big]$ (where $5/4$ may be any constant exceeding one). The polygon crosses the strip $\big[ -n^\nu, n^\nu \big]$ at least twice. Let $\psi^1$ and $\psi^2$ denote the subpaths of $\phi$ consisting of the uppermost and lowermost of these crossings. Suppose  that the origin, which is $\northeast(\phi)$, lies in $\psi^1$. 
(This event plausibly has a probability that is uniformly positive in $n$. In the ensuing argument, it is certainly valid to question what happens in the opposite case. We will comment briefly on questions of this type having presented our heuristics.)

We now  conduct a {\em resampling experiment} that realizes the law $\psap_{2n}$ in a way that is useful for deriving~(\ref{e.mjpbound}). 
In the first instance, we simply sample the law $\psap_{2n}$.
Call the resulting polygon the {\em original sample}.
We record two pieces of information concerning this sample: 
\begin{itemize}
\item the subpath $\psi^1$; and the subpath $\psi^2$, up to vertical translation.
\end{itemize}  
Call this information the {\em retained data}.

  \begin{figure}
    \begin{center}
      \includegraphics[width=0.9\textwidth]{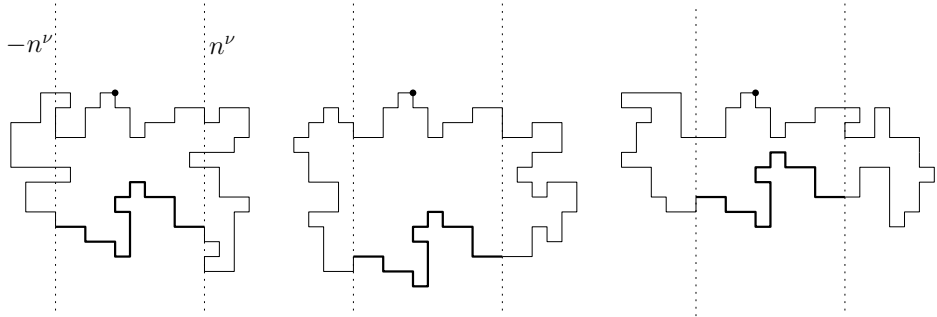}
    \end{center}
    \caption{An instance of the original sample $\phi$ of $\psap_{2n}$ on the left, with $\psi^2$ drawn in bold and the vertex $\northeast(\phi) = 0$ marked with a dot. Two reconstructed samples of $\psap_{2n}$ appear alongside. In the right sketch, the $\psi^2$-translate appears in its highest possible location and creates a macroscopic join plaquette.
}\label{f.vertshift}
  \end{figure}

We then form a new polygon, the {\em reconstructed sample}, whose law is given by the conditional distribution of the law $\psap_{2n}$ subject to the constraint that the retained data associated to this reconstructed sample is equal to that 
given in the bullet-point for the original sample.
Note that since the original sample has law $\psap_{2n}$, so does the reconstructed sample.

 We argue for the bound~(\ref{e.mjpbound}) by deriving it for the reconstructed sample copy of the law~$\psap_{2n}$. In a given but typical instance of retained data, what happens to the subpaths $\psi^1$ and $\psi^2$ during the formation of the reconstructed sample? The path $\psi^1$ stays in place, while $\psi^2$ undergoes a random vertical shift, as Figure~\ref{f.vertshift} illustrates. It is consistent with the hypothesis~(\ref{e.nu}) that the height of typical polygons under $\sap_{2n}$ is of order~$n^\nu$. It is also a natural hypothesis, akin to the Russo-Seymour-Welsh theorem for critical percolation (which has been recently recounted in~\cite{Werner09}), that any macroscopic subpath of such a polygon has comparable probability of moving in one or other of the four compass directions by an amount of the order of the polygon's diameter. These suppositions point to the conclusion that the random vertical shift experienced by $\psi^2$ during the reconstruction step has order of magnitude $n^{\nu}$ and that any admissible shift in this range has probability of order $n^{-\nu}$. The highest such shift is the only one that may bring the translate of $\psi^2$ within distance one of the upper subpath~$\psi^1$. This translate may result in a macroscopic join plaquette with $x$-coordinate in $\big[ - n^\nu, n^\nu \big]$, but all the others may not. Thus, there is probability at most a constant multiple of $n^{-\nu}$ that the reconstructed sample polygon has a macroscopic join plaquette.  This rests the case for~(\ref{e.mjpbound}) and thus for~(\ref{e.thetonesided}).

(One might object to this derivation on the grounds that it is only with positive probability that the reconstruction mechanism has any effect. For example, a positive probability set of original samples stay in the strip $[-n^\nu,n^\nu] \times \R$. When this happens, the original and resampled copies are equal, and the argument tells us nothing. 
A non-rigorous solution to this objection is offered by shortening the width of the strip by a factor of that is a large power of $\log n$. It is then plausible that all but a small set of original sample polygons, whose $\psap_{2n}$-probability has super-polynomial decay in $n$, cross at least one such strip. The vertical mobility between the original and resampled polygons may then be exploited as described above. A union bound is then needed involving a sum over a polylogarithmic number of such strips. This leads to~(\ref{e.mjpbound}) with an exponent of $\nu - o(1)$ on the right-hand side.)

\subsection{The polygon number deficit exponent as the sum of three terms}\label{s.sumthree}

We now formulate some notation that will permit us to describe the upcoming proofs in terms of the preceding three step derivation. We informally specify three exponents $\eone$, $\etwo$ and $\ethree$ (and suppose that they exist).
\begin{itemize}
\item
 Let the number of ways of joining a typical pair of length $n$ polygons in step one scale as $n^{\eone}$. 
 \item
 Let the number of lengths $k \in [n/2,n]$ for which typically polygon pairs of lengths $k$ and $2n-k$ may be joined in about $n^{\eone}$ ways scale as $n^{\etwo}$. 
 \item Finally, suppose that
$$
\psap_{2n} \big( \Gamma \in \sap_{2n}^{\rm mj} \big) \,\, \textrm{scales as}  \, \, n^{-\ethree} \, .
$$
\end{itemize}
Then we have seen that we should expect that $\thet \geq \eone + \etwo + \ethree$ (and in fact equality should hold), and that
 $\eone = \nu$, $\etwo = 1$ and $\ethree = \nu$ in $d=2$. Madras has proved $\theta \geq 1/2$ in this dimension by arguing that $\eone$ is indeed  at least $\nu$ and using the trivial bound $\nu \geq 1/2$; we may say that he carried out a $(1/2,0,0)$-argument, in that these are the proven lower bounds on the respective exponents.


Theorems~\ref{t.polydev},~\ref{t.thetexist} and~\ref{t.threed}
will be proved with the aid of $(\eone,\etwo,\ethree)$-arguments.
This language carries the caveat that it is merely a mnemonic, and we make no attempt to specify the three exponents precisely.

\section{Polygon number bounds via joining: a rigorous treatment}\label{s.three}

In this section, we prove Theorem~\ref{t.polydev}, (so that $d=2$ is assumed throughout the section).
The derivation of the theorem involves a $(1/2,1,0)$-argument. Indeed, the next proposition is a key step. It offers a rigorous interpretation of~(\ref{e.joinstronger}).


\begin{definition}\label{d.highpolynum}
For $\maceta > 0$, the set  $\highpolynum_\maceta \subseteq 2\N$ of $\maceta$-{\em high polygon number} indices is given by
$$
 \highpolynum_\maceta = \Big\{ n \in 2\N : p_n \geq n^{- \maceta} \mu^n \Big\} \, .
$$
\end{definition}



\begin{proposition}\label{p.polyjoin.newwork}
For any $\maceta > 0$, there is a constant $\Cpoly = \Cpoly(\maceta) > 0$ such that, for $n \in 2\N \cap \highpolynum_{\Cpcp}$,
\begin{equation}\label{e.polyjoin.newwork}
 p_n \geq \frac{1}{\Cpoly \log n} \sum_{j \in 2\N \, \cap \, [2^{i-1},2^i]} (n-j)^{1/2} p_j \, p_{n-j}  \, ,
\end{equation}
where $i \in \N$ is chosen so that $n \in 2\N \cap [2^i,2^{i+1}]$. 
\end{proposition}

In essence, the meaning of the proposition is that when  $n \in 2\N \cap \highpolynum_{\Cpcp}$, $p_n$ is at least a small constant multiple of $n^{1/2} (\log n)^{-1} \sum p_j p_{n-j}$, where the sum is over an interval of order $n$ indices $j$ around the value $n/2$ (though such a statement does not follow directly from the proposition).

are both positive and even. As such, the number of contributing summands reaches capacity only for choices of $n$ far from the two endpoints of the given dyadic scale.

Perhaps our proof of  Theorem~\ref{t.polydev} via Proposition~\ref{p.polyjoin.newwork}
can be refined to quantify the rate of convergence of the density one index set in the theorem.
However, the proposition in isolation is inadequate for proving the conclusion~$\thet_n \geq 3/2 - o(1)$ for all $n \in 2\N$, because this tool permits occasional spikes in the value of the $p_n$, as the sequence
$$
 \thet_n = \left\{
  \begin{array}{l l}
    1 & \quad \text{if $n$ is a power of $2$,}\\
    3/2 + \tfrac{1}{100} & \quad \text{if $n \in 2\N$ is otherwise,}
  \end{array} \right.
$$
demonstrates.

This section has five subsections. The first four present tools needed for the proof of Proposition~\ref{p.polyjoin.newwork}, with the proof of this result in the fourth. The fifth proves Theorem~\ref{t.polydev} as a consequence of the proposition.

Proposition~\ref{p.polyjoin.newwork} will be proved by
rigorously implementing a version of step two of the hyperscaling relation lower bound derivation in the preceding section. 
Of course, this used a polygon joining technique. 
In the first subsection, we specify the local details of the technique due to Madras that we will use. Recall that the heuristic argument in favour of~(\ref{e.joinstronger})
depended on the near-injectivity of the multi-valued map~$\Psi$, which was argued by making a case for the sparsity of macroscopic join plaquettes.
In the second subsection, we present
Proposition~\ref{p.globaljoin} and Corollary~\ref{c.globaljoin},
our rigorous versions of this sparsity claim. In the third subsection, we set up the apparatus needed to specify our joining mechanism~$\Psi$, and in the fourth, we define and analyse the mechanism (and so obtain Proposition~\ref{p.polyjoin.newwork}).

\subsection{Madras' polygon joining procedure}\label{s.madrasjoin}

When a pair of polygons is close, there may not be a plaquette whose vertical edges are divided between the two elements of the pair.  We now recall Madras' joining technique which works in a general way for such pairs. 

\medskip

Consider two polygons $\tau$ and $\sigma$ of lengths $n$ and $m$ for which the intervals 
\begin{equation}\label{e.intint}
\textrm{$\big[ \ymin(\tau) - 1 , \ymax(\tau) + 1 \big]$ and  $\big[ \ymin(\sigma) - 1 ,\ymax(\sigma) + 1 \big]$ intersect} \, .
\end{equation}
 Madras' procedure joins  $\tau$ and $\sigma$ to form a new polygon of length $n + m + 16$ in the following manner. 

First translate $\sigma$ to the right by far enough that the $x$-coordinates of the vertices of this translate are all strictly greater than all of those of $\tau$. Now shift $\sigma$ to the left step by step until the first time at which there is a pair of vertices, one in $\tau$ and the other in the $\sigma$-translate, that share an $x$-coordinate and whose $y$-coordinates differ by at most two; such a moment necessarily occurs, by the assumption~(\ref{e.intint}). Write $\sigma' = \sigma + T_1 e_1$ (with $T_1 \in \Z$) for this particular horizontal translate of $\sigma$. There is at least one vertex $z \in \Z^2$ such that the set $\big\{ z - e_2,z, z + e_2 \big\}$ contains a vertex of $\tau$ and a vertex of $\sigma'$. The set of such vertices contains at most one vertex with any given $y$-coordinate. Denote by $Y$ the vertex $z$ with the maximal $y$-coordinate. 

\begin{figure}
    \begin{center}
      \includegraphics[width=0.6\textwidth]{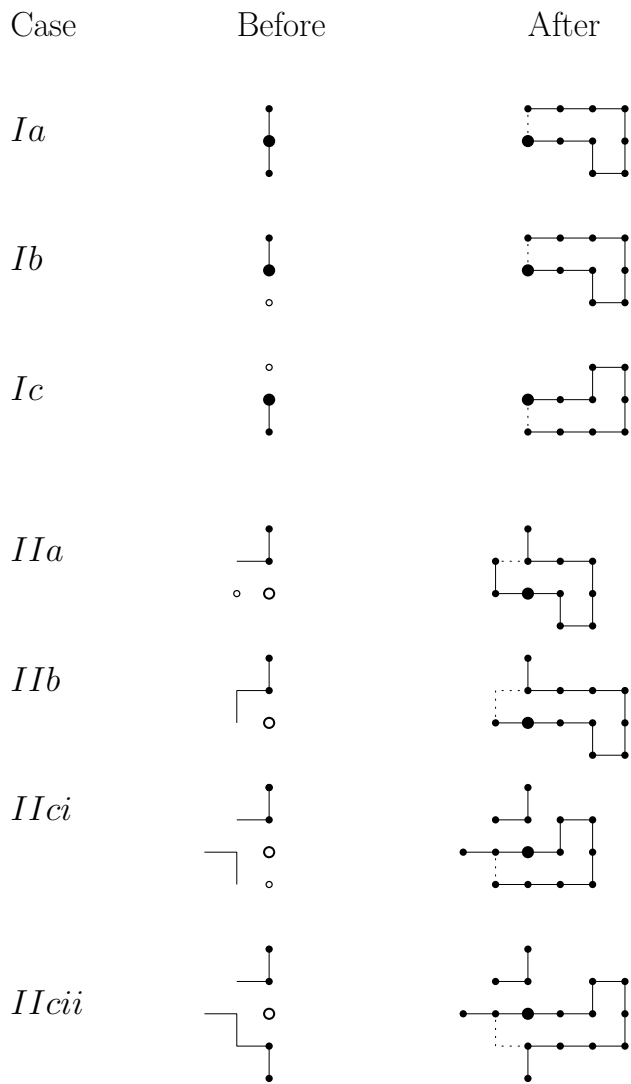}
    \end{center}
    \caption{Changes made near the vertex $Y$ in a polygon $\tau$ to produce $\modify{\tau}$ are depicted. The second column depicts $\tau$ around~$Y$, with $Y$ indicated by a large disk or circle; a disk denotes a vertex that belongs to $\tau$ and a circle one that does not; the line segments denote edges of $\tau$. In the third column, the modified polygon $\modify{\tau}$ is shown in the locale of $Y$. The vertex $Y$ is shown as a black disk. Black line segments are edges in $\tau$ or $\modify{\tau}$ and dotted line segments on the right are edges in $\tau$ that are removed in the formation of~$\modify{\tau}$. Several cases are not depicted. These may be labelled cases IIIa, IIIb, IIIci and IIIcii. In each case, the picture of $\tau$ and $\modify{\tau}$ is formed by reflecting the counterpart case II picture horizontally through~$Y$.}\label{f.madrascases}
  \end{figure}

Madras now defines a modified polygon $\modify{\tau}$, which is formed from $\tau$ by changing its structure in a neighbourhood of $Y \in \Z^2$. Depending on the structure of $\tau$ near~$Y$, either two edges are removed and ten edges added to form $\modify{\tau}$ from $\tau$, or one edge is removed and nine are added. As such, $\modify{\tau}$ has length $n + 8$. The rule that specifies $\modify{\tau}$ is recalled from~\cite{Madras95} in Figure~\ref{f.madrascases}. 

A modified polygon formed from $\sigma'$ is also defined. Rotate $\sigma'$ about the vertex~$Y$ by $\pi$ radians to form a new polygon $\sigma''$. Form $\modify{\sigma''}$ according to the same rules, recalled in Figure~\ref{f.madrascases}. Then rotate back the outcome by $\pi$ radians about $Y$ to produce the modification of 
$\sigma'$, which to simplify notation we denote by $\modify{\sigma}$. 

Writing $Y = (Y_1,Y_2)$ in place of $\big( x(Y) , y(Y) \big)$, note that no vertex of $\tau$ belongs to the {\em right corridor}
$\big\{ Y_1 + 1 , Y_1 + 2, \cdots \big\} \times  \big\{ Y_2 - 1 , Y_2 , Y_2 + 1 \big\}$, the region
that lies strictly to the right of $\big\{ Y-e_2,Y,Y+e_2 \big\}$. (Indeed, it is this fact that implies that $\modify{\tau}$ is a polygon.)
Equally, no vertex of $\sigma'$ belongs to the {\em left corridor}  
$\big\{ \cdots , Y_1 - 2 , Y_1 -1 \big\} \times \big\{ Y_2 - 1 , Y_2 , Y_2 + 1 \big\}$.

Note that the polygon $\modify{\tau}$ extends $\tau$ to the right of $Y$ by either two or three units inside the right corridor (by two in case IIa, IIci or IIIci and by three otherwise). Likewise $\modify{\sigma}$ extends $\sigma'$ to the left of $Y$ by either two or three units in the left corridor (by two when $\sigma''$ satisfies case IIa, IIci or IIIci and by three otherwise).

Note from Figure~\ref{f.madrascases} that, in each case, $\modify{\tau}$ contains two vertical edges that cross the right corridor at the maximal $x$-coordinate adopted by vertices in $\modify{\tau}$ that lie in this corridor. Likewise, $\modify{\sigma}$ contains two vertical edges that cross the left corridor at the minimal $x$-coordinate adopted by vertices in $\modify{\sigma}$ that lie in the left corridor.

Translate $\modify{\sigma}$ to the right by $T_2$ units, where $T_2$ equals
\begin{itemize}
\item five when one of cases IIa, IIci and IIIci obtains for both $\tau$ and $\sigma''$;
\item six  when one of these cases obtains for exactly one of these polygons;
\item seven when none of these cases holds for either polygon. 
\end{itemize}
Note that $\modify\tau$ and $\modify\sigma + T_2 e_1$ are disjoint polygons such that, for some pair of vertically adjacent plaquettes $(P^1,P^2)$ in the right corridor (whose left sides have $x$-coordinate either $Y_1 + 2$ or $Y_1 + 3$), the edge-set of $\modify\tau$ intersects the plaquette pair on the two left sides of $P^1$ and $P^2$, while the edge-set of $\modify\sigma + T_2 e_1$ intersects this pair on the two right sides of $P^1$ and $P^2$. Let $P^1$ denote the upper element of this plaquette pair.

The polygon that Madras specifies as the join of $\tau$ and $\sigma$ is given by $\big( \modify{\tau} \cup (\modify\sigma + T_2 e_1) \big) \, \Delta \, P^1$. 
Note that, to form the join polygon, $\sigma$ is first horizontally translated by $T_1$ units to form $\sigma'$, modified locally to form $\modify\sigma$, and then further horizontally translated by $T_2$ units to produce the polygon $\modify\sigma + T_2 e_1$ that is joined onto $\modify\tau$. Thus, $\sigma$ undergoes a horizontal shift by $T_1 + T_2$ units as well as a local modification before being joined with~$\modify\tau$.

  \begin{figure}
    \begin{center}
      \includegraphics[width=0.75\textwidth]{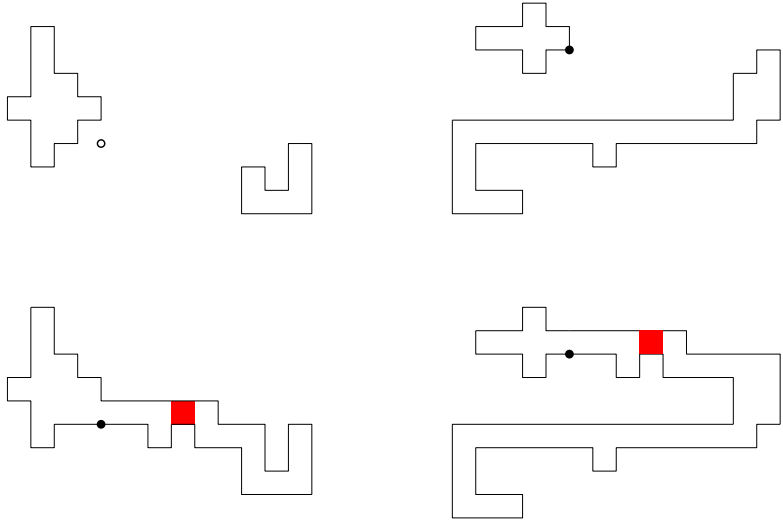}
    \end{center}
    \caption{Two pairs of Madras joinable polygons and the Madras join polygon of each pair. The vertex $Y$ is indicated by a circle or a disk in each picture. The junction plaquette is shaded red. On the left, $\tau$ satisfies case IIb and $\sigma''$, case IIa; on the right, Ib and~IIb.}\label{f.madrasjoinable}
  \end{figure}

\begin{definition}
For two polygons $\tau \in \sap_n$ and $\sigma \in \sap_m$ satisfying~(\ref{e.intint}), define the {\em Madras join polygon} 
$$
J(\tau,\sigma) = 
\Big( \modify{\tau} \cup \big( \modify\sigma + T_2 \, e_1 \big) \Big) \, \Delta \, P^1 \, \in \,  \sap_{n + m + 16} \, .
$$

The plaquette $P^1$ will be called the {\em junction} plaquette.

Such polygons $\tau$ and $\sigma$ are called {\em Madras joinable}  
if $T_1 + T_2 = 0$: that is, no horizontal shift is needed so that $\sigma$ may be joined to $\tau$ by the above procedure. Note that the modification made is local in this case: $J(\tau,\sigma) \Delta \big( \tau \cup \sigma \big)$ contains at most twenty edges. See Figure~\ref{f.madrasjoinable}.
\end{definition}

\subsection{Global join plaquettes are few}

Recall that in step two of the derivation of~(\ref{e.thetonesided}), the near injectivity of the multi-valued map $\Psi$ was argued as a consequence of the sparsity of macroscopic join plaquettes. We now present in~Corollary~\ref{c.globaljoin} a rigorous counterpart to this sparsity assertion. In the rigorous approach, we use a slightly different definition to the notion of macroscopic join plaquette. 


\begin{definition}\label{d.globaljoin}
For $n \in 2\N$, let $\phi \in \sap_n$. 
A join plaquette $P$ of $\phi$ is called {\em global} if the two polygons 
comprising $\phi \, \Delta \, P$ may be labelled $\phi^\ell$ and $\phi^r$ in such a way that 
\begin{itemize}
\item every rightmost vertex in $\phi$ is a vertex of $\phi^r$;
\item and  $\northeast(\phi)$ is a vertex of $\phi^\ell$.
\end{itemize}
Write $\globaljoin_\phi$ for the set of global join plaquettes of the 
polygon~$\phi$.
\end{definition}

\begin{proposition}\label{p.globaljoin}
There exists $c > 0$ such that, for $n \in 2\N$ and any $k \in \N$,  
$$
  \# \Big\{ \phi \in \sap_{n} : \big\vert \globaljoin_\phi \big\vert \geq k \Big\} \leq  c^{-1} 2^{-k \mu^{-2}/2}    \mu^n   \, .
$$
\end{proposition}
\begin{corollary}\label{c.globaljoin}
For all $\maceta > 0$ and $\Cgjone > 0$, there exist $\Cgjtwo, \Cgjthree > 0$ such that, for $n \in \highpolynum_\maceta$,
$$
\psap_{n} \Big( \, \big\vert \globaljoin_\Gamma \big\vert \geq \Cgjtwo \log n \Big)
 \leq \Cgjthree n^{- \Cgjone} \, .
$$
\end{corollary}

The next lemma will be used in the proof of Proposition~\ref{p.globaljoin}.
\begin{lemma}\label{l.globaloneedge}
Let $n \in 2\N$ and $\phi \in \sap_n$.  Writing $j \in [0,n]$ so that $\eastsouth(\phi) = \phi_j$, consider the two subpaths
$\phi_{[0,j]}$ and $\phi_{[j,n]}$, the first starting at $\northeast(\phi) = 0$ and the second ending there. Each of these paths contains precisely one of the two horizontal edges of any element in $\globaljoin_\phi$.
\end{lemma} 
\noindent{\bf Proof.} For given $\phi \in \sap_n$, let $P \in \globaljoin_\phi$. We may decompose $\phi \, \Delta \, P$ as $\phi^\ell \cup \phi^r$ in accordance with Definition~\ref{d.globaljoin}. We then have that $\northeast(\phi)$ is a vertex of $\phi^{\ell}$ and $\eastsouth(\phi)$ a vertex of $\phi^{r}$. The path $[0,n] \to \Z: j \to \phi_j$ leaves $0$ to follow  $\phi^{\ell}$ until passing through an edge in $P$ to arrive in~$\phi^{r}$, tracing this polygon until passing back through the other horizontal edge of $P$ and following $\phi^{\ell}$ until returning to $\northeast(\phi)$. It is during the trajectory in $\phi^{r}$ that the visit to $\eastsouth(\phi)$ is made. \qed    

\medskip

\noindent{\bf Proof of Proposition~\ref{p.globaljoin}.}
For $n \in \N$, an element $\gamma \in \saw_n$
is called a {\em half-space} walk if $y(\gamma_i) \leq 0$ for each $i \in [0,n]$. We call a half-space walk $\ga$ {\em returning} if,
after the last visit that $\ga$ makes to the lowest $y$-coordinate that this walk attains, $\ga$ makes a unique visit to the $x$-axis, with this occurring at its endpoint $\ga_n$.
Let $\rhssaw_n \subset \saw_n$ denote the set of length~$n$ returning half-space walks. We will first argue that, for any $n \in \N$, 
\begin{equation}\label{e.rhswalkbound}
\vert \rhssaw_n \vert \leq \mu^n \, .   
\end{equation}
To see this, we consider a map $\rhssaw_n \to \saw_n$. Let $\ga \in \rhssaw_n$, and let $j \in [0,n]$
be the index of the final visit of $\ga$ to the lowest $y$-coordinate visited by $\ga$. The image of $\ga$ is defined to be the concatenation of $\ga_{[0,j]}$ and the reflection of $\ga_{[j,n]}$ through the horizontal line that contains $\ga_j$.
(We have not defined concatenation but hope that the meaning is evident.) Our map is  injective because, given an element in its image, the horizontal coordinate of the line used to construct the image walk may be read off from that walk; (the coordinate is one-half of the $y$-coordinate of the image walk's non-origin endpoint). Its image lies in the set of {\em bridges} of length~$n$, where a bridge is a walk whose starting point has maximal $y$-coordinate and whose endpoint uniquely attains the minimal $y$-coordinate. 
The set of bridges of length~$n$ beginning at the origin 
has cardinality at most $\mu^n$: this classical fact, which follows from (1.2.17) in~\cite{MS93} by symmetries of the Euclidean lattice, is proved by a superadditivity argument with similarities to the proof of Proposition~\ref{p.basic}(2). Thus, considering this map proves~(\ref{e.rhswalkbound}).

Noting these things allow us to reduce the proof of the proposition to verifying the following assertion. 
There exists $c > 0$ such that, for $\delta \in (0,1/2)$, $n \in 2\N$ and any $k \in \N$,  
\begin{equation}\label{e.globaljoin.reduce}
    \#  \, \rhssaw_{n + 2\lfloor \delta k \rfloor}  \, \geq \, c \, \delta^{-\delta k} \cdot \# \Big\{ \phi \in \sap_{n} : \big\vert \globaljoin_\phi \big\vert \geq k \Big\}  \, .
\end{equation}
Indeed, applying (\ref{e.rhswalkbound}) with the role of $n$ played by $n + 2\lfloor \delta k \rfloor$, we see from~(\ref{e.globaljoin.reduce}) that
$$
\# \Big\{ \phi \in \sap_{n} : \big\vert \globaljoin_\phi \big\vert \geq k \Big\} \leq c^{-1} \big( \mu^2 \delta \big)^{\delta k} \mu^n \, .
$$
Setting $\delta = \mu^{-2}/2$, we obtain Proposition~\ref{p.globaljoin}.

To complete the proof of the proposition, we must prove~(\ref{e.globaljoin.reduce}), and this we now do.
Let $\phi \in \sap_n$. Setting $j \in [0,n]$ so that 
$\eastsouth(\phi) = \phi_j$ (as we did in Lemma~\ref{l.globaloneedge}),
write $\phi^1 =  \phi_{[0,j]}$ and $\phi^2 = \phi_{[j,n]}$. 
Writing $\reflect^2_z$ for reflection in the vertical ($e_2$-directed) line that passes through $z \in \Z^2$,
define a map $\mathscr{S}: \sap_n \to \rhssaw_n$ to be the concatenation
$$
\mathscr{S}(\phi) = \phi^1 \circ \reflect^2_{\eastsouth(\phi)}(\phi^2) \, .
$$
By Lemma~\ref{l.globaloneedge}, each of $\phi^1$ and $\phi^2$ traverses precisely one horizontal edge of each of $\phi$'s global join plaquettes. Set $r = \# \globaljoin_\phi$ and enumerate $\globaljoin_\phi$ by the sequence $\big( P^1,\cdots,P^r \big)$ (in an arbitrary order; for example, in the order in which $\phi^1$ traverses an edge of each plaquette). 
For each $j \in \{1,\cdots,r \}$, let $(s_j,f_j)$ denote the unique edge in $P^j$ traversed by $\phi^2$. 
Consider the path formed by modifying $\phi^2$ so that the one-step subpath $(s_j,f_j)$ is replaced by a three-step subpath from $s_j$ to $f_j$ that traverses the plaquette $P^j$ using its three edges other than $(s_j,f_j)$. The modification may be made iteratively for several choices of $j \in \{ 1,\cdots,r \}$, and the outcome is independent of the order in which the modifications are made. In this way, we may define a modified path $\phi^{2,\kappa}$ for each $\kappa \subseteq \{1,\cdots,r\}$, under which the modified route is taken along plaquettes $P^j$ precisely when $j \in \kappa$. Note that $\phi^{2,\kappa}$ is a self-avoiding walk whose length exceeds $\phi^2$'s by $2\vert \kappa \vert$: it is self-avoiding because this walk differs from $\phi^2$
by several disjoint replacements of one-step subpaths by three-step alternatives, and, in each case, the two new vertices visited in the alternative route are vertices in $\phi^1$, and, as such, cannot be vertices in $\phi^2$.

Note further that the intersection of the edge-sets of $\phi^1$ and $\phi^{2,\kappa}$ equals $\cup_{j \in \kappa} E(\phi^1) \cap E(P^j)$ (where the sets in the union are each singletons).

For each $\kappa \subseteq \{ 1,\cdots,r\}$, define $\mathscr{S}_\kappa(\phi) \in \rhssaw_{n + 2 \vert \kappa \vert}$,
$$
\mathscr{S}_\kappa(\phi) = \phi^1 \circ \reflect^2_{\eastsouth(\phi)}(\phi^{2,\kappa}) \, .
$$

Recall that $\delta \in (0,1/2)$ and that $k \in \N$ is given. Consider the multi-valued map 
$$
\Psi = \Psi_k : \Big\{ \phi \in \sap_n: \vert \globaljoin_\phi \vert \geq k \Big\}
 \to \mathcal{P} \big( \rhssaw_{n + 2 \lfloor \delta k \rfloor} \big)
$$
that associates to each $\phi \in \sap_n$ with $\vert \globaljoin_\phi \vert \geq k$ the set 
$$
\Psi(\phi) = 
\Big\{ \mathscr{S}_\kappa(\phi): \kappa \subseteq \globaljoin_\phi \, , \,  \vert \kappa \vert = \lfloor \delta k \rfloor \Big\} \, ,
$$ 
where here we abuse notation and identify a subset of $\globaljoin_\phi$ with its set of indices under the given enumeration of $\globaljoin_\phi$.

Note that, for some constant $c > 0$ that is independent of $\delta \in (0,1/2)$ and for all $k \in \N$, 
$$
\big\vert \Psi(\phi) \big\vert =  {\vert \globaljoin_\phi \vert \choose \lfloor \delta k \rfloor} \geq  {k \choose \lfloor \delta k \rfloor} \geq c \, \delta^{-\delta k} \, .
$$
(In the latter displayed inequality, the factor of $k^{-1/2}$
that appears via Stirling's formula has been cancelled against other omitted terms. This detail is inconsequential and the argument is omitted.)

Note that, for any $\ga \in \saw_{n + 2 \lfloor \delta k \rfloor}$, the preimage $\Psi^{-1} (\gamma)$ is either the empty-set or a singleton. Indeed, if $\ga \in \Psi(\phi)$ for some $\phi \in \sap_n$ with $\vert \globaljoin_\phi \vert \geq k$, then $\phi$ may be recovered from $\ga$ as follows:
\begin{itemize}
\item the coordinate $x\big(\eastsouth(\phi)\big)$ equals $x\big( \ga_{n + 2 \lfloor \delta k \rfloor} \big)/2$;
\item the vertex $\eastsouth(\phi)$ is the lowest among the vertices of $\ga$ having the above $x$-coordinate;
\item setting $j \in [0,n]$ so that $\ga_j$ is this vertex, consider the non-self-avoiding walk $\ga_{[0,j]} \circ \reflect_{\ga_j}^2\big( \ga_{[j,n]} \big)$. This walk begins and ends at $0$. There are exactly $\lfloor \delta k \rfloor$ instances where the walk traverses an edge twice. In each case, the three-step journey that the walk makes in the steps preceding, during and following the second crossing of the edge follow three edges of a plaquette. Replace this journey by the one-step journey across the remaining edge of the plaquette, in each instance. The result is~$\phi$.
\end{itemize}

We may thus use the multi-valued map Lemma~\ref{l.mvm} to find a lower bound on the cardinality of $\rhssaw_{n + 2 \lfloor \delta k \rfloor}$. The resulting bound is precisely~(\ref{e.globaljoin.reduce}). This completes the derivation of Proposition~\ref{p.globaljoin}.   \qed

\medskip

\noindent{\bf Proof of Corollary~\ref{c.globaljoin}.}
Using
\begin{equation}\label{e.globalk}
\psap_n \Big(  \vert \globaljoin_\phi \vert \geq k \Big) = p_n^{-1} \cdot \# \Big\{ \phi \in \sap_n: \vert \globaljoin_\phi \vert \geq k \Big\} \, ,
\end{equation}
and Proposition~\ref{p.globaljoin},
we find that if $n \in \highpolynum_{\maceta}$ then 
$$
\psap_n \big(  \vert \globaljoin_\phi \vert \geq k \big) \leq c^{-1} \, 2^{- k \mu^{-2}/2}  n^{\Cpcp} 
$$
for $k \in \N$.
Set $\Cgjtwo = 2 \mu^2 \big( \log 2 \big)^{-1} \big( \Cgjone + \Cpcp + 1 \big)$. Then $n \geq c^{-1}$ implies that, for $k \geq \Cgjtwo \log n$,
$$
\psap_n  \big(  \vert \globaljoin_\phi \vert \geq k \big) \leq   c^{-1} n^{-\Cgjone} \, .
$$
Setting $\Cgjthree = c^{-1}$, we obtain Corollary~\ref{c.globaljoin}. 
 \qed

\subsection{Preparing for joining surgery: left and right polygons}\label{s.refpert.newwork}

When we prove Proposition~\ref{p.polyjoin.newwork} in the next subsection, polygons from a certain set $\saprward_{n-j}$
will be joined to others in another set $\sapellward_j$.
We think of the former as being joined to the latter on the right, so that the superscripts indicate a handedness associated to the joining.

Anyway, in this subsection, we specify these polygon sets, give a lower bound on their size in Lemma~\ref{l.polysetbound.newwork}, and then explain in Lemma~\ref{l.strongjoin} 
how there are plentiful opportunities for joining pairs of such polygons in a surgically useful way (so that the join may be detected by virtue of its global nature).

Let $\phi$ be a polygon. Recall from Definition~\ref{d.corners} the notation $\ymax(\phi)$ and $\ymin(\phi)$, as well as the height~$h(\phi)$ and width~$w(\phi)$.

\begin{definition}\label{d.leftright}
For $n \in 2\N$, let $\sapellward_n$ denote the set of {\em left} polygons $\phi \in \sap_n$
such that
\begin{itemize}
\item $h(\phi) \geq w(\phi)$ (and thus, by a trivial argument, $h(\phi) \geq n^{1/2}$),
\item and $y\big(\righttip(\phi) \big) \leq \tfrac{1}{2} \big( \ymin(\phi) + \ymax(\phi) \big)$. 
\end{itemize}
Let $\saprward_n$ denote the set of {\em right}  polygons $\phi \in \sap_n$
such that
\begin{itemize}
\item $h(\phi) \geq w(\phi)$.
\end{itemize}
\end{definition}
\begin{lemma}\label{l.polysetbound.newwork}
For $n \in 2\N$,
$$
\big\vert \sapellward_n \big\vert \geq \tfrac{1}{4} \cdot \big\vert \sap_n \big\vert \, \, \, 
\textrm{and}  \, \, \,
\big\vert \saprward_n \big\vert \geq  \tfrac{1}{2} \cdot \big\vert \sap_n \big\vert \, .
$$
\end{lemma}
\noindent{\bf Proof.} An  element $\phi \in \sap_n$
not in $\saprward_n$ is brought into this set
by right-angled rotation. 
If, after the possible rotation, it is not in $\sapellward_n$, it may 
brought there by reflection in the $x$-axis. \qed

\begin{definition}\label{d.globallymj}
A Madras joinable polygon pair $(\phi^1,\phi^2)$
is called globally Madras joinable if the junction plaquette of the join polygon $J(\phi^1,\phi^2)$ is a global join plaquette of  $J(\phi^1,\phi^2)$.
\end{definition}
Both polygon pairs in the upper part of Figure~\ref{f.madrasjoinable} are globally Madras joinable.

\begin{lemma}\label{l.strongjoin}
Let $n,m \in 2\N$ and let $\phi^1 \in \sapellward_n$ and $\phi^2 \in \saprward_m$.

Every value 
\begin{equation}\label{e.kvalues}
k \in \Big[ y \big(\righttip(\phi^1)\big)  \, , \,  y\big(\righttip(\phi^1)\big)  +  \min \big\{ n^{1/2}/2,m^{1/2} \big\}  - 1 \Big]
\end{equation}
is such that $\phi^1$ and some horizontal shift of $\phi^2 + ke_2$ is globally Madras joinable. 

Write $\strongjoin_{(\phi^1,\phi^2)}$ for the set of $\vec{u} \in \Z^2$ such that 
the pair $\phi^1$ and $\phi^2 + \vec{u}$ is globally Madras joinable.
Then
$$
  \big\vert \strongjoin_{(\phi^1,\phi^2)} \big\vert \geq   \min \big\{n^{1/2}/2,m^{1/2} \big\}  \, .
$$
\end{lemma}

  \begin{figure}
    \begin{center}
      \includegraphics[width=0.75\textwidth]{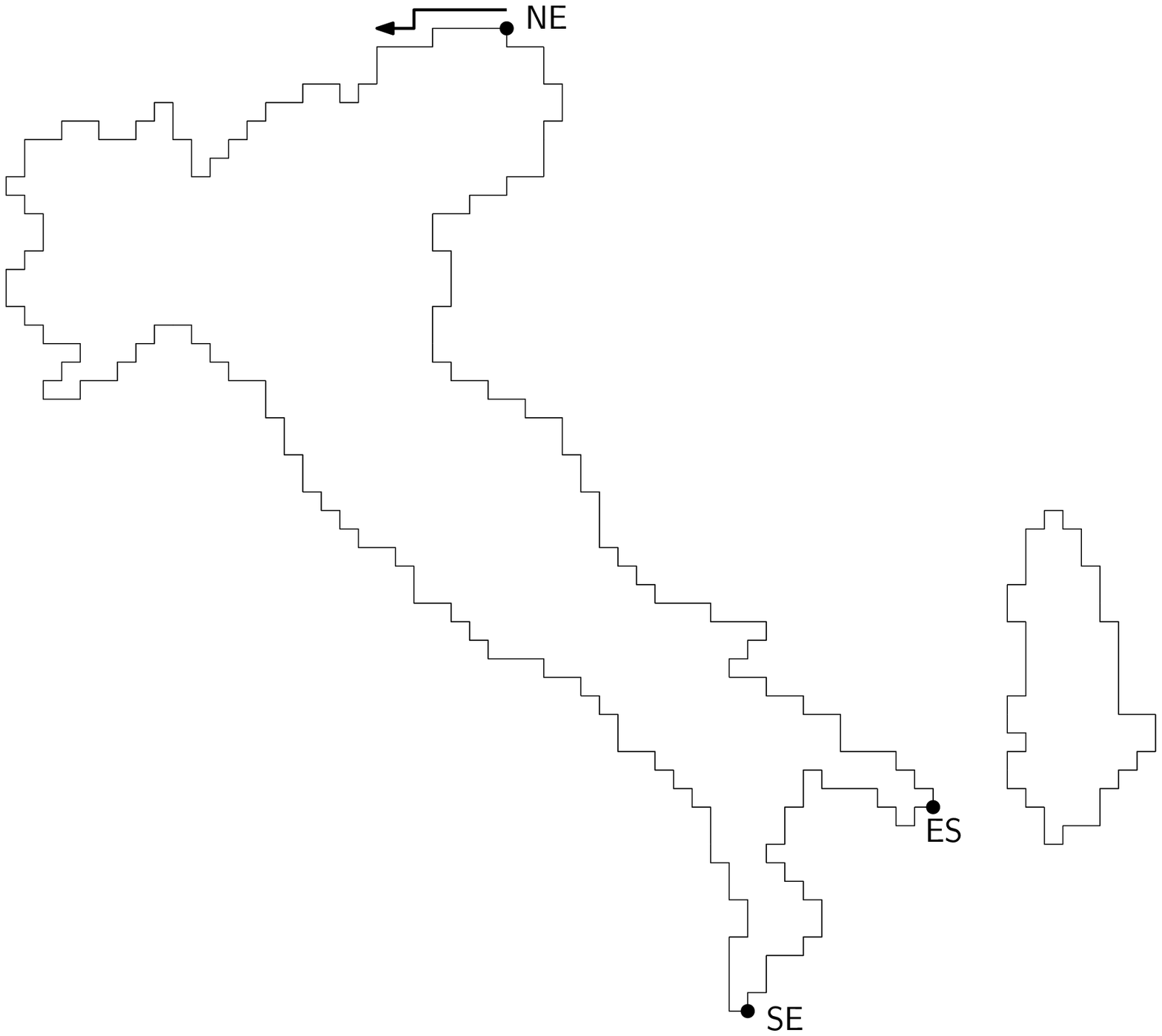}
    \end{center}
    \caption{Illustrating Lemma~\ref{l.strongjoin} with two polygons, $\mathsf{ITA}$ and $\mathsf{ALB}$. Denoting the polygons' lengths by $n$ and $m$, note that the translates of  $\mathsf{ITA}$ and $\mathsf{ALB}$ with $\northeast = 0$ belong to $\sapellward_n$ and $\saprward_m$. 
The start of the counterclockwise tour of $\mathsf{ITA}$ from $\northeast(\mathsf{ITA})$ 
dictated by convention is marked by an arrow.    
    The depicted polygons are globally Madras joinable after a horizontal shift of $\mathsf{ALB}$. This remains the case if this polygon is shifted upwards by at most two units or downwards by at most sixteen.}
    \label{fig:albita}
  \end{figure}

\noindent{\bf Proof.} Recall that since $\phi^1 \in \sap_n$
and $\phi^2 \in \sap_m$, $\ymax(\phi^1) = \ymax(\phi^2) = 0$.
Note that whenever $k \in \Z$ is such that the two intervals 
$$
\big[ \ymin(\phi^1) , 0 \big] \, \, \, \textrm{and} \, \, \, k + \big[ \ymin(\phi^2) , 0 \big]
$$ 
intersect, there is some horizontal displacement $j \in \Z$ such that $\phi^1$ and $\phi^2 + (j,k)$ are Madras joinable. Note also that $\phi^1 \in \sapellward_n$ satisfies $\ymin(\phi^1)  \leq - n^{1/2}$. Choices of $k \in [ - n^{1/2}  ,  - 1 ]$ thus produce polygon pairs $\big( \phi^1 , \phi^2 + k e_2 \big)$ whose first element contains vertices that are more northerly than any of the second, and which are Madras joinable after a horizontal shift of the second element of the pair. Moreover,  since $\phi^1 \in \sapellward_n$, $y\big( \righttip(\phi^1) \big)$ is at most $- n^{1/2}/2$.  Thus, any value of $k$ in~(\ref{e.kvalues}) lies in $[ - n^{1/2}  ,  - 1 ]$. In view of the Definition~\ref{d.globaljoin} of global join plaquette, in order to verify the lemma's first assertion, it remains to verify that the Madras joinable polygon pair associated by a suitable horizontal shift to any given value of $k$ in~(\ref{e.kvalues}) has the property that its second polygon contains all of the most easterly vertices in either of the two polygons in the pair. Here is the reason that this property holds: 
$\phi^2 \in \saprward_m$ implies that $h(\phi^2) \geq m^{1/2}$, and thus,
for such values of $k$, the $y$-coordinate of $\righttip(\phi^1)$
is shared by a vertex in $\phi^2 + ke_2$; from this, we see that, when this right polygon is shifted horizontally to be Madras joinable with $\phi^1$, it is forced to intersect the half-plane bordered on the left by the vertical line through $\righttip(\phi^1)$ and thus to verify the claimed property. 

The second assertion of the lemma is an immediate consequence of the first.  \qed

\subsection{Polygon joining is almost injective: deriving Proposition~\ref{p.polyjoin.newwork}}

We begin
by reducing this result to the next lemma.

\begin{lemma}\label{l.rgjp.newwork}
For any $\maceta > 0$, there is a constant $\Cpoly = \Cpoly(\maceta) > 0$ such that, for $n  \in 2\N$ satisfying $n + 16 \in \highpolynum_{\Cpcp}$, 
$$
  p_{n + 16}  \geq \frac{1}{\Cpoly \log n} \sum_{j \in 2\N \, \cap \, [2^{i-1},2^i]} p_j p_{n-j} (n-j)^{1/2} \, ,
$$
where $i \in \N$ is chosen so that $n \in 2\N \cap [2^{i},2^{i+1}]$. 
\end{lemma}

\noindent{\bf Proof of Proposition~\ref{p.polyjoin.newwork}.} 
The lemma coincides with the proposition when the two instances of $n + 16$ in its statement are replaced by $n$. To infer the proposition from the lemma (with a relabelling of the value of $\Cpoly(\maceta)$), it is thus enough to establish two things.
First, for some constant $c > 0$, $p_n \geq c \, p_{n+16}$ whenever $n \in 2\N$.  Second, that for $\maceta,\varphi > 0$, there exists $n_0(\maceta,\varphi)$ such that if $n \in 2\N$ satisfies  $n \geq n_0$ and $n + 16 \in \highpolynum_\maceta$, then
$n \in \highpolynum_{\maceta + \varphi}$.
To verify the first statement, recall from \cite[Theorem 7.3.4(c)]{MS93} that $\lim_{n \in 2\N} p_{n+2}/p_n = \mu^2$. Thus, $p_n \geq \mu^{-16} p_{n+16}/2$ for all $n$ sufficiently high; the value of $c > 0$ may be decreased from $\mu^{-16}/2$ if necessary in order that all $n \in 2\N$ satisfy the stated bound.
Given that the first statement holds, 
the second is confirmed by noting that if $n+16 \in \highpolynum_{\maceta}$, then $p_{n+16} \geq (2n)^{-\maceta} \mu^{n+16}$ whenever $n \in 2\N$ satisfies $n \geq 16$.  \qed

\medskip

\noindent{\bf Proof of Lemma~\ref{l.rgjp.newwork}.}
Set 
$A = \bigcup \, \sapellward_j \times \saprward_{n-j}$
where the union is taken over indices 
$j \in 2\N \, \cap \, [2^{i-1},2^i]$ for which $n - j$ is positive. Set $B = \sap_{n + 16}$.

We construct a multi-valued map $\Psi:A \to \mathcal{P}(B)$.
Consider a generic domain point 
 $(\phi^1,\phi^2) \in \sapellward_j \times \saprward_{n-j}$ where $j \in 2\N \, \cap \, [2^{i-1},2^i]$ satisfies $n - j > 0$. 
This point's image under $\Psi$ is defined to be 
 the collection of length-$(n+16)$ polygons 
formed by Madras joining $\phi^1$ and $\phi^2 + \vec{u}$ as $\vec{u}$ ranges over the set $\strongjoin_{(\phi^1,\phi^2)}$ specified in Lemma~\ref{l.strongjoin}.

 Since $j \geq 2^{i-1}$ and $n \leq 2^{i+1}$, we have that $j \geq (n-j)/3$.
By Lemma~\ref{l.strongjoin},
$$
  \Big\vert  \Psi\Big( (\phi^1,\phi^2) \Big) \Big\vert   =  \big\vert \strongjoin_{(\phi^1,\phi^2)} \big\vert \geq \min \big\{ j^{1/2}/2 , (n-j)^{1/2}  \big\} \geq \tfrac{1}{2} \, 3^{-1/2} (n-j)^{1/2} \, .
$$

Applying Lemma~\ref{l.polysetbound.newwork}, we learn that the number of arrows in $\Psi$ is at least
\begin{equation}\label{e.totalarrow}
3^{-1/2}  \, 2^{-4} \sum_{j = 2^{i-1}}^{2^i} p_j p_{n-j}  (n-j)^{1/2} \, ,
\end{equation}
where note that the summand is zero if $j$ is odd.


For a constant $\Ch > 0$, denote by 
$$
H_{n + 16} = \Big\{ \phi \in \sap_{n + 16} : \big\vert \globaljoin_\phi \big\vert \geq  \Ch \log (n+16) \Big\} 
$$
the set of length-$(n + 16)$ polygons with a {\em high} number of global join plaquettes.

We now fix the value of $\Ch$ with which we will work.
Recalling that $n + 16 \in \highpolynum_\maceta$, we may apply Corollary~\ref{c.globaljoin} with $\Cgjone = 1$ 
and the role of $n$ played by $n + 16$.
Set $\Ch > 0$ equal to the value of $\Cgjtwo$ determined via the corollary by this choice of $\Cgjone$ and the value of~$\maceta$.  Since $\psap_{n+16} \big( \vert \globaljoin_\Gamma \vert \geq \Ch \log (n+16) \big) =  \vert H_{n+16} \vert  p_{n + 16}^{-1}$, we find that
\begin{equation}\label{e.hnc}
 \big\vert H_{n + 16} \big\vert \leq \Cgjthree (n + 16)^{-1}  p_{n + 16}  \, .
\end{equation}

The set of preimages under $\Psi$ of a given 
$\phi \in \sap_{n+16}$ may be indexed by the junction plaquette 
associated to the Madras join polygon $J(\phi^1,\phi^2)$ that equals $\phi$. Recalling Definition~\ref{d.globallymj}, this plaquette is a global join plaquette of $\phi$. Thus,  
\begin{equation}\label{e.phisap}
\vert \Psi^{-1}(\phi) \vert \leq \vert \globaljoin_\phi \vert
 \, \, \,
\textrm{for any $\phi \in \sap_{n + 16}$} \,  
.
\end{equation}
 Since $\vert \globaljoin_\phi \vert \leq n + 16$ for all $\phi \in \sap_{n + 16}$, we find that
\begin{equation}\label{e.maxsap}
 \max \big\{ \big\vert \Psi^{-1}(\phi) \big\vert : \phi \in \sap_{n + 16} \big\} \leq n + 16 \, .
\end{equation}

The proof of Lemma~\ref{l.rgjp.newwork} will be completed by considering two cases.

\medskip

 In the first case,
\medskip
\begin{itemize}
\item at least one-half of the arrows in $\Psi$ point to elements of $H_{n+16}$ \, ;
\end{itemize}
\medskip
in the second case, then, at least one-half of these arrows point to elements of
$\sap_{n + 16} \setminus H_{n + 16}$.

\medskip

\noindent{{\bf The first case.}} The inequality
\begin{equation}\label{e.firstcase}
 \big\vert H_{n + 16} \big\vert \, \cdot \, \max \Big\{ \vert \Psi^{-1}(\phi) \vert: \phi \in H_{n + 16} \Big\} \geq \tfrac{1}{2} \cdot  3^{-1/2}  \, 2^{-4} \sum_{j = 2^{i-1}}^{2^i} p_j \, p_{n-j}  (n-j)^{1/2} 
\end{equation}
holds because the left-hand side is an upper bound on the number of arrows in $\Psi$ arriving in $H_{n + 16}$, which is at least one-half of the total number of arrows; and the latter quantity is at least the right-hand side. Applying~(\ref{e.hnc}) and~(\ref{e.maxsap}),
$$
 p_{n + 16} \geq  3^{-1/2} \, 2^{-5} \Cgjthree^{-1}     \sum_{j = 2^{i-1}}^{2^i}  p_j p_{n-j}(n-j)^{1/2}   \, ,
$$ 
so that Lemma~\ref{l.rgjp.newwork} is obtained in the first case.

\medskip

\noindent{{\bf The second case.}}
The quantity
$$
 \Big\vert \, \sap_{n + 16}  \setminus H_{n+16} \, \Big\vert \, \cdot \, \max \Big\{  \big\vert \Psi^{-1}(\phi) \big\vert : \phi \in  \sap_{n + 16}  \setminus H_{n+16} \Big\} $$
is an upper bound on the number of arrows incoming to $\sap_{n + 16}  \setminus H_{n+16}$. In the case that we now consider, the displayed quantity  is thus in view of the total arrow number lower bound~(\ref{e.totalarrow})
 at least $3^{-1/2}  \, 2^{-5} \sum_{j = 2^{i-1}}^{2^i} p_j p_{n-j}  (n-j)^{1/2}$.

 By~(\ref{e.phisap}), for $\phi \in \sap_{n + 16} \setminus H_{n+16}$, $\vert \Psi^{-1}(\phi) \vert \leq \Ch \log (n+16)$.
 Thus,
\begin{equation}\label{e.casetwo}
 p_{n + 16} \geq \Ch^{-1} \big( \log (n+16) \big)^{-1}  3^{-1/2}  \, 2^{-5} \sum_{j = 2^{i-1}}^{2^i} p_j p_{n-j}  (n-j)^{1/2}
 \, ,
\end{equation}
so that Lemma~\ref{l.rgjp.newwork} is proved in the second case also. \qed

\subsection{Inferring Theorem~\ref{t.polydev} from Proposition~\ref{p.polyjoin.newwork}}

In this subsection, we prove Theorem~\ref{t.polydev}.
The key step is the following, a consequence of Proposition~\ref{p.polyjoin.newwork}.

\begin{proposition}\label{p.hpn.propagate}
For each $\varphi > 0$, there exists $i_0(\varphi) > 0$ such that, for any $\delta > 0$,  $a \in (0,1)$ and $i \in \N$ for which $i \geq i_0(\varphi) + 2 (\log 2)^{-1} \big( 3\varphi^{-1} + 1 \big)  \log a^{-1}$,
the condition that
\begin{equation}\label{e.hyphigh}
 \Big\vert \, \highpolynum_{3/2 - \delta} \cap \big[2^{i-1},2^i \big] \, \Big\vert \, \geq \, a \cdot 2^{i-2}
\end{equation}
implies that
$$
 \Big\vert \, \highpolynum_{3/2 - 2\delta + \varphi} \cap \big[2^i,2^{i+1} \big] \, \Big\vert \, \geq \,  \tfrac{1}{8} a^2 \cdot 2^{i-1} \, .
$$
\end{proposition}
The proposition stipulates an unstable state of affairs should Theorem~\ref{t.polydev} fail, as we now outline.
Let $\delta > 0$ and assume that infinitely many dyadic scales are occupied by $\highpolynum_{3/2 - \delta}$ at a small but uniform fraction $a$. The proposition implies that the next scale up from any one of these scales (excepting finitely many) is occupied at a fraction at least $a^2/8$ by the more stringently specified set $\highpolynum_{3/2 - 3\delta/2}$.
We may iterate this inference over several scales and find that, a few dyadic scales higher from where we began, there is a tiny but positive fraction of indices that actually belong to  $\highpolynum_{\maceta}$ for a negative value of $\maceta$. This index set is of course known to be empty. Thus, the  $\highpolynum_{3/2 - \delta}$-occupancy assumption is found to be invalid, and Theorem~\ref{t.polydev} is proved. 

First, we prove Proposition~\ref{p.hpn.propagate} and then, in Proposition~\ref{p.hpn.contradict} and its proof, we implement a rigorous version of the argument just sketched.

\medskip 

\noindent{\bf Proof of Proposition~\ref{p.hpn.propagate}.}
We will prove the result with the choice 
\begin{equation}\label{e.varphiformula}
  i_0(\varphi) = \max \Big\{ 6, \varphi^{-1} \big( 18 + \tfrac{2}{\log 2} \log \Cpoly \big) , \tfrac{\log(2\pi)}{2\log 2} + \tfrac{1}{\log 2} \big( 4\varphi^{-1} + 3/2 \big) \log \big( 4 \varphi^{-1} + 1 \big)  \Big\} \, ,
\end{equation}
where the constant $\Cpoly$ is determined by  Proposition~\ref{p.polyjoin.newwork} with $\maceta = 3$.

Consider the map ${\rm sum}$ from the product of two copies of
$2\N \cap [2^{i-1},2^i]$
to $2\N \cap [2^i,2^{i+1}]$
that sends $(j,k)$ to $j+k$. 
Borrowing the language that we use for multi-valued maps, we think of the function as being specified by a collection of arrows $(j,k) \to j+k$ with target $j+k$. 
Fixing $\delta > 0$, we call an arrow $(j,k) \to j+k$ {\em high}
if both $j$ and $k$ are elements of 
$\highpolynum_{3/2 - \delta}$. 
We also identify the set $\targetmany \subseteq 2\N \cap [2^i,2^{i+1}]$
of {\rm sum} image points {\em targeted by many high arrows},
$$
 \targetmany = \Big\{ n \in  2\N \cap [2^i,2^{i+1}]: n \, \, \textrm{is the target of at least $\tfrac{1}{8} a^2 2^{i-2}$ high arrows}    \Big\} \, .
$$
We now argue that $\targetmany$ occupies a proportion of at least $a^2/8$ of the {\rm sum} image points.
\begin{lemma}\label{l.targetmany}
Suppose that the hypothesis~(\ref{e.hyphigh})
 holds for a given $a \in (0,1)$.
Then
$$
 \big\vert \targetmany \big\vert \geq \tfrac{1}{8} a^2 \cdot 2^{i-1} \, .
$$ 
\end{lemma}
\noindent{\bf Proof.} 
The  {\rm sum} image set  $2\N \cap [2^i,2^{i+1}]$
has $2^{i-1} + 1 \leq 2^i$ elements. Thus, the number of high arrows targeting elements of  $\targetmany^c$
is at most $2^i \cdot \tfrac{1}{8} a^2 2^{i-2} = a^2 \cdot 2^{2i - 5}$. 

By~(\ref{e.hyphigh}), the total number of high arrows is at least $a^2 \cdot 2^{2i- 4}$.
We learn then that at least half of these arrows target elements of $\targetmany$. However, no image point is the target of more than $2^{i-2} + 1 \leq 2^{i-1}$ arrows. Thus, there must be at least $a^2 \cdot 2^{2i - 5} \cdot 2^{1-i}$ elements of $\targetmany$. This completes the proof. \qed
 
\medskip

The next lemma and its proof concern a useful additive stability property of the high polygon number sets: for $\maceta > 0$, $\highpolynum_{\maceta} \times \highpolynum_{\maceta}$
is mapped by {\rm sum} into  $\highpolynum_{2\maceta}$. This fact is useful when we seek to apply Proposition~\ref{p.polyjoin.newwork}: the proposition will be useful when the right-hand side in~(\ref{e.polyjoin.newwork}) is in a suitable sense large, but, in this circumstance, the $\highpolynum$-stability property  verifies the proposition hypothesis that $n$ is an element of some $\highpolynum$-set.
 
\begin{lemma}\label{l.target.hpn}
If $n \in 2\N \cap [2^i,2^{i+1}]$ is the target of a high arrow, then~$n \in \highpolynum_{3 - 2\delta}$.
\end{lemma}
\noindent{\bf Proof.}
Let $(j,n-j)$ be the high arrow that targets $n$. 
By Proposition~\ref{p.basic}(2), $p_n \geq p_j p_{n-j}$.
Since $j$ and $n-j$ are elements of $\highpolynum_{3/2 - \delta}$, we have that 
\begin{equation}\label{e.mintwop}
\min \{ p_j \mu^{-j}, p_{n-j} \mu^{-(n-j)} \} \geq n^{-3/2 + \delta}  \, ,
\end{equation}
where we also used $\delta \leq 3/2$, which bound follows from $\highpolynum_{3/2 - \delta} \not= \emptyset$ and $p_n \leq \mu^n$ from Proposition~\ref{p.basic}(3).
 Thus, $p_n \geq n^{-3 + 2\delta} \mu^{2n}$. \qed

\medskip

Fix an element $n \in \targetmany$. Lemma~\ref{l.target.hpn}
shows that Proposition~\ref{p.polyjoin.newwork} is applicable with $\maceta = 3$. Thus,
\begin{equation}\label{e.pnlb}
 p_n \geq \frac{1}{\Cpoly \log n} \sum (n-j)^{1/2} p_j \, p_{n-j}  \, ,
\end{equation}
where the sum is taken over $j \in 2\N \, \cap \, [2^{i-1},2^i]$
such that $(j,n-j) \to n$ is a high arrow.

By its definition and the structure of the map ${\rm sum}$,  the set $\targetmany$ is in fact a subset of the even integers in the interval $[2^i + A, 2^{i+1} - A]$, where $A = \tfrac{1}{4} a^2 2^{i-2} - 2$.

Recalling (\ref{e.varphiformula}),
note that the hypotheses of the proposition entail that 
$i \geq 6 + 2(\log 2)^{-1} \log a^{-1}$.
Thus, $a \geq 2^{3 - i/2}$.
Let $j \in [2^{i-1},2^i]$.   
Since $n \in \targetmany$, we see from $n - j \geq A$
and this lower bound on $a$ that $n-j$ is at least
$a^2 2^{i-5}$.

Note then that
\begin{eqnarray*}
 p_n & \geq &  \frac{1}{\Cpoly \log n} a 2^{(i-5)/2} \sum  p_j \, p_{n-j} \\
  & \geq &     \frac{1}{\Cpoly \log n} a 2^{(i-5)/2}   \cdot  \tfrac{1}{8} a^2 2^{i-2}  \cdot n^{-3 + 2 \delta} \mu^{n} \\
  & \geq & 
   \frac{1}{2^{5/2+3+2} \Cpoly \log n} a^3 \cdot 2^{-3/2} n^{3/2} \cdot n^{-3 + 2 \delta} \mu^{n} = \frac{1}{2^{9} \Cpoly \log n} a^3  \cdot n^{-3/2 + 2 \delta} \mu^{n} \, .  
\end{eqnarray*}
In the first line, we sum over the same set of values of $j$ as we did in~(\ref{e.pnlb}). The inequality in this line is due to
(\ref{e.pnlb}) and $n - j \geq a^2 2^{i-5}$.
The second inequality is due to $n \in \targetmany$ and the bound~(\ref{e.mintwop}).  The final line arises because  $2^{i} \geq n/2$.

Finally, we argue that 
\begin{equation}\label{e.varphiargument}
 2^{-9} \Cpoly^{-1} a^3 \big(\log n \big)^{-1} \geq n^{-\varphi} \, .
\end{equation}
Combining with the last display, we learn that $p_n \geq n^{-3/2 + 2\delta - \varphi} \mu^n$, which is to say, $\thet_n \leq 3/2 - 2\delta + \varphi$. Recalling that this inference has been made for any element $n \in \targetmany$, we note that Lemma~\ref{l.targetmany}
completes the proof of Proposition~\ref{p.hpn.propagate}. 

It remains to confirm~(\ref{e.varphiargument}). This follows from two estimates: first, $\log n \leq n^{\varphi/2}$;
second, $2^{-9} \Cpoly^{-1} a^3 \geq n^{-\varphi/2}$.
The first is these follows by taking $k = \lceil 4\varphi^{-1} \rceil$ in $e^{n^{\varphi/2}} \geq \tfrac{1}{k!} n^{k \varphi/2}$
using $k! \leq (2\pi)^{1/2} \big( 4 \varphi^{-1} + 1 \big)^{4 \varphi^{-1} + 3/2}$ as well as $n \geq 2^i$ and $i \geq i_0(\varphi)$. The second follows from $n \geq 2^i$
and the hypothesised lower bound on~$i$. \qed

\begin{proposition}\label{p.hpn.contradict}
Let $\delta > 0$ and  $a \in (0,1)$. Let $i_0(\delta/2) > 0$ be  specified by Proposition~\ref{p.hpn.propagate}. 
Suppose that $i \in \N$ satisfies 
$$
i \geq 
 i_0(\delta/2) +  f(\delta,a) + 2 (\log 2)^{-1} \big( 6 \delta^{-1} + 1 \big) \log a^{-1} \, , 
$$ 
where  
$$
f(\delta,a) = 4 \big( \log 2 \big)^{-1}  \big( 6 \delta^{-1} + 1 \big)  \delta^{-( \log 3/2 )^{-1} \log 2}   \log \big( 8 a^{-1} \big) \, .
$$
Then
$$
 \Big\vert \, \highpolynum_{3/2 - \delta} \cap \big[2^{i-1},2^i \big] \, \Big\vert \, < \, a \cdot 2^{i-2} \, .
$$
\end{proposition}
\noindent{\bf Proof of Theorem~\ref{t.polydev}.} This is an immediate consequence of the proposition. \qed

\medskip

\noindent{\bf Proof of Proposition \ref{p.hpn.contradict}.}
Suppose the contrary. Let $\delta > 0$, $a \in (0,1)$ and 
\begin{equation}\label{e.imin}
i \geq i_0 \big( \delta/2 \big) \, + \, f(\delta,a) \, + \, 2 (\log 2)^{-1} \big( 6 \delta^{-1} + 1 \big)  \log a^{-1} 
\end{equation}
be such that~(\ref{e.hyphigh}) holds.

For $k \in \N$, set $h_k(a) = 8^{1 - 2^k} a^{2^k}$, and note for future reference that the recursion $h_{k+1}(a) = h_k(a)^2/8$ is satisfied with initial condition $h_0(a) = a$.

Note that the inequality 
\begin{equation}\label{e.hkdelta}
  2 (\log 2)^{-1} \big( 6 \delta^{-1} + 1 \big) \log h_k(a)^{-1} > f(\delta,a) + 2 (\log 2)^{-1}  \big( 6 \delta^{-1} + 1 \big)  \log a^{-1}  
\end{equation} 
  is satisfied for all sufficiently high $k \in \N$. Let $K \in \N$
  denote the smallest value of $k \in \N$ that satisfies this bound. Our choice of the function $f$ entails that $K > 1 + \big( \log 3/2 \big)^{-1}  \log \delta^{-1}$.

We claim that, for any $k \in [0,K]$,
\begin{equation}\label{e.k.hpn}
 \Big\vert \, \highpolynum_{3/2 - (3/2)^k \delta} \cap \big[2^{i+k-1},2^{i+k} \big] \, \Big\vert \, \geq \,  h_k(a) \cdot 2^{i-k-2} \, .
\end{equation}
This bound may be proved by induction on $k$, where the
base case $k=0$ is validated by our assumption that~(\ref{e.hyphigh}) holds.
Proposition~\ref{p.hpn.propagate} yields the bound for $k \in [1,K]$ as we now explain.
Suppose that the bound holds for index $k-1$. We seek to apply Proposition~\ref{p.hpn.propagate}
with the roles of $\delta$, $\varphi$, $i$ and $a$ being played by $(3/2)^{k-1} \delta$, $(3/2)^{k-1} \delta/2$,  $i + k -1$
and $h_{k-1}(a)$. The conclusion of the proposition 
 is then indeed the bound at index $k$, but we must check that the proposition's hypotheses are valid. 
Beyond the inductively supposed hypothesis, we must confirm that
\begin{equation}\label{e.ikhyp}
 i + k -1 \geq  i_0 ( \varphi ) \, + \, 
  2 (\log 2)^{-1} \big( 3 \varphi^{-1} + 1 \big) \log h_{k-1}(a)^{-1}   \, ,
\end{equation}
with $\varphi =  (3/2)^{k-1} \delta/2$. Note that since $k \in [0,K]$,
the inequality~(\ref{e.hkdelta}) is violated at index $k-1$.
Note that $i + k - 1 \geq i$ may be bounded below by~(\ref{e.imin}) and then, by this violation, by
$i_0(\delta/2) +   2 (\log 2)^{-1} \big( 6 \delta^{-1} + 1 \big) \log h_{k-1}(a)^{-1}$. Since $\varphi \geq \delta/2$ and $i_0$ is a decreasing function of $\varphi > 0$, we do indeed verify~(\ref{e.ikhyp}) and so complete the derivation of~(\ref{e.k.hpn}) for each $k \in [0,K]$. 

Consider now (\ref{e.k.hpn}) with $k = K$:
we see that 
$\big\vert \highpolynum_{\zeta} \cap \big[2^{i+
 K-1},2^{i+K} \big]  \big\vert  \geq   h_K(a) \cdot 2^{i-K-2} > 0$ where $\zeta = 3/2 - (3/2)^K \delta$.
However, by recalling that $K$ exceeds $\big( \log 3/2 \big)^{-1}  \log \delta^{-1}$ by more than one, we find $\zeta$ 
 to be negative, so that $p_n \leq \mu^n$ implies that $\highpolynum_\zeta = \emptyset$.
  By obtaining this contradiction, we have completed the proof of Proposition~\ref{p.hpn.contradict}. \qed

\section{Polygon joining in preparation for the proof of Theorem~\ref{t.thetexist}}\label{s.pjprep}

In the final Chapter~\ref{c.aboveonehalf}, we will prove Theorem~\ref{t.thetexist} by an approach that combines the polygon joining technique with other ideas.
The approach for polygon joining used in the proof is an adaptation of the constructs just used to prove Theorem~\ref{t.polydev}. In this section, we present the notation and results used for this adaptation. We emphasise that these ideas are used only in the proof of Theorem~\ref{t.thetexist}. The reader may choose to
read this section after the preceding one, when the ideas being adapted are fresh, or may wish to
 defer the reading until the start of Section~\ref{s.almost} when this content will be needed.

In outline, we revisit the joining of left and right polygons used in the proof of Lemma~\ref{l.rgjp.newwork}.
We refine the definitions of these polygon types (in the first subsection), and also of the joining mechanism for them. The polygons that are formed under the new, stricter, mechanism will be called 
{\em regulation global join} polygons: in the second subsection, these polygons are defined and their three key properties stated.
The remaining three subsections prove these properties in turn.

\subsection{Left and right polygons revisited}\label{s.refpert}
Recall the left and right polygon sets $\sapellward_n$
and $\saprward_n$ from Definition~\ref{d.leftright}.
We will refine these notions now, specifying new left and right polygon sets $\sapell_n$
and $\sapr_n$. Each of a subset of its earlier counterpart. 

In order to specify $\sapell_n$, we now define a {\em left-long} polygon.
Let $\phi \in \sap_n$. 
Again we employ the notationally abusive parametrization $\phi:[0,n] \to \Z^2$ such that $\phi_0 = \phi_n = \northeast(\phi)$ and $\phi_1 = \northeast(\phi) - e_1$. If $j \in [0,n]$ is such that $\phi_j = \base(\phi)$, note that $\phi$ may be partitioned into two paths $\phi_{[0,j]}$ and $\phi_{[j,n]}$. (This is not the division into two paths used in the proof of Proposition~\ref{p.globaljoin}: the outward journey is now to $\base(\phi)$, not $\eastsouth(\phi)$.) The two paths are edge-disjoint, and the first lies to the left of the second: indeed, writing~$H$ for the horizontal strip $\big\{ (x,y) \in \R^2: y\big(\base(\phi)\big) \leq y \leq y\big(\northeast(\phi)\big) \big\}$, the set $H \setminus \cup_{i=0}^{j-1} [\phi_i,\phi_{i+1}]$ has two connected components; one of these, the right component, contains $H \cap \big\{ (x,y) \in \R^2: x \geq C \big\}$ for large enough $C$; and the union of the edges $[\phi_i,\phi_{i+1}]$, $j \leq i \leq n-1$, excepting the points $\phi_j = \base(\phi)$ and $\phi_n = \northeast(\phi)$, lies in the right component. It is thus natural to call $\phi_{[0,j]}$ the left path, and $\phi_{[j,n]}$, the right path. We call $\phi$ {\em left-long}  if $j \geq n/2$ and {\em right-long} if $j \leq n/2$.  

For $n \in 2\N$, let $\sapell_n$ denote the set of left-long elements of
$\sapellward_n$.
We simply set $\sapr_n$ 
equal to $\saprward_n$.
Note that $\mathsf{ITA}$ in Figure~\ref{fig:albita} is an element of $\sapell_n$.

\begin{lemma}\label{l.polysetbound}
For $n \in 2\N$,
$$
\big\vert \sapell_n \big\vert \geq \tfrac{1}{8} \cdot \big\vert \sap_n \big\vert \, \, \, 
\textrm{and}  \, \, \, 
\big\vert \sapr_n \big\vert \geq  \tfrac{1}{2} \cdot \big\vert \sap_n \big\vert \, .
$$
\end{lemma}
\noindent{\bf Proof.} The second assertion has been proved in Lemma~\ref{l.polysetbound.newwork}. In regard to the first, write the three requirements for an element $\phi \in \sap_n$ to satisfy $\phi \in \sapell_n$
in the order: $h(\phi) \geq w(\phi)$; $h$ is left-long; and  $y\big(\righttip(\phi) \big) \leq \tfrac{1}{2} \big( \ymin(\phi) + \ymax(\phi) \big)$.
 These conditions may be satisfied as follows.
\begin{itemize}
\item The first property may be ensured by a right-angled counterclockwise rotation if it does not already hold. 
\item It is easy to verify that, when a polygon is reflected in a vertical line, the right path is mapped to become a sub-path of the left path of the image polygon. Thus, a right-long polygon maps to a left-long polygon under such a reflection.
(We note incidentally the asymmetry in the definition of left- and right-long: this last statement is not always true {\em vice versa}.) A polygon's height and width are unchanged by either horizontal or vertical reflection, so the first property is maintained if a reflection is undertaken at this second step.
\item The third property may be ensured if necessary by reflection in a horizontal line.  The first property's occurrence is not disrupted for the reason just noted. Could the reflection disrupt the second property? 
When a polygon is reflected in a horizontal line, $\northeast$ and $\base$ in the domain  map to $\base$ and $\northeast$ in the image. The left path maps to the reversal of the left path, (and similarly for the right path). Thus, any left-long polygon remains left-long when it is reflected in such a way. We see that the second property is stable under the reflection in this third step.
\end{itemize}
Each of the three operations has an inverse, and thus  $\# \sapell_n \geq \tfrac{1}{8} \# \sap_n$. \qed

\subsection{Regulation global join polygons: three key properties}\label{s.threekey}

\begin{definition}\label{d.rgjp}
Let $k,\ell \in 2\N$
satisfy $k/2 \leq \ell \leq 35 k$. 
Let $\dubbub_{k,\ell}$ denote the set of {\em regulation global join} polygons (with length pair $(k,\ell)$), whose elements are formed by
Madras joining the polygon pair $\big(\phi^1,\phi^2 + \vec{u}\big)$, where
\begin{enumerate}
\item  $\phi^1 \in \sapell_k$;
\item $\phi^2 \in \sapr_\ell$;
\item $\ymax\big( \phi^2 + \vec{u} \big) \in \big[ y\big(\eastsouth(\phi^1)\big),y\big(\eastsouth(\phi^1)\big) +  \lfloor k^{1/2}/{10} \rfloor - 1 \big]$;
\item and $\big(\phi^1,\phi^2 + \vec{u} \big)$ is Madras joinable.
\end{enumerate}
Let $\dubbub$ denote the union of the sets $\dubbub_{k,\ell}$ over all such choices of $(k,\ell) \in 2\N \times 2\N$.
\end{definition}

Note that  Lemma~\ref{l.strongjoin} implies that whenever a polygon pair is joined to form an element of $\dubbub$, this pair is globally Madras joinable. 

The three key properties of regulation polygons are now stated as propositions.

\subsubsection{An exact formula for the size of $\dubbub_{k,\ell}$} 

\begin{proposition}\label{p.rgjpnumber}
Let $k,\ell \in 2\N$ satisfy $k/2 \leq \ell \leq 35 k$.  For any $\phi \in \dubbub_{k,\ell}$, there is a unique choice of 
$\phi^1 \in \sapell_k$, $\phi^2 \in \sapr_\ell$ and $\vec{u} \in \Z^2$ for which $\phi = J\big(\phi^1,\phi^2 + \vec{u} \big)$. 
We also have that
$$
\big\vert \dubbub_{k,\ell} \big\vert = \lfloor  k^{1/2}/{10} \rfloor \cdot  \big\vert \sapell_k \big\vert \, \big\vert \sapr_\ell \big\vert \, . 
$$
\end{proposition}

\subsubsection{The law of the initial path of a given length in a random regulation polygon}

We show that initial subpaths in regulation polygons are distributed independently of the length of the right polygon.

\begin{definition}
Write $\psapleft{n}$ for the uniform law on $\sapell_n$.
\end{definition}

\begin{proposition}\label{p.leftlong}
Let $k,\ell \in 2\N$ satisfy $k/2 \leq \ell \leq 35 k$. 
Let $j \in \N$ satisfy $j \leq k/2$. For any $\phi \in \sapell_k$,
$$
\psap_{k+\ell + 16} \Big( \Gamma_{[0,j]} = \phi \, \Big\vert \, \Gamma \in \dubbub_{k,\ell} \Big)
  =   \psapleft{k} \Big( \Gamma_{[0,j]} = \phi \Big) \, .
$$
\end{proposition}

\subsubsection{Regulation polygon joining is almost injective}

We have a counterpart to Proposition~\ref{p.polyjoin.newwork}.

\begin{proposition}\label{p.polyjoin.aboveonehalf}
For any $\Theta > 0$, there exist  $c > 0$ and $n_0 \in \N$ such that the following holds. Let  $i \geq 4$ be an integer and let $n \in 2\N \cap \big[ 2^{i+2},2^{i+3} \big]$ satisfy $n \geq n_0$.
Suppose that 
 $\mathsf{R}$
is a subset of $2\N \cap [2^i,2^{i+1}]$
that contains an element  $\kmac$ such that
$\max \big\{ \theta_{\kmac} , \theta_{n - 16 - \kmac} \big\} \leq \Theta$.
Then
$$
 \bigg\vert \, \bigcup_{j \in \mathsf{R}} \dubbub_{j,n-16-j} \, \bigg\vert \geq c \, \frac{n^{1/2}}{\log n} \sum_{j  \in \mathsf{R}} p_j \, p_{n-16-j}  \, .
$$
\end{proposition}

Section~\ref{s.pjprep} concludes with three subsections consecutively devoted to the proofs of these propositions.

\subsection{Proof of Proposition~\ref{p.rgjpnumber}} 
\begin{figure}
    \begin{center}
      \includegraphics[width=0.6\textwidth]{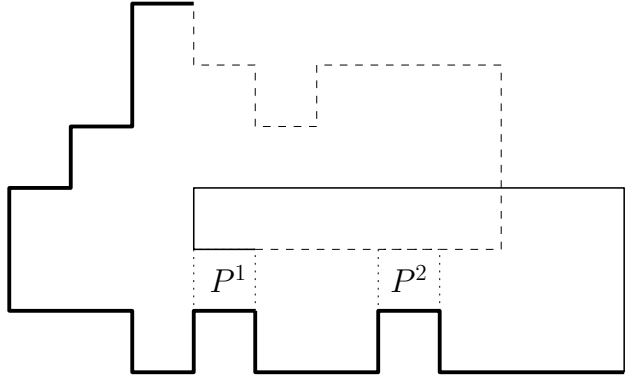}
    \end{center}
    \caption{Suppose that two global join plaquettes $P^1$ and $P^2$ of a polygon $\phi$ are such that, in the outward journey along $\phi$, from $\northeast(\phi)$ to $\eastsouth(\phi)$, $P^1$ is encountered before $P^2$, and that the order is maintained on the return. Then the return journey must visit some vertex twice. The outward journey is depicted in bold. The crossing of an edge of $P^1$ by the return forces the dashed future of the return into a bounded component of the complement of the existing path.}\label{f.planar}
  \end{figure}

We need only prove the first assertion, the latter following directly. To do so, it is enough to determine the junction plaquette associated to the Madras join that forms $\phi$. 
By Lemma~\ref{l.strongjoin}, any such polygon pair $(\phi^1,\phi^2 + \vec{u})$ is globally Madras joinable. Recalling Definition~\ref{d.globallymj}, the associated junction plaquette is thus a global join plaquette of~$\phi$. By Lemma~\ref{l.globaloneedge}, each global join plaquette has one horizontal edge traversed in the outward journey along $\phi$ from $\northeast(\phi)$ to $\eastsouth(\phi)$, and one traversed on the return journey. A short exercise discussed in Figure~\ref{f.planar} that  invokes planarity shows that the set of $\phi$'s global join plaquettes is totally ordered under a relation in which the upper element is both reached earlier on the outward journey and later on the return. From this, we infer that the map that sends a global join plaquette $P$ of $\phi$ to the length of the polygon in $\phi \, \Delta \, P$ that contains $\northeast(\phi) = 0$ is injective. Since this length must equal $k+8$ for any admissible choice of 
$\phi^1 \in \sapell_k$, $\phi^2 \in \sapr_\ell$ and $\vec{u} \in \Z^2$, this choice is unique, and we are done. \qed

\medskip 

Note further that the set of join locations stipulated by the third condition in Definition~\ref{d.rgjp} is restricted to an interval of length of order~$k^{1/2}$. The restriction causes the formula in Proposition~\ref{p.rgjpnumber} to hold. It may be that many elements of~$\sapell_k$ and~$\sapr_\ell$ have heights much exceeding $k^{1/2}$, so that, for pairs of such polygons, there are many more than an order of  $k^{1/2}$ choices of translate for the second element that result in a globally Madras joinable polygon pair. The term {\em regulation} has been attached to indicate that a specific rule has been used to produce elements of $\dubbub$ and to emphasise that such polygons do not exhaust the set of polygons that we may conceive as being globally joined.

\subsection{Deriving Proposition~\ref{p.leftlong}}

The next result is needed for this derivation.

\begin{lemma}\label{l.leftlongbasic}
Consider a globally Madras joinable polygon pair $(\phi^1,\phi^2)$.
Let $j \in \N$ be such that $\phi^1_j$ equals $\southeast(\phi^1)$.  The join polygon $J(\phi^1,\phi^2)$ has the property that the initial subpaths $\phi^1_{[0,j]}$ and $J(\phi^1,\phi^2)_{[0,j]}$ coincide.
\end{lemma}
\noindent{\bf Proof.}
Since  $(\phi^1,\phi^2)$ is globally Madras joinable, the northeast vertex of the join polygon coincides with $\northeast(\phi^1)$. As such, it suffices to confirm that $\phi^1$'s left path $\phi^1_{[0,j]}$ forms part of $J(\phi^1,\phi^2)$.

 Review Figure~\ref{f.madrascases}
and suppose that some case other than IIcii or IIIcii occurs. Note by inspection that, when $J(\phi^1,\phi^2)$ is formed, either one or two edges are removed from $\phi^1$, and that every endpoint~$v$ of these edges has the property that the horizontal line segment extending to the right from~$v$ intersects the vertex set of~$\phi^1$ only at~$v$. 
Note however that every vertex $\phi^1_i$, $0 \leq i \leq j-1$, (i.e., the vertices of the left path of $\phi^1$ except its endpoint $\southeast(\phi^1)$) has the property that the line segment extending rightwards from the vertex does have a further intersection with the vertex set of~$\phi^1$. Thus, the above vertex $v$ must equal $\phi^1_i$ for some $i \geq j$ (i.e., must lie in the right path of $\phi^1$). The initial subpath $\phi^1_{[0,j]}$ is thus unchanged by the joining operation and is shared by $J(\phi^1,\phi^2)$.

  \begin{figure}
    \begin{center}
      \includegraphics[width=0.7\textwidth]{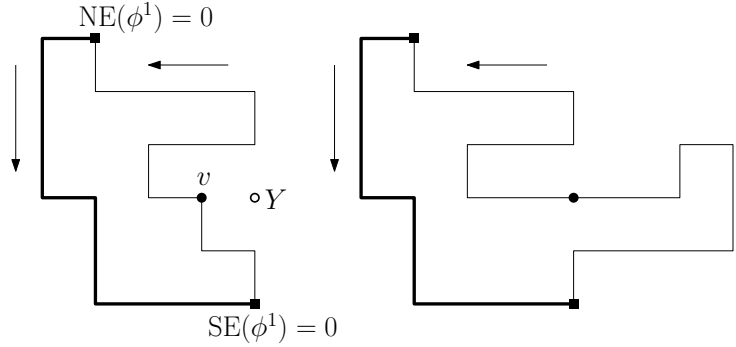}
    \end{center}
    \caption{A case in the proof of Lemma~\ref{l.leftlongbasic}. The polygon $\phi^1$ depicted on the left lies in $\sap_n$ for some $n \in 2\N$. The polygon is indexed counterclockwise from $\northeast(\phi^1)$ as the arrows indicate.
   Recalling Figure~\ref{f.madrascases}, the circle in the left sketch, and the disk in the right, marks the point $Y$ used in the formation of the right sketch $\modify{\tau}$ with $\tau = \phi^1$, (during the process of Madras joining with an undepicted~$\phi^2$). 
    Case~IIcii obtains for $\phi^1$. The left sketch disk marks a vertex $v$ as discussed in the proof.  Emboldened is the left path of $\phi^1$, which passes undisturbed by the surgery from the left sketch to the right to form an initial subpath of the postsurgical polygon.}\label{f.caseiicii}
  \end{figure}

Cases IIcii and IIIcii remain. Figure~\ref{f.caseiicii} depicts an example of the first of these. In either case, note that at least one endpoint~$v$ of one of the two removed edges enjoys the above property and thus belongs to the right path of $\phi^1$. Note by inspecting Figure~\ref{f.madrascases} that $\southeast(\phi^1)$ cannot lie in the vicinity of~$v$ that is altered by the joining operation. Thus, once again,   
 $\phi^1_{[0,j]}$ is shared by $J(\phi^1,\phi^2)$.  \qed

\medskip

\noindent{\bf Proof of Proposition~\ref{p.leftlong}.} By the first assertion in Proposition~\ref{p.rgjpnumber}, we find that the conditional distribution of 
$\psap_{k+\ell+16}$ given that~$\Gamma \in \dubbub_{k,\ell}$ has the law of the output in this procedure:
\begin{itemize}
\item select a polygon according to the law $\psapleft{k}$;
\item independently select another uniformly among elements of $\sapr_\ell$;
\item choose one of the $\lfloor \litnum \, k^{1/2} \rfloor$ vertical displacements specified in Definition~\ref{d.rgjp} uniformly at random, and shift the latter polygon by this displacement;
\item translating the displaced polygon by a suitable horizontal shift, Madras join the first and the translated second polygon to obtain the output.
\end{itemize} 
Denote by $\phi^1$ the element of $\sapell_k$ selected in the first step. This polygon is left-long; let $m \geq k/2$ satisfy $\phi^1_m = \southeast(\phi^1)$. By Lemma~\ref{l.leftlongbasic}, the output polygon has an initial subpath that coincides with the left path $\phi^1_{[0,m]}$ of~$\phi^1$. Since $j \leq k/2 \leq m$, we see that the initial length~$j$ subpath of the output polygon coincides with $\phi^1_{[0,j]}$. Thus, its law is that of $\Gamma_{[0,j]}$ where the polygon $\Gamma$ is distributed according to~$\psapleft{k}$. \qed

\subsection{Proof of Proposition~\ref{p.polyjoin.aboveonehalf}}

This is a variation on the argument for Lemma~\ref{l.rgjp.newwork}, which is the principal component of Proposition~\ref{p.polyjoin.newwork}.
As such, our argument will need to invoke the key hypothesis of Corollary~\ref{c.globaljoin} that  $n \in \highpolynum_\maceta$ for a suitable choice of $\maceta > 0$. We begin by pointing out that this hypothesis  is satisfied for any choice of $\maceta$ exceeding $2\Theta$ (provided that $n$ exceeds a level $n_0$ determined by~$\maceta$). That this is true follows from 
the straightforward $p_n \geq p_{n-16}$, Proposition~\ref{p.basic}(2)'s $p_{n-16} \geq p_{\kmac} p_{n-16-\kmac}$ and the hypothesised $\max \big\{ \thet_{\kmac},\thet_{n-16-\kmac} \big\} \leq \Theta$. (In fact, this inference is in essence the $\highpolynum$-additive stability property that is discussed before Lemma~\ref{l.target.hpn}.)

Recall that $\mathsf{R}$ (alongside $i \geq 4$, $n$, $\Theta$, and $\kmac$) is specified in Proposition~\ref{p.polyjoin.aboveonehalf}'s statement.
Adapting the proof of Lemma~\ref{l.rgjp.newwork}, we consider  a multi-valued map $\Psi:A \to \mathcal{P}(B)$. 
We take
$A = \bigcup_{j \in \mathsf{R}} \, \sapell_j \times \sapr_{n-16-j}$  and   $B =  \cup_{j \in \mathsf{R}} \,  \dubbub_{j,n-16-j}$ (which is a subset of $\sap_n \cap \dubbub$).
Note  the indices concerned in specifying $A$ are non-negative  because $\max \mathsf{R} \leq 2^{i+1}$, $n \geq 2^{i+2}$ and $i \geq 3$ ensure that $n - 16 -j \geq 0$ whenever $j \in \mathsf{R}$. (Moreover,  the choice $(k,\ell) = (j,n-16-j)$ in Definition~\ref{d.rgjp} does satisfy the bounds stated there. We use $i \geq 4$ to verify that $k/2 \leq \ell$.)
We specify $\Psi$ so that  the generic domain point 
 $(\phi^1,\phi^2) \in \sapell_j \times \sapr_{n-16-j}$, $j \in \mathsf{R}$, is mapped to  
 the collection of length-$n$ polygons 
formed by Madras joining $\phi^1$ and $\phi^2 + \vec{u}$ as $\vec{u}$ ranges over the subset of $\strongjoin_{(\phi^1,\phi^2)}$ specified for the given pair $(\phi^1,\phi^2)$ by conditions $(3)$ and $(4)$ in Definition~\ref{d.rgjp}. (That this set is indeed a subset of $\strongjoin_{(\phi^1,\phi^2)}$ is noted after Definition~\ref{d.rgjp}.)

We have that $\big\vert  \Psi\big( (\phi^1,\phi^2) \big) \big\vert = j^{1/2}/10$. 
 Since $j \geq 2^{i}$ and $n \leq 2^{i+5}$, we have that $j \geq  2^{-5}n$, so that $\big\vert  \Psi\big( (\phi^1,\phi^2) \big) \big\vert \geq \tfrac{1}{10} 2^{-5/2} n^{1/2}$.

Applying Lemma~\ref{l.polysetbound}, we learn that the number of arrows in $\Psi$ is at least
\begin{equation}\label{e.totalarrow.late}
 2^{-4} \cdot \tfrac{1}{10} 2^{-5/2} n^{1/2} \sum_{j \in \mathsf{R}} p_j p_{n-j}  \, .
\end{equation} 

Respecify the set $H_n = H_n(\Ch)$ from the proof of Lemma~\ref{l.rgjp.newwork}
so that it is given by  $\big\{ \phi \in B : \vert  \globaljoin_\phi \vert \geq  \Ch \log n  \big\}$.
We now set the value of $\Ch > 0$.
Specify a parameter $\maceta = 2\Theta + 1$, and recall that 
$n \in \highpolynum_\maceta$ is known provided that $n \geq n_0(\maceta)$. We set the parameter $\Ch$ equal to the value  $\Cgjtwo$ furnished by Corollary~\ref{c.globaljoin} applied with $\Cgjone$ taking the value $2\Theta + 1$ and for the given value of 
$\maceta$ (which happens to be the same value). This use of the corollary also provides a value of $\Cgjthree$, and,  
since $\psap_n \big( \vert \globaljoin_\Gamma \vert \geq \Ch \log n \big) = p_n^{-1} \vert H_n \vert$, we learn that
$\vert H_n \vert \leq  \Cgjthree n^{-(2\Theta + 1)}  p_n$.

As in the proof of Lemma~\ref{l.rgjp.newwork}, there are two cases to analyse, the first of which is specified by at least one-half of the arrows in $\Psi$ pointing to $H_n$.
In the first case, we have
\begin{eqnarray*}
 & &  n \times  \Cgjthree n^{-(2\Theta + 1)}  p_n  \, \geq \,  \max \big\{ \vert \Psi^{-1}(\phi) \vert: \phi \in H_{n} \big\} \times  \big\vert H_{n} \big\vert 
  \, \geq \, \tfrac{1}{2} \# \big\{ \,  \textrm{arrows in $\Psi$} \, \big\} \\
  & \geq & 
 \tfrac{1}{2} \cdot \tfrac{1}{10} 2^{-13/2} n^{1/2} \sum_{j \in \mathsf{R}} p_j p_{n-16-j} \, \geq \, \alpha \, n^{1/2} p_{\kmac} p_{n-16- \kmac} \, \geq \, \alpha \, n^{1/2} n^{-2\Theta} \mu^{n -16} \, .  
\end{eqnarray*} 
where $\alpha = \tfrac{1}{10} 2^{-15/2}$.
The consecutive inequalities are due to: $\big\vert \Psi^{-1}(\phi) \big\vert \leq n$ for  $\phi \in \sap_n$; the first case obtaining; the lower bound
(\ref{e.totalarrow.late}); the membership of $\mathsf{R}$ by the hypothesised index $\kmac$; and the hypothesised properties of $\kmac$. We have found that $p_n \geq c \, n^{1/2} \mu^n$ with $c = \mu^{-16} \alpha \Cgjthree^{-1}$. Since $p_n \leq \mu^n$ by Proposition~\ref{p.basic}(3), the first case ends in contradiction if 
we stipulate that the hypothesised lower bound $n_0(\Theta)$ on $n$ is at least $c^{-2} = 100 \cdot 2^{15} \mu^{32} \Cgjthree$ (whose $\Theta$-dependence is transmitted via $\Cgjthree$).

The argument for~(\ref{e.phisap}) equally yields $\vert \Psi^{-1}(\phi)\vert \leq \vert \globaljoin_\phi \vert$ for $\phi \in \sap_n$. We find then that, when the second case obtains,
$$
 \vert B \vert \geq  \big\vert B \setminus H_n \big\vert  \geq    \Ch^{-1} \big( \log n \big)^{-1} \times
\tfrac{1}{2} \cdot \tfrac{1}{10} 2^{-13/2} n^{1/2} \sum_{j \in \mathsf{R}} p_j p_{n-16-j}  \, .
$$
Recalling the definition of the set $B$ completes the proof of Proposition~\ref{p.polyjoin.aboveonehalf}. \qed

\section{General dimension: deriving Theorem~\ref{t.threed}}\label{s.threed}

Theorem~\ref{t.threed} is a general dimensional counterpart to Theorem~\ref{t.polydev}. 
The result is expressed in terms of
 $3$-edge self-avoiding walk model
in order to eliminate technicalities. We now adapt the proof of Theorem~\ref{t.polydev} from Section~\ref{s.three}
to prove Theorem~\ref{t.threed}. 
 
First we recall that in the introduction we asserted the existence of $\muhat =  \lim_{n \in 2\N} \barp_n^{1/n}$. 
In fact, we also have that $ \lim_{n \in \N} \carp_n^{1/n}$
exists and equals this value.  Moreover, we asserted that
$c_n \geq \muhat^n$ for $n \in \N$. 
These claims follow from a slight variation of Proposition~\ref{p.basic}  for the $3$-edge self-avoiding walk model.  The proof given in this article, and the Hammersley-Welsh unfolding argument needed to show Proposition~\ref{p.basic}(4), can be adapted with only minor changes, a discussion of which we omit.

We respecify the index set from Definition~\ref{d.highpolynum}
so that,
for $\maceta > 0$, 
$\highpolynum_\maceta = \big\{ n \in 2\N : \barp_n \geq n^{- \maceta} \muhat^n \big\}$.
With this change made, we may state an analogue of Proposition~\ref{p.polyjoin.newwork}.
\begin{proposition}\label{p.polyjoin3d}
Let $d \geq 2$. For any $\maceta > 0$, there is a constant $\Cpoly = \Cpoly(\maceta) > 0$ such that, for $n \in 2\N \cap \highpolynum_{\Cpcp}$,
$$
 \barp_n \geq \frac{1}{\Cpoly \log n} \sum_{j \in 2\N \, \cap \, [2^{i-1},2^i]} (n-j)^{1 - 1/d} \barp_j \, \barp_{n-j}  \, ,
$$
where $i \in \N$ is chosen so that $n \in 2\N \cap [2^i,2^{i+1}]$. 
\end{proposition}

Our proof follows closely the structure of Theorem~\ref{t.polydev}'s, with this section having five subsections each of which plays a corresponding role  to its counterpart in Section~\ref{s.three}. As such, the first four subsections lead to the proof of Proposition~\ref{p.polyjoin3d}, and, in the final one, Theorem~\ref{t.threed} is proved from the proposition.

\subsection{Specifying the joining of two polygons}

\begin{definition}
Let $\phi$ be a polygon.  A {\em join} edge of $\phi$ is a nearest neighbour edge in $\Z^d$ that is traversed precisely two times by $\phi$, with one crossing in each direction. 
\end{definition}
Note that when a join edge is removed from a polygon, two polygons result.

\begin{definition}
For $I$ any $d-1$ element subset of $\{1,2,\cdots,d\}$, let $\proj_I:\Z^d \to \Z^d$ denote the projection onto the axial plane containing the vectors $e_i$ for $i \in I$. Write $\proj = \proj_{2,\cdots,d}:\Z^d \to \{ 0 \} \times \Z^{d-1}$. 
\end{definition}
We present an analogue of Madras' local surgery procedure for polygon joining, respecifying the join $J(\phi,\phi')$ of two polygons. 
For $n,m \in 2\N$, let $\phi$ and $\phi'$ be polygons of lengths $n$ and $m$. Suppose that $\proj(\phi) \cap \proj(\phi') \not= \emptyset$. 
If $\phi'$ is translated to a sufficiently high coordinate in the $e_1$-direction,   the vertex sets of $\phi$ and the $\phi'$-translate will be disjoint. Make such a translation and then translate $\phi'$ backwards in the negative $e_1$-direction  until the final location at which the translate's vertices remain disjoint from those of $\phi$. When two polygons have such a relative position, they are in a situation comparable to being Madras joinable; we call them {\em simply} joinable. Note that there exists an $e_1$-oriented nearest neighbour edge of~$\Z^d$
with one endpoint in one of the polygon vertex sets and the other in the other such set.
Let $e$ denote the maximal edge among these (according to some fixed ordering of nearest neighbour edges of~$\Z^d$). Define $J(\phi,\phi')$ to be the length~$n + m + 2$ polygon formed by following the trajectory of $\phi$ until an endpoint of $e$ is encountered, crossing $e$, following the whole trajectory of $\phi'$, recrossing $e$, and completing the trajectory of~$\phi$. The edge $e$ will be called the {\em junction} edge in this construction.

Let the left vertex $\leftv(\phi)$ of a polygon $\phi$ be the lexicographically minimal vertex in $\phi$ (so that $\leftv(\phi)$ has minimal $e_1$-coordinate). For our purpose, the left vertex will play a counterpart role  to that of the northeast vertex  in the setting of the proof of Theorem~\ref{t.polydev}. For $n \in 2\N$, we define $\saptd_n$ to be the set of polygons of length~$n$ whose left vertex is the origin. Note that $\barp_n = \# \saptd_n$. We adopt a convention for parametrizing any $\phi \in \saptd_n$, so that we may write $\phi:[0,n] \to \Z^d$.
The polygon $\phi$ visits the origin, possibly more than once.
For each such visit, we may record two lists of length $n$ of the vertices consecutively visited by $\phi$, with one having the opposite orientation of time to the other. 
Among these finitely many lists, we select the one of minimal lexicographical order, and take $\phi:[0,n] \to \Z^d$ to be this list.

We also define the right vertex $\rightv(\phi)$ of a polygon $\phi$ be the  lexicographically maximal vertex in $\phi$ (which is thus one of maximal $e_1$-coordinate).

\subsection{Global join edges are sparse}

Here we present analogues to Definition~\ref{d.globaljoin}, Proposition~\ref{p.globaljoin} and this result's proof.

\begin{definition}
Let $n \in 2\N$ and let $\phi \in \saptd_n$.  A join edge $e$ of $\phi$ is called {\em global} if the two polygons comprising $\phi$ without $e$ may be labelled  $\phi^\ell$ and $\phi^r$ in such a way that 
\begin{itemize}
\item every vertex of maximal $e_1$-coordinate in $\phi^\ell \cup \phi^r$ belongs to $\phi^r$;
\item and every vertex of minimal $e_1$-coordinate in $\phi^\ell \cup \phi^r$ belongs to $\phi^\ell$.
\end{itemize}
Write $\globaljoin_\phi$ for the set of global join edges of the 
polygon~$\phi$.
\end{definition}

Proposition~\ref{p.globaljoin}'s most direct analogue holds.
\begin{proposition}\label{p.globaljoin3d}
There exists $c > 0$ such that, for $n \in 2\N$ and any $k \in \N$,
$$
  \# \Big\{ \phi \in \saptd_{n} : \big\vert \globaljoin_\phi \big\vert \geq k \Big\} \leq  c^{-1} 2^{-k \muhat^{-2}/2}    \muhat^n   \, .
$$
\end{proposition}
\noindent{\bf Proof.} Let $\phi \in \saptd_n$. The right vertex $\rightv(\phi)$ replaces $\eastsouth(\phi)$ in this proof. 
Set $j \in [0,n]$ so that $\phi_j = \rightv(\phi)$; write $\phi^1 =  \phi_{[0,j]}$ and $\phi^2 = \phi_{[j,n]}$, and denote by $\reflect_z$ reflection in the hyperplane with normal vector~$e_1$ that passes through $z \in \Z^d$. Then set
$$
\mathscr{S}(\phi) = \phi^1 \circ \reflect_{\rightv(\phi)}(\phi^2) \, .
$$
The proof proceeds as Proposition \ref{p.globaljoin}'s did.
In the present case, $\mathscr{S}$ is mapping $\saptd_n$
into the set of bridges of length~$n$, namely into the set of $3$-edge length-$n$ self-avoiding walks whose minimal and maximal $e_1$-coordinate is attained at the start and the endpoint. In fact, bridges are usually defined so that this attainment is unique at one of the extremes. The addition of an extra edge at the end can assure this uniqueness. Using the classical upper bound on bridge number discussed in the proof of~(\ref{e.rhswalkbound}) (which is easily adapted to the $3$-edge case), we see that the image set $\mathscr{S}(\saptd_n)$
has cardinality at most $\muhat^{n+1}$. 
 
Another modification is made in specifying the alternative walks
$\mathscr{S}_\kappa(\phi)$, where $\kappa \subseteq \{ 1,\cdots,r\}$. Denote the $i\textsuperscript{th}$ global join edge by $g_i \in \globaljoin_\phi$. Because $g_i$ is a join edge, it is traversed twice by~$\phi$. The removal of $g_i$ from $\phi$
results in two polygons $\phi^\ell$ and $\phi^r$. Since $g_i$ is a global join edge, we have that 
$\leftv(\phi) \in \phi^\ell$ and $\rightv(\phi) \in \phi^r$.
Thus, $\phi_{[0,j]}$ contains one traversal of the edge $g_i$ made by $\phi$,
and $\phi_{[j,n]}$ the other. 
In light of this, we may describe the local modication for index~$i$, analogous to the three step subpath around $P^i$ in the original proof. 
Consider the edge in  $\mathscr{S}(\phi)$ that is the reflected image of $g_i$.
The modified walk is specified by insisting that it follows the course of  $\mathscr{S}(\phi)$ until this particular edge is traversed; after the traversal, the new walk performs a two-step move, crossing back over the edge, and then back again; after this, it pursues the remaining trajectory of   $\mathscr{S}(\phi)$.
  Note that no edge is traversed more than three times in the resulting definition of $\mathscr{S}_\kappa(\phi)$, and also that, when the post-concatenation part of this walk is reflected back, the outcome is a walk with an edge being traversed four times associated to any surgery, rendering the location of these surgeries detectable. (This property serves to explain our use of $3$-edge self-avoiding walks.) \qed  
 
\subsection{Left-right polygon pairs}

We now specify a counterpart  $\lrpp_{(k,\ell)}$  to the set $\sapellward_k \times \saprward_\ell$ of pairs of left and right polygons. The definition-lemma-definition-lemma structure of the counterpart Subsection~\ref{s.refpert.newwork} is maintained.


For $i \in [1,d]$, the $e_i$-span of a polygon is defined to be the difference between the maximum and the minimum $e_i$-coordinate of elements of the vertex set of $\phi$. Write  $\xspan(\phi)$ for the $e_1$-span of $\phi$.
 
\begin{definition} 
Let $k,\ell \in 2\N$.
A $(k,\ell)$ left-right polygon pair $(\phi,\phi')$
is an element of the union of $\saptd_k \times \saptd_\ell$
and $\saptd_\ell \times \saptd_k$ such that
\begin{itemize}
\item $\xspan(\phi) \geq \xspan(\phi')$, 
\item and $\vert \proj(\phi') \vert \geq (3d)^{-(1 - 1/d)} l(\phi')^{1 - 1/d}$, where $l(\phi')$ denotes the length of~$\phi'$.
\end{itemize}
Write $\lrpp_{(k,\ell)}$ for the set of  $(k,\ell)$ left-right polygon pairs.
\end{definition}

\begin{lemma}\label{l.polysetbound3d}
For $k,\ell \in 2\N$,
$$
\Big\vert \, \lrpp_{(k,\ell)} \,  \Big\vert \, \geq \, \tfrac{1}{2} d^{-2} \cdot \big\vert \saptd_k \big\vert \cdot  \big\vert \saptd_\ell \big\vert   
\, .
$$
\end{lemma}
\noindent{\bf Proof.}
Begin with a polygon pair $(\phi,\phi') \in \saptd_k \times \saptd_\ell$. By relabelling and reordering the pair, we may ensure that $\xspan(\phi)$
is at least the span of $\phi'$ in any of the~$d$ coordinates.
This relabelling is responsible for a factor of $(2d)^{-1}$ on the right-hand side. Now further relabel the coordinate axes in regard only to $\phi'$ in order that the cardinality of $\proj(\phi')$ for the relabelled $\phi'$
is at least as large as the cardinality associated to the other $d-1$ directions. This entails the appearance of a further factor of $d^{-1}$ on the right-hand side. It is easily seen that $\vert V(\phi') \vert \geq \tfrac{1}{3d}l(\phi')$ (since $\phi'$ is a $3$-edge self-avoiding polygon), and, as such, the lemma will be proved if we can show that the maximum cardinality among the axial projections of $\proj(\phi')$ is at least $\vert V(\phi') \vert^{1 - 1/d}$.
 This is a consequence of the Loomis-Whitney inequality~\cite{LoomisWhitney}: for any $d \geq 2$ and any finite $A \subseteq \Z^d$, the product of the cardinality of the projections of $A$ onto each of the $d$ axial hyperplanes is at least $\vert A \vert^{d-1}$. Also using the arithmetic-geometric mean inequality, one of the projections has cardinality at least $\vert A \vert^{1 - 1/d}$. \qed

\medskip

For polygons $\phi$ and $\phi'$, write $\strongtdjoin_{(\phi,\phi')}$ for the set of $\vec{u} \in \Z^d$ such that the pair 
$\big( \phi,\phi' + \vec{u} \big)$
is {\em globally joinable}, which is to say that 
\begin{itemize}
\item
$\big( \phi,\phi' + \vec{u} \big)$ is simply joinable;
\item any element of minimal $e_1$-coordinate  among the vertices of $\phi$ or $\phi' + \vec{u}$ is a vertex in $\phi$;
\item and any element of maximal $e_1$-coordinate  among the vertices of $\phi$ or $\phi' + \vec{u}$ is a vertex in $\phi' + \vec{u}$.
\end{itemize}

\begin{lemma}\label{l.claimmin}
For 
$k,\ell \in 2\N$, if 
$(\phi,\phi') \in \lrpp_{(k,\ell)}$, 
then 
$$
\big\vert \strongtdjoin_{(\phi,\phi')} \big\vert \geq (3d)^{-(1 - 1/d)} \min \big\{ k^{1-1/d} , \ell^{1-1/d} \big\} \, .
$$
\end{lemma}
\noindent{\bf Proof.}
Suppose that  $\vec{u} \in \{ 0 \} \times \Z^{d-1}$ is such that  $\proj(\phi' + \vec{u})$ contains the vertex $\proj\big(\rightv(\phi)\big)$. Due to our choice of $\phi'$, such $\vec{u}$ number 
at least the right-hand side in the inequality in the lemma's statement.  It is thus sufficient to argue that, for any such $\vec{u}$, there exists $k \in \Z$ for which the pair $(\phi,\phi' + \vec{u} + k e_1)$ is globally joinable. Recall that there is a maximal choice of $k \in \Z$ such that the polygon pair   $\big(\phi,\phi' + \vec{u} + (k-1) e_1 \big)$ is not vertex disjoint but the pair  $\big(\phi,\phi' + \vec{u} + k e_1 \big)$ is vertex disjoint, and that the pair  $\big(\phi,\phi' + \vec{u} + k e_1 \big)$ is said to be simply joinable. Specifying $k \in \Z$ in this way, we must also check that  $\vec{u} + k e_1$ verifies the second and third conditions for membership of~$\strongtdjoin_{(\phi,\phi')}$. To do  this, note that 
the polygon $\phi' + \vec{u} + k e_1$ contains a vertex with an $e_1$-coordinate strictly exceeding that of any vertex in $\phi$, so that the third condition is confirmed. Moreover, this same fact implies the second condition, because   the $e_1$-span of $\phi$ exceeds that of $\phi'$. We have proved Lemma~\ref{l.claimmin}. \qed

\subsection{Proof of Proposition~\ref{p.polyjoin3d}}
We begin with Lemma~\ref{l.rgjp.newwork}'s counterpart.
\begin{lemma}\label{l.rgjp.newwork.hd}
For any $\maceta > 0$, there is a constant $\Cpoly = \Cpoly(\maceta) > 0$ such that, for $n  \in 2\N$ satisfying $n + 2 \in \highpolynum_{\Cpcp}$, 
$$
  \hat{p}_{n + 2}  \geq \frac{1}{\Cpoly \log n} \sum_{j \in 2\N \, \cap \, [2^{i-1},2^i]} \hat{p}_j \hat{p}_{n-j} (n-j)^{1 - 1/d} \, ,
$$
where $i \in \N$ is chosen so that $n \in 2\N \cap [2^{i},2^{i+1}]$. 
\end{lemma}
Lemma~\ref{l.rgjp.newwork.hd}  implies Proposition~\ref{p.polyjoin3d} as Lemma~\ref{l.rgjp.newwork} did Proposition~\ref{p.polyjoin.newwork}, with the assertion counterpart to \cite[Theorem 7.3.4(c)]{MS93} that $\lim_{n \in 2\N} \hat{p}_{n+2}/\hat{p}_n = \muhat^2$ being used. 
The cited result has a non-trivial proof, but the changes needed to the proof are trivial. We do not provide details, but mention that constructs in Section~\ref{s.six} would be used: the reader may wish to consider how changing one uniformly chosen type $II$ pattern to a type~$I$ pattern in an element of $\psaptd_{n+2}$ generates a law on $\saptd_n$ which is provably close to the uniform law $\psaptd_n$ in light of a counterpart to Kesten's pattern theorem~\cite[Theorem 1]{kestenone} for $3$-edge walks.

\noindent{\bf Proof of Lemma~\ref{l.rgjp.newwork.hd}.} This is similar to the proof of Lemma~\ref{l.rgjp.newwork}. We consider the multi-valued map $\Psi:A \to \mathcal{P}(B)$, where 
$$
A = \bigcup_{j \in 2\N \, \cap \, [2^{i-1},2^i]} \lrpp_{(n-j,j)} 
$$
and $B = \saptd_{n+2}$. The map $\Psi$
associates to each $(\phi^1,\phi^2) \in  \lrpp_{(n-j,j)}$, $j \in 2\N \, \cap \, [2^{i-1},2^i]$, the set of length-$(n+2)$ polygons 
formed by simply joining $\phi^1$ and $\phi^2 + \vec{u}$ for choices of $\vec{u}$ in $\strongtdjoin_{(\phi^1,\phi^2)}$.
That the image set may be chosen to be $B$ is due to the left vertex of any formed polygon being the origin (as the second property in the definition of globally joinable indicates).

The proof follows the earlier one. In place of the assertion leading to~(\ref{e.phisap}) that the junction plaquette associated to the join polygon of a globally Madras joinable polygon pair is a global join plaquette, we instead assert that  in the join polygon $J(\phi^1,\phi^2)$ of two concerned polygons $\phi^1$ and $\phi^2$, the junction edge is a global join edge (which statement is a trivial consequence of the definition of globally joinable). Note that Corollary~\ref{c.globaljoin} may be invoked with obvious notational changes, because Proposition~\ref{p.globaljoin3d} replaces Proposition~\ref{p.globaljoin}.  \qed

\subsection{Proof of Theorem~\ref{t.threed}}

Since Proposition~\ref{p.polyjoin3d}
varies from 
Proposition~\ref{p.polyjoin.newwork} 
only by the relabelling of some notation,
Theorem~\ref{t.threed} 
follows from the former proposition as Theorem~\ref{t.polydev}
did from the latter after evident notational changes are made.
There are two manifestations of the value of dimension in these arguments, in the proof of Lemma~\ref{l.target.hpn} and at the end of the proof of Proposition~\ref{p.hpn.contradict}, when $p_n \leq (d-1)\mu^n$ is applied for $d=2$. The proofs now work provided that the dyadic scale $i$ is sufficiently high. \qed

\part{Closing probability upper bounds via the snake method}\label{c.thesnakemethod}

We saw in the last part that closing probability upper bounds may be inferred from polygon number upper bounds via the relation~(\ref{e.closepc}). 
In this part, we leave aside this combinatorial perspective and instead study a probabilistic approach, the snake method, for finding closing probability upper bounds.  
After some general notation and definitions in Section~\ref{s.four}, we present the general apparatus of the snake method in Section~\ref{s.five} and then use the  method via Gaussian pattern fluctuation in Section~\ref{s.six} to prove the $n^{-1/2 + o(1)}$-closing probability  Theorem~\ref{t.closingprob}. 
Dimension $d$ takes any value at least two throughout Part~\ref{c.thesnakemethod}.

\section{Some generalities}\label{s.four}

\subsection{Notation}

\subsubsection{Path reversal}
For $n \in \N$ and a length~$n$ walk $\ga:[0,n] \to \Z^d$, the reversal $\reverse\gamma:[0,n] \to \Z^d$ of $\ga$ is given by $\reverse\ga_j = \ga_{n-j}$ for $j \in [0,n]$.

\subsubsection{Walk length notation} We write $\vert \ga \vert = n$ for the length of any $\ga \in \saw_n$.

\subsubsection{Maximal lexicographical vertex} Definition~\ref{d.corners} specified the northeast vertex $\northeast(\gamma)$ of a walk~$\gamma$ when $d=2$. We now extend the definition to $d \geq 3$. For such~$d$ and 
any finite $V \subset \Z^d$, we write $\northeast(V)$ for the lexicographically maximal vertex in $V$. The notation is extended to any walk  $\ga:[0,n] \to \Z^d$ by setting $\northeast(\ga)$ equal to $\northeast(V)$ with $V$ equal to the image of~$\gamma$. 

In order that the definitions when $d=2$ and $d \geq 3$ coincide, it is necessary to adopt the convention that when~$d=2$ the lexicographical ordering is specified with the second coordinate having priority over the first. 

Our choice of notation $\northeast$ certainly emphasises the case $d=2$. This emphasis is made because when arguments in this article are made in general dimension, there have almost no dimensional dependence, and the two-dimensional case is then worth focussing on, because 
 visualizing  constructions in this dimension may aid understanding.



\subsubsection{The two-part decomposition}

In the snake method, we represent any given walk $\gamma$
in a {\em two-part} decomposition. This consists of an ordered pair of walks $(\gamma^1,\gamma^2)$ that emanate from
a certain common vertex and that are disjoint except at that vertex.
The two walks are called the {\em first part} and the {\em second part}.
To define the decomposition, consider any walk $\gamma$ of length~$n$. We first mention that the common vertex is chosen to be the lexicographically maximal vertex $\northeast(\gamma)$ on the walk.
Choosing $j \in [0,n]$  
so that $\gamma_j = \northeast(\gamma)$, the walk $\gamma$ begins at $\ga_0$ and approaches $\northeast(\ga)$ along the subwalk $\gamma_{[0,j]}$,
and then continues to its endpoint $\ga_n$ along the subwalk $\gamma_{[j,n]}$. The reversal $\reverse\gamma_{[0,j]}$ of the first walk, and the second walk $\gamma_{[j,n]}$, form a pair of walks that emanate from $\northeast(\gamma)$. (When $j$ equals zero or $n$,
one of the walks is the length zero walk that only visits $\northeast(\gamma)$; in the other cases, each walk has positive length.) The two walks will be the two elements in the two-part decomposition; all that remains is to order them, to decide which is the first part. Associated to each walk is the list of vertices consecutively visited by the walk $\big( \gamma_j,\cdots,\gamma_0 \big)$ and  $\big( \gamma_j,\cdots,\gamma_n \big)$.
These lists may be viewed as elements in $\Z^{d(j+1)}$ and $\Z^{d(n+1-j)}$. 
The first part, $\gamma^1$, is chosen to be the walk in the pair whose list is lexicographically the larger; the second part, $\gamma^2$, is the other walk.

We use square brackets to indicate the two-part decomposition, writing $\ga = [\ga^1,\ga^2]$.

As a small aid to visualization, it is useful to note that if the first $\ga^1$ part of a two-dimensional walk $\ga$ for which $\northeast(\gamma) = 0$ has length~$j$, then $\gamma^1:[0,j] \to \Z^2$ satisfies 
\begin{itemize}
\item $\gamma^1_0 = 0$ and $\gamma^1_1 = - e_1$;
\item $y \big( \gamma^1_i \big) \leq 0$ for all $i \in [0,j]$;
\item $\gamma^1_i \not\in \N \times \{ 0 \}$ for any $i \in [1,j]$.
\end{itemize}

  \begin{figure}
    \begin{center}
      \includegraphics[width=0.3\textwidth]{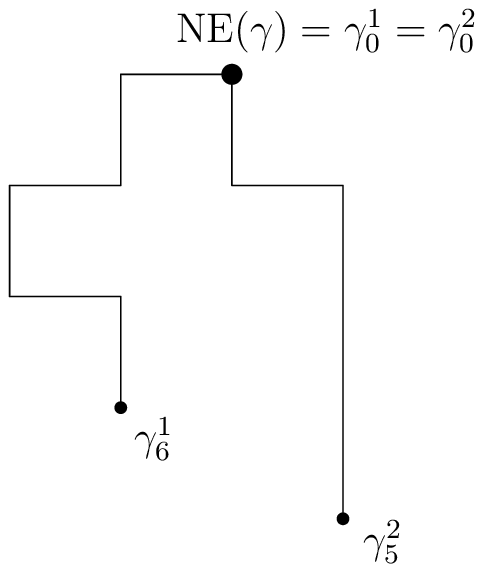}
    \end{center}
    \caption{The two-part decomposition of a walk of length eleven.
}\label{f.firstpart}
  \end{figure}

\subsubsection{Polygonal invariance}

The following trivial lemma will play an essential role. It is an important indication as to why polygons can be more tractable than walks.

\begin{lemma}\label{l.polyinv}
  For $n \in 2\N + 1$ and $j \in [1,n-1]$,  let $\chi:[0,n] \to \Z^d$ be a closing walk, and let $\chi'$ be the closing walk obtained from $\chi$ by the cyclic shift $\chi'(i) = \chi\big( j+i  \mod n+1 \big)$, $i \in [0,n]$. Then
  $$
  \psaw_n \Big(\Ga \, \textrm{is a translate of $\chi$} \Big) = 
  \psaw_n \Big(\Ga   \, \textrm{is a translate of $\chi'$} \Big) \, .
  $$
\end{lemma}
\subsubsection{Notation for walks not beginning at the origin}\label{s.sawfree}
Let $n \in \N$. We write $\sawfree_n$ for the set of self-avoiding walks~$\ga$ of length~$n$ (without stipulating the location~$\ga_0$). We further write $\saw_n^0$ for the subset of $\sawfree_n$ whose elements have northeast vertex at the origin. Naturally, an element $\ga \in \saw_n^0$ is said to close (and be closing) if $\dist \ga_n - \ga_0 \dist = 1$. The uniform law on $\saw_n^0$ will be denoted by $\psaw_n^0$. The sets $\saw_n^0$ and $\saw_n$ are in bijection via a clear translation; we will use this bijection implicitly.

\subsection{First parts and closing probabilities}
\subsubsection{First part lengths with low closing probability are rare}

\begin{lemma}\label{l.closecard}
Let $n \in 2\N + 1$ be such that, for some $\alphamac > 0$, $\psaw_n \big( \Gamma \closes \big) \geq n^{-\alphamac}$. For any $\deltamac > 0$,  the set of $i \in [0,n]$ for which
$$
 \# \Big\{ \gamma \in \saw^0_n: \vert \gamma^1 \vert = i \Big\}
 \geq n^{\alphamac + \deltamac}  \cdot \# \Big\{ \gamma \in \saw^0_n: \vert \gamma^1 \vert = i  \, , \, \gamma \closes \Big\}
$$
has cardinality at most $2 n^{1- \deltamac}$. 
\end{lemma}
\noindent{\bf Proof.} Note that $\psaw_n \big( \Gamma \closes \big) \geq n^{-\alphamac}$ implies that 
$$
 \big\vert \saw_n \big\vert \leq  n^{\alphamac}  \cdot \# \Big\{ \gamma \in \saw_n: \gamma \closes \Big\} \, .
$$
Note also that this inequality holds when $\saw_n$ is replaced by $\saw_n^0$.

We have that
$$
\# \Big\{ \gamma \in \saw^0_n: \gamma \closes \Big\} = \sum_{j=0}^n \# \Big\{ \gamma \in \saw^0_n: \gamma \closes \, ,  \, \vert \gamma^1 \vert = j \Big\}
$$
where, by Lemma~\ref{l.polyinv}, each term on the right-hand side has equal cardinality. Write $Q = Q_{\deltamac} \subseteq \{ 0,\ldots,n \}$ for the index set in the lemma's statement.
Note that 
\begin{eqnarray*}
  \big\vert \saw^0_n \big\vert & \geq & \vert Q \vert \cdot
 n^{\alphamac + \deltamac}  \cdot \tfrac{1}{n+1} \, \# \Big\{ \gamma \in \saw^0_n: \gamma \closes \Big\}  \\
 & \geq & \vert Q \vert \cdot  \tfrac{1}{2} \,
 n^{\alphamac -1 + \deltamac}   \, \# \Big\{ \gamma \in \saw^0_n: \gamma \closes \Big\} \, .
\end{eqnarray*}
Thus, $ \vert Q \vert \cdot \tfrac{1}{2} n^{\alphamac - 1 + \deltamac} \leq n^{\alphamac}$. \qed

\subsubsection{Possible first parts and their conditional closing probabilities}\label{s.possparts}

Let $d \geq 2$. For $n \in \N$, let $\fpart_n \subseteq \saw_n$ denote  the set of walks  $\gamma:[0,n] \to \Z^d$
whose lexicographically maximal vertex is $\gamma_0 = 0$.
We wish to view $\fpart_n$ as the set of possible first parts of walks $\phi \in \saw^0_m$ of some length $m$ that is at least $n$.
(In the two-dimensional case, we could be more restrictive in specifying $\fpart_n$, stipulating if we wish that any element $\ga$ satisfies $\ga_1 = -e_1$. What matters, however, is only that $\fpart_n$ contains all possible first parts.)

\begin{figure}
    \begin{center}
      \includegraphics[width=0.4\textwidth]{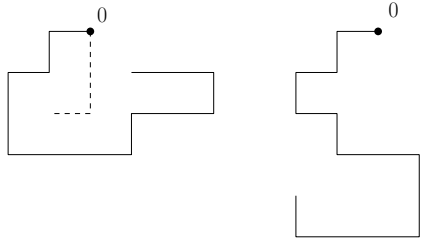}
    \end{center}
    \caption{{\em Left}: the bold $\phi \in \saw_{14}$ and dashed $\ga \in \saw_3$ are such that $[\phi,\ga]$ is a two-part decomposition. Note that $\phi \in \fpart_{14,14+3} \cap \fpart_{14,14+4}^c$. {\em Right}: An element of $\cap_{m = 1}^\infty\fpart_{14,14+m}$.}\label{f.twophi}
  \end{figure}

Note that, as Figure~\ref{f.twophi} illustrates, for given $m > n$, only some elements of $\fpart_n$ appear as such first parts, and we now record notation for the set of such elements (whatever the value of $d \geq 2$). Write $\fpart_{n,m} \subseteq \fpart_n$ for the set of $\gamma \in \fpart_n$ for which there exists an element $\phi \in \saw_{m-n}$ (necessarily with $\northeast(\phi) = 0$) such that $[\gamma,\phi]$ is the two-part decomposition of some element $\chi \in \saw_m$ with $\northeast(\chi) = 0$. 

 In this light, we now define the conditional closing probability 
$$
 q_{n,m}:\fpart_{n,m} \to [0,1] \,   , \, \, \,  \, q_{n,m}(\gamma) = \psaw^0_m \Big( \Gamma \closes \, \Big\vert \, \Gamma^1 = \gamma \Big) \, ,
$$ 
where here $m,n \in \N$ satisfy $m > n$; note that since $\gamma \in \fpart_{n,m}$, the event in the conditioning on the right-hand side occurs for some elements of $\saw_m$, so that the right-hand side is well-defined. 

We also identity a set of first parts with {\em high} conditional closing probability: for $\alpha > 0$, we write
$$
 \hPhi_{n,m}^\alpha = \Big\{ \gamma \in \fpart_{n,m}: q_{n,m}(\gamma) > m^{-\alpha} \Big\} \, .
$$

\section{The snake method: general elements}\label{s.five}

In this section, we present in a general form a proof-by-contradiction technique, which we call the {\em snake method}, that was already employed in~\cite{ontheprob}. We will later use it in two different ways, to prove Theorem~\ref{t.closingprob} and Theorem~\ref{t.thetexist}.

The snake method is used to prove upper bounds on the closing probability, and assumes to the contrary that to some degree this probability has slow decay. For the technique to be used, two ingredients are needed.
\begin{enumerate}
\item A charming snake is a walk or polygon $\ga$ many of whose subpaths beginning at $\northeast(\ga)$  have high conditional closing probability, when extended by some common length. It must be shown that charming snakes are not too atypical.
\item A general procedure is then specified in which 
a
charming snake is used to manufacture huge numbers of alternative self-avoiding walks. These alternatives overwhelm the polygons in number and show that the closing probability is very small, contradicting the assumption.
\end{enumerate}

The first step certainly makes use of the assumption of slow decay on the closing probability. It is not however a simple consequence of this assumption. We will carry out the first step in two very different ways in this article: via Gaussian pattern fluctuation to prove Theorem~\ref{t.closingprob} later in this part, and via polygon joining to prove Theorem~\ref{t.thetexist} in Part~\ref{c.aboveonehalf}.

In contrast to the different approaches used to implement the first step, the second step is performed using a general tool, valid in any dimension $d \geq 2$, that we present in this section. Indeed, the snake method including this second step tool was used in~\cite{ontheprob} to show that the closing probability is bounded above by~$n^{-1/4 + o(1)}$. The second step draws inspiration from the notion  that reflected walks offer alternatives to closing (or near closing) ones that appears in Madras' derivation~\cite{Madras14} of lower bounds on moments of the endpoint distance under~$\psaw_n$.  

\subsection{The general apparatus of the method}

\subsubsection{Parameters}

The snake method has three exponent parameters:
\begin{itemize}
\item the inverse charm $\alpha > 0$;
\item the snake length $\beta \in (0,1]$;
\item and the charm deficit $\esm \in (0,\beta)$.
\end{itemize}
It has two index parameters:
\begin{itemize}
\item $n \in 2\N + 1$ and $\ell \in \N$, with $\ell \leq n$.
\end{itemize}

\subsubsection{Charming snakes}

Here we define these creatures.

\begin{definition}
Let $\alpha > 0$, $n \in 2\N + 1$,  $\ell \in [0,n]$,
$\gamma \in \fpart_{\ell,n}$, and $k \in [0,\ell]$ with $\ell - k \in 2\N$.  We say that $\ga$ is $(\alpha,n,\ell)$-{\em charming} at (or for) index $k$ if 
\begin{equation}\label{e.charming}
\psaw^0_{k + n - \ell} \Big( \Gamma \closes \, \Big\vert \, \vert \Gamma^1 \vert = k \, , \, \Gamma^1 = \gamma_{[0,k]} \Big) > n^{-\alpha} \, .
\end{equation}
\end{definition}
(The event  that $\vert \Gamma^1 \vert = k$ in the conditioning is redundant and is recorded for emphasis.)
Note that an element $\ga \in \fpart_{\ell,n}$  is $(\alpha,n,\ell)$-charming at  index~$k$ if in selecting a length~$n-\ell$ walk beginning at $0$ uniformly among those choices that form the second part of a walk whose first part is $\ga_{[0,k]}$, the resulting $(k + n - \ell)$-length walk closes with probability exceeding $n^{-\alpha}$. 
(Since we insist that $n$ is odd and that $\ell$ and $k$ share their parity, the length $k + n - \ell$ is odd; the condition displayed above could not possibly be satisfied if this were not the case.)
Note that, for any $\ell \in [0,n]$,  $\ga \in \fpart_{\ell,n}$ is $(\alpha,n,\ell)$-charming at the special choice of index $k = \ell$ precisely when $\ga \in \hPhi_{\ell,n}^{\alpha}$. When $k < \ell$ with $k+n - \ell$ of order~$n$, the condition that 
 $\ga \in \fpart_{\ell,n}$ is $(\alpha,n,\ell)$-charming at  index $k$ is {\em almost} the same as  $\ga_{[0,k]} \in \hPhi_{k,k + n - \ell}^{\alpha}$; (the latter condition would involve replacing $n$ by $k + n -\ell$ in the right-hand side of~(\ref{e.charming})). 

For $n \in 2\N + 1$, $\ell \in [0,n]$, $\alpha,\beta > 0$ and $\esm \in (0,\beta)$, define the charming snake set 
\begin{eqnarray*}
\mathsf{CS}_{\beta,\esm}^{\alpha,\ell,n}  & =  & \Big\{ \ga \in \fpart_{\ell,n}: \textrm{$\ga$ is $(\alpha,n,\ell)$-charming} \\
 & & \quad \textrm{for at least $n^{\beta- \esm}/4$ indices belonging to the interval $\big[ \ell - n^{\beta}, \ell \big]$} \Big\}
 \, .
\end{eqnarray*}
For any element of $\ga \in \fpart_{\ell,n}$, think of an extending snake consisting of $n^{\beta} + 1$ terms
$\big( \gamma_{[0,\ell - n^{\beta}]},\ga_{[0,\ell - n^{\beta} + 1]}, \cdots , \ga_{[0,\ell]} \big)$. If $\ga \in 
\mathsf{CS}_{\beta,\eps}^{\alpha,\ell,n}$, then there are many charming terms in this snake, for each of which there is a high conditional closing probability for extensions of a {\em shared} continuation length~$n - \ell$.

\subsubsection{A general tool for the method's second step}

For the snake method to produce results, we must work with a choice of parameters for which $\beta - \esm - \alpha > 0$. (The method could be reworked to handle equality in some cases.) Here we present the general tool for carrying out the method's second step. This technique was already presented in \cite[Lemma 5.8]{ontheprob}, and our treatment differs only by 
using notation adapted for the snake method with general parameters. 

The tool asserts that, if $\beta - \esm - \alpha > 0$ and even a very modest proportion of snakes are charming, then the closing probability drops precipitously.

\begin{theorem}\label{t.snakesecondstep}
Let $d \geq 2$. Set $c = 2^{\tfrac{1}{5(4d + 1)}} > 1$ 
and set $K = 20 (4d+1) \tfrac{\log(4d)}{\log 2}$.
Suppose that the exponent parameters satisfy $\delta = \beta - \esm - \alpha > 0$.
If the index parameter pair $(n,\ell)$ satisfies $n \geq K^{1/\delta}$ and 
\begin{equation}\label{e.afewcs}
 \psap_{n+1} \Big( \Gamma_{[0,\ell]} \in \mathsf{CS}_{\beta,\esm}^{\alpha,\ell,n}  \Big)  \geq c^{- n^{\delta}/2} \, ,
\end{equation}
then
$$
\psaw_n \Big( \Gamma \closes \Big) \leq 2 (n+1) \,  c^{-n^{\delta}/2} \, .
$$ 
\end{theorem}

Note that since the closing probability is predicted to have polynomial decay, the hypothesis~(\ref{e.afewcs}) is never satisfied in practice. For this reason, the snake method will always involve argument by contradiction, with~(\ref{e.afewcs}) being established under a hypothesis that the closing probability is, to some degree, decaying slowly.

\subsection{A charming snake creates huge numbers of reflected walks}
Here is the principal component of the proof of Theorem~\ref{t.snakesecondstep}. 
\begin{proposition}\label{p.avoidance}
Let $d \geq 2$. Set $\delta = \beta - \esm - \alpha$ and suppose that $\delta > 0$. 
Let $\phi \in \mathsf{CS}_{\beta,\esm}^{\alpha,\ell,n}$. With $c  > 1$ and $K > 0$ specified in Theorem~\ref{t.snakesecondstep}, 
we have that, if $n \geq K^{1/\delta}$, then 
$$ 
\# \Big\{ \gamma \in \sawfree_n:  \gamma_{[0,\ell]} = 
 \reverse\phi \Big\} \geq c^{n^{\delta}} \cdot
 \# \Big\{ \gamma \in \sawfree_n:  \northeast(\gamma) = 0 ,  \gamma^1 = \phi \Big\} \, .
$$
\end{proposition}
Note here that walks beginning with the reversal of an element $\phi \in \saw_\ell$ will necessarily not begin at the origin, and thus we employ the notation introduced in Subsection~\ref{s.sawfree}.  

\medskip

\noindent{\bf Proof of Proposition~\ref{p.avoidance}.}
 Let $ \Wtmp$ denote the set of walks $\ga$ of length~$n-\ell$
that originate at $0$ and for which $\northeast(\ga) = 0$. This set is not the same as $\fpart_{n-\ell}$ when $d=2$ (though it is when $d \geq 3$), because we do not stipulate that $\ga_1 = -e_1$; in fact, whatever the value of $d \geq 2$, we will be using that $\Wtmp$ contains all possible length~$n-\ell$ walks that form the second (rather than the first) part of the two-part decomposition of some walk of at least this length.

  Let $\Ptmp$ denote the uniform measure on the set $\Wtmp$. We will denote by
  $\Gamma$ a random variable distributed according to $\Ptmp$. 
  In particular, $\Gamma$ is contained in the lower half-space including the origin. (When $d=2$, we mean the region on or below the $x$-axis, and, when $d \geq 3$, the region of non-positive $e_1$-coordinate.)

  We now extend the notion of closing walk by saying that  $\ga'$ {\em closes} $\ga$ if $\ga_0 = \ga'_0$ 
  and the endpoints of  $\ga$ and $\ga'$ are adjacent. 
  We say that $\ga'$ {\em avoids} $\ga$ if no vertex except $\ga_0$ is visited by both $\ga'$ and $\ga$.

We are given $\phi \in \fpart_{\ell,n}$ such that $\phi \in \mathsf{CS}_{\beta,\esm}^{\alpha,\ell,n}$. By definition, we may find indices $j_1 < j_2 < \ldots < j_{\lceil n^{\beta - \esm}/4 \rceil}$ lying in $\big[ \ell - n^{\beta},\ell \big]$ at each of which $\phi$ is $(\alpha,n,\ell)$-charming.

  For $1 \leq i \leq \lceil n^{\beta - \esm}/4 \rceil$, define the events
  \begin{align*}
    A_i & = \Big\{ \Gamma \textrm{ avoids } \, \phi_{[0,j_i]} \Big\} \quad \text{ and } \quad C_i=\Big\{ \Gamma \textrm{ closes } \, \phi_{[0,j_i]}\Big\} \, .
  \end{align*}
  Also, define the set $\Atmp = \big\{ \ga \in \Wtmp : \ga \textrm{ avoids } \phi_{[0,\ell]} \big\}$.
   
  Since $\phi$ is $(\alpha,n,\ell)$-charming at index $j_i$, 
  \begin{align}\label{e.closing.avoiding}
    \Ptmp \Big( \Gamma \closes \, \, \phi_{[0,j_i]} \, \Big\vert \,  \Gamma \text{ avoids } \phi_{[0,j_i]} \Big)
    = \Ptmp \big( C_i \, \big\vert \,  A_i \big) 
    > n^{ - \alpha} \, .
  \end{align}
  Write $k = \lceil  4d \, n^{\alpha} \rceil$. (Note that $k \leq n^{\beta - \esm}/4$ holds for $n$ high enough since $\delta$ is supposed positive.) 
  Any realization $\Ga \in \Wtmp$ is in at most $2d$ events $C_i$.
  Hence, by \eqref{e.closing.avoiding} and the~$A_j$ being decreasing, 
  \begin{align*}
    2d 
    \geq \sum_{i=1}^{k} \Ptmp (C_i)
    \geq \sum_{i=1}^{k} \Ptmp \big(C_i \cond A_i\big) \cdot \Ptmp(A_k)
    \geq 4d \, \Ptmp(A_k) \, .
  \end{align*}
  Therefore, $\Ptmp (A_k) \leq \frac{1}{2}$.
  If the procedure is repeated for indices between $k + 1$ and $2k$, one obtains
  \begin{align*}
    2d
    \geq \sum_{i= k+1}^{2k} \Ptmp (C_i \cond A_k)
    \geq \sum_{i=k+1}^{2k} \Ptmp \big( C_i \cond A_i \big) \cdot \Ptmp \big(A_{2k} \cond A_k \big)
    \geq 4d \, \Ptmp \big( A_{2k} \cond A_k \big) \, ,
  \end{align*}
  and thus $P(A_{2k} \cond A_k ) \leq 1/2$.
  Since $A_{2k}\subset A_k$, we find 
  $$
  \Ptmp (A_{2k}) = \Ptmp (A_k)\Ptmp (A_{2k} \cond A_k) \leq \tfrac14 \, .
  $$
  In these inequalities, we see the powerful bootstrap mechanism at the heart of the snake method, demonstrating that $P(A_{(i+1)k})$ is at most one-half of~$P(A_{ik})$. The mechanism works because
  the method's definitions imply that all walk extensions are of common length~$n-\ell$, and the avoidance conditions are monotone (i.e., the events $A_i$ are decreasing).
  
Indeed,  the procedure may be repeated $ \lfloor \frac{n^{\beta - \esm}}{4k} \rfloor \geq \tfrac{n^{\beta - \esm - \alpha}}{4(4d+1)} - 1$ times. 
Recalling that $\phi = \phi_{[0,\ell]}$, we obtain
  \begin{align*}
    \frac{|\Atmp|}{|\Wtmp|}=\Ptmp\Big(\Gamma\text{ avoids } \phi \Big) \le \Ptmp \big(A_{\lceil n^{\beta - \esm}/4 \rceil} \big) 
    \leq 2^{- \lfloor \tfrac{n^{\beta - \esm}}{4k}\rfloor} 
    \leq 2^{1 - \tfrac{n^{\beta - \esm -  \alpha}}{4(4d+1)}} \, ,
  \end{align*}
whatever the value of $n \in 2\N + 1$.

  \begin{figure}
    \begin{center}
      \includegraphics[width=0.75\textwidth]{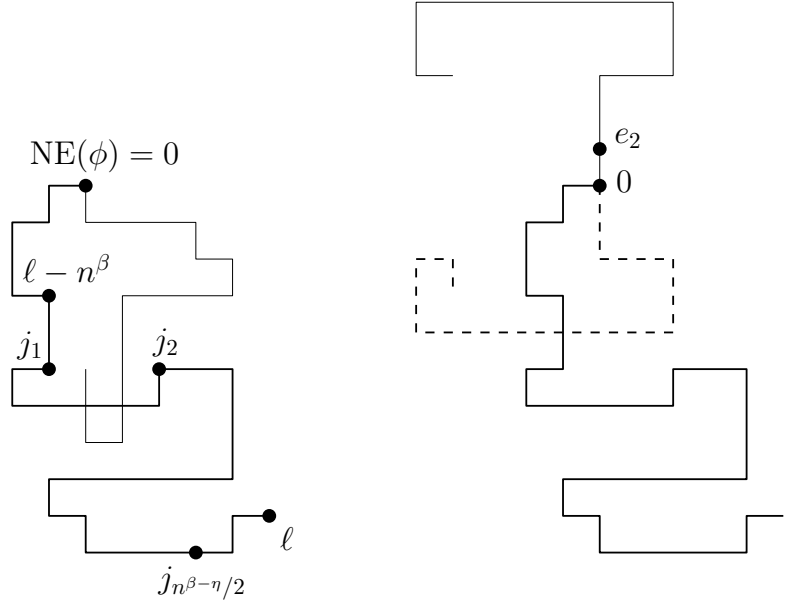}
    \end{center}
    \caption{On the left, $\phi \in \fpart_{\ell,n}$ is depicted in bold. The labels along $\phi$ are indices. A walk in $W \cap A_1 \cap C_1 \cap A_2^c$ is also shown. On the right, another element $\ga \in W$ is depicted with dashed lines. The output of the three-part concatenation procedure used at the end of the proof of Proposition~\ref{p.avoidance} is shown above~$0$.}\label{f.avoided}
  \end{figure}

  The set $\Atmp$ contains all length $n - \ell$ walks $\ga$  
  for which $\big[ \phi , \ga \big]$ is the two-part decomposition of  a walk of length $n$ with $\northeast = 0$. Thus,
$$
 \big\vert \Atmp \big\vert = \# \Big\{ \gamma \in \sawfree_n:  \northeast(\gamma) = 0 \, , \, \vert \gamma^1 \vert = \ell \, , \, \gamma^1 = \phi \Big\} \, .
$$

  On the other hand, for $\ga \in \Wtmp$, 
  consider the walk obtained by concatenating three paths (and illustrated in Figure~\ref{f.avoided}). When $d=2$, these are: the reversal $\reverse{\phi}$ of $\phi$; 
  the edge $e_2$; and 
the $e_2$-translation of the reflection of $\ga$ in the horizontal axis. When $d \geq 3$, we substitute $e_1$ for $e_2$, and the zero-coordinate hyperplane in the $e_1$-direction for the horizontal axis, to specify these paths.
The walk that results has length $n + 1$ and is self-avoiding. 
  By deleting the last edge of such walks, we obtain at least $|\Wtmp| / 2d$
  walks of length $n$, each of which follows $\reverse\phi$ in its first $\ell$ steps.
  Thus, 
$$
\# \Big\{ \gamma \in \sawfree_n:  \gamma_{[0,\ell]} = \reverse\phi \Big\} \geq |\Wtmp| / 2d \, .
$$
The three preceding displayed equations combine 
prove that the ratio of the cardinalities on the left and right-hand sides in Proposition~\ref{p.avoidance}
is at least $(2d)^{-1} 2^{\tfrac{n^{\beta - \esm -  \alpha}}{4(4d+1)} - 1}$.
Noting the lower bound on $n$ stated in Proposition~\ref{p.avoidance} completes the proof of this result. \qed

\medskip

\noindent{\bf Proof of Theorem~\ref{t.snakesecondstep}.}
Write $\mathsf{CS}=  \mathsf{CS}_{\beta,\esm}^{\alpha,\ell,n}$.
Now using Proposition~\ref{p.avoidance} in the second and~(\ref{e.afewcs}) in the fourth inequalities,
\begin{eqnarray*}
 \big\vert \saw_n \big\vert & \geq &   \sum_{\phi \in \mathsf{CS}} \# \Big\{ \gamma \in \saw_n: \gamma_{[0,\ell]} = -\phi_{\ell} + \reverse\phi \Big\} \\
& = & \sum_{\phi \in \mathsf{CS}} \# \Big\{ \gamma \in \sawfree_n: \gamma_{[0,\ell]} = \reverse\phi \Big\} \\
 & \geq &  c^{n^\delta}  \cdot \# \Big\{ \gamma \in \sawfree_n:  \northeast(\gamma) = 0 \, , \, \gamma^1 \in \mathsf{CS} \Big\} \\
 & \geq &  c^{n^\delta}  \cdot \# \Big\{ \gamma \in \sawfree_n:  \northeast(\gamma) = 0 \, , \, \gamma^1 \in \mathsf{CS}  \, , \, \gamma \closes \Big\} \\
 & = &  c^{n^\delta}  \cdot \# \Big\{ \gamma \in \sap_{n+1}:   \gamma_{[0,\ell]} \in \mathsf{CS} \Big\} \geq  c^{n^{\delta}/2}  \cdot    \big\vert  \sap_{n+1} \big\vert \, .  \\
\end{eqnarray*}
Thus, $p_{n+1}/c_n \leq c^{-n^{\delta}/2}$. We find that $$\psaw_n \big( \Gamma \closes \big) = 2 (n+1) p_{n+1}/c_n \leq 2 (n+1) c^{-n^{\delta}/2} \, . $$ \qed  

\section{The snake method applied via Gaussian pattern fluctuation}\label{s.six}

In \cite[Theorem 1.1]{ontheprob}, it is proved that $\psaw_n \big( \Ga \closes \big) \leq n^{-1/4 + o(1)}$. The technique used, which is the snake method with the first step carried out using Gaussian pattern fluctuation, clearly cannot prove faster decay on the closing probability than $n^{-1/2}$. We now rework the method to prove the assertion Theorem~\ref{t.closingprob} that this faster decay rate can be obtained for any dimension $d \geq 2$. 

This application of the snake method is probabilistic, with the closing probability being discussed without reference to the polygon and walk numbers and the relating equation~(\ref{e.closepc}). As such, the argument does not fit in the framework of $(\eone,\etwo,\ethree)$-arguments introduced in Subsection~\ref{s.sumthree}.


As we have stated, we will prove Theorem~\ref{t.closingprob} by assuming that its conclusion fails and seeking a contradiction. By a relabelling of $\eps > 0$, we may express the premise that the conclusion fails in the form that, for some $\eps > 0$ and infinitely many $n \in 2\N + 1$,
\begin{equation}\label{e.neass}
\psaw_n\big( \Gamma \closes \big) \geq n^{-1/2 + 4\eps} \, .
\end{equation}


We fix the three snake method exponent parameters. Fixing a given choice of $\eps \in (0,1/4)$, we set them equal to  $\alpha = 1/2 - 2\eps$, $\beta = 1/2$ and $\esm  = 0$.
 We will argue that the hypothesis~(\ref{e.afewcs}) is comfortably satisfied, with charming snakes being the~norm.
\begin{proposition}\label{p.snakeg}
There exists a positive constant $C$ such that, if $n \in 2\N + 1$ and $\eps \in (0,1/4)$ satisfy (\ref{e.neass}) as well as $n \geq \max \{ 4^{1/\eps} , C \}$, then there exists  $\ell \in \big[ n/4 , 3n/4   \big]$ for which
$$
 \psap_{n+1} \Big( \Gamma_{[0,\ell]} \in \mathsf{CS}_{1/2,0}^{1/2 - 2\eps,\ell,n} \Big) \geq 1 \, - \, 
 n^{-\eps/7}  \, . 
$$
\end{proposition}

\noindent{\bf Proof of Theorem~\ref{t.closingprob}.} 
Given our choice of the three exponent parameters, note that $\delta = \beta - \esm - \alpha$ equals $2\eps$ and is indeed positive. The conclusion of Theorem~\ref{t.snakesecondstep} contradicts~(\ref{e.neass}) if $n$ is high enough. Thus~(\ref{e.neass}) is false for all $n$ sufficiently high. Since $\eps \in (0,1/4)$ may be chosen arbitrarily small, we are done. \qed

\medskip

It remains only to prove Proposition~\ref{p.snakeg}, and the rest of the section is devoted to this proof.
The plan in outline has two steps. In the first, implemented in Lemma~\ref{l.polyclose},
we will infer from the closing probability lower bound~(\ref{e.neass}) a deduction that for typical indices $\ell$ close to $n/2$, the initial (or first part) segment $\Gamma_{[0,\ell]}$
of a $\psap_n$-distributed polygon typically has a not-unusually-high conditional closing probability for a second part extension of length $n - \ell$. This step is a straightforward, Fubini-style statement. In the second step, we work with this deduction to build a charming snake. For this, we want an inference of the same type, but with the second part length remaining fixed, even as the first part length varies over an interval of length close to $n^{1/2}$. Beginning in Subsection~\ref{sec:pattern}, the mechanism of Gaussian fluctuation in local patterns along the polygon is used for this second step. Crudely put, typical deviations in the number of local configurations (type $I$ and type $II$ patterns) in the initial and final ten-percent-length segments of the polygon create an ambiguity in the location marked by the index $\ell \approx n/2$ of order~$n^{1/2}$.
This square-root fuzziness yields charming snakes.  

\subsection{Setting the snake method index parameters}\label{s.indexset}
We now specify the values of $n$ and $\ell$.
The value of $n \in 2\N + 1$ is supposed to satisfy~(\ref{e.neass}) for our given $\eps \in (0,1/4)$,  as well as the bound $n \geq 4^{1/\eps}$.
Applying Lemma~\ref{l.closecard} with $\alphamac = 1/2 - 4\eps$ and $\deltamac = \eps$, and noting that this lower bound on $n$ ensures that
$2 n^{1-\eps} < \# \big[ n/4 , 3n/4  \big]$, we find that we may select $\ell$ to lie in $\big[ n/4 , 3n/4   \big]$ and to satisfy 
$$
 \# \Big\{ \gamma \in \saw^0_n: \vert \gamma^1 \vert = \ell \Big\}
 < n^{1/2 - 3\eps}  \cdot \# \Big\{ \gamma \in \saw^0_n: \vert \gamma^1 \vert = \ell  , \gamma \closes \Big\} \, ,
$$
or equivalently
$$
\psaw^0_n \Big(  \Gamma \closes \, \Big\vert \, \vert \Gamma^1 \vert = \ell   \Big) > n^{-1/2 + 3\eps} \, .
$$
The value of $\ell$ is so fixed henceforth in the proof of Proposition~\ref{p.snakeg}.

\begin{lemma}\label{l.polyclose}
$$
 \psap_{n+1} \Big( \Gamma_{[0,\ell]} \not\in \hPhi_{\ell,n}^{1/2 - 2\eps} \Big) \leq n^{-\eps} \, .
$$
\end{lemma}
\noindent{\bf Proof.} 
Note that $\Gamma_{[0,\ell]}$ under $\psap_{n+1}$ shares its law with the first part $\Gamma^1$ under $\psaw^0_n \big( \cdot \, \big\vert \, \Gamma \closes \, , \, \vert \Gamma^1 \vert = \ell \big)$. For this reason, the statement may be reformulated
\begin{equation}\label{e.reform}
 \psaw^0_n \Big(  \Gamma^1 \not\in \hPhi_{\ell,n}^{1/2 - 2\eps} \, \Big\vert \, \vert \Gamma^1 \vert = \ell \, , \, \Gamma \closes \Big) \leq n^{-\eps} \, .
\end{equation}
To derive (\ref{e.reform}), set $p$ equal to its left-hand side. Note that
\begin{eqnarray*}
  & & \# \Big\{ \gamma \in \saw^0_n: \vert \gamma^1 \vert = \ell   \Big\} \\
 & \geq &   \# \Big\{ \gamma \in \saw^0_n: \vert \gamma^1 \vert = \ell \, , \, \gamma^1 \not\in \hPhi_{\ell,n}^{1/2 - 2\eps}   \Big\} \\
 & = & \sum_{\phi \in \fpart_{\ell,n} \setminus \hPhi_{\ell,n}^{1/2 -  2\eps}}  \# \Big\{ \gamma \in \saw^0_n:  \gamma^1  = \phi   \Big\} \\
 & \geq & \sum_{\phi \in \fpart_{\ell,n} \setminus \hPhi_{\ell,n}^{1/2 - 2\eps}} n^{1/2 - 2\eps} \cdot  \# \Big\{ \gamma \in \saw^0_n:  \gamma^1  = \phi \, , \, \gamma \closes  \Big\} \\
 & = &   n^{1/2 -  2\eps} \cdot  \# \Big\{ \gamma \in \saw^0_n: \vert \gamma^1 \vert = \ell \, , \, \gamma^1 \not\in \hPhi_{\ell,n}^{1/2 - 2\eps} \, , \, \ga \closes  \Big\} \\ 
& = &   n^{1/2 - 2\eps} \cdot  p \cdot \# \Big\{ \gamma \in \saw^0_n: \vert \gamma^1 \vert = \ell  \, , \, \ga \closes  \Big\} \, ,
\end{eqnarray*}
whence 
$$
 n^{1/2 - 2\eps} \cdot  p \leq \psaw^0_n \Big( \Ga \closes \, \Big\vert \, \vert \Gamma^1 \vert = \ell  \Big)^{-1} \, , 
$$
whose right-hand side we know to be at most $n^{1/2 - 3\eps}$. Thus, $p \leq n^{-\eps}$, and we have verified~(\ref{e.reform}). \qed

\medskip

\noindent{\em Remark.}
Lemma~\ref{l.polyclose} may be compared to \cite[Lemma 5.5]{ontheprob} with $k=0$.
The former result states that $\hPhi_{\ell,n}^{1/2 - 2\eps}$ membership by~$\Gamma_{[0,\ell]}$ is the norm under~$\psap_{n+1}$, while the latter merely asserts that a comparable membership is not unlikely (having probability at least $n^{-1/4 + 2\eps}$).
It may be possible to improve the inequality in~\cite[(5.7)]{ontheprob}
to reflect the fact that the conditioning on closing leading from the law $\psaw_n$
to $\psap_{n+1}$ reweights the measure on first parts proportionally in accordance
with the conditional closing probability given the first part. Analysing this reweighting may lead to a replacement of the right-hand side of \cite[(5.5)]{ontheprob} by a term of the form $1 - n^{-o(1)}$ and, alongside other suitable changes, permit a derivation of Theorem~\ref{t.closingprob}.

\subsection{Patterns and shells}\label{sec:pattern}

Patterns are local configurations in self-avoiding walks that are the subject of a famous theorem~\cite{kestenone} due to Kesten that we will shortly state. For our present purpose, we identify two particular patterns.

\begin{figure}
  \begin{center}
    \includegraphics[width=0.4\textwidth]{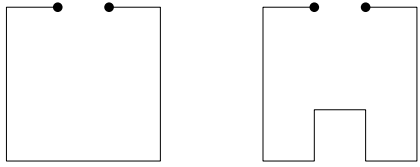}
  \end{center}
  \caption{An example of type I and II patterns for $d =2$.}
  \label{f.patternstwo}
\end{figure}

\begin{definition}[Type I/II patterns] 
  A pair of type I and II patterns is a pair of self-avoiding walks $\chi^I$, $\chi^{II}$,  
  both contained in the cube $[0,3]^d$, with the properties that
  \begin{itemize*}
  \item $\chi^I$ and $\chi^{II}$ both visit all vertices of the boundary of $[0,3]^d$,
  \item $\chi^I$ and $\chi^{II}$ both start at   $\big( 1,3,1, \cdots , 1 \big)$ and end at $\big( 2,3, 1, \cdots, 1 \big)$,
  \item the length of $\chi^{II}$ exceeds that of $\chi^{I}$ by two.
  \end{itemize*}
\end{definition}
Figure \ref{f.patternstwo} depicts examples of such patterns for $d =2$. 
The existence of such pairs of walks for any dimension $d \geq 2$ may be easily checked, and no details are given here. 
Fix a pair of type $I$ and $II$ patterns henceforth. 

A pattern $\chi$ is said to occur at step $k$ of a walk $\ga$ 
if $\ga_{[k, k + |\chi|]}$ is a translate of $\chi$ (where recall that $\vert \chi \vert$ is the length of~$\chi$).
A {\em slot} of $\ga$ is any translate of $[0,3]^d$ containing $\ga_{[k, k + |\chi|]}$ where a pattern~$\chi$ of type $I$ or $II$ occurs at step $k$ of $\gamma$.
Note that the slots of $\ga$ are pairwise disjoint. 

In~\cite{ontheprob}, the notion of shell was introduced. A shell is an equivalence class of self-avoiding walks under the relation that two walks are identified if one may be obtained from the other by changing some patterns of type~$I$ to be of type~$II$ and {\em vice versa}. The walks in a given shell share a common set of slots, but they are of varying lengths. The shell of a given walk $\ga$ is denoted~$\vas(\ga)$.

Consider two walks $\ga^1$ and $\ga^2$ with $\ga^1_0 = \ga^2_0$. Recall that $\ga^2$ {\em avoids} $\ga^1$ if  no other vertex is visited by both walks. 
The next fact is crucial to our reasons for considering shells; its almost trivial proof is omitted.
\begin{lemma}\label{l.sharedavoidance}
For some $m \in \N$, let $\ga \in \saw_m$ and let  $\ga' \in \vas(\ga)$.
A walk beginning at $\ga_0$ avoids $\ga$ if and only if it avoids $\ga'$. 
\end{lemma}

\begin{figure}
  \begin{center}
    \includegraphics[width=1.0\textwidth]{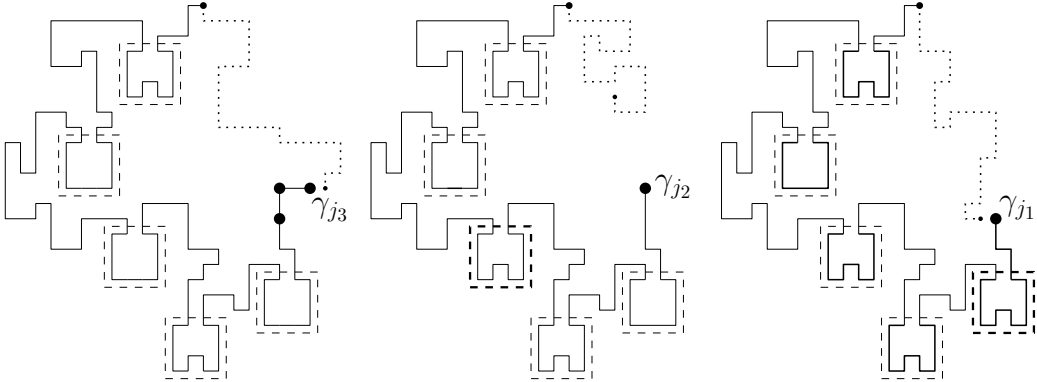}
  \end{center}
  \caption{In this figure, we explain in outline the method. Given $n \in 2\N + 1$, the index $\ell \in [n/4,3n/4]$ has been fixed so that the bound~(\ref{e.reform}) holds. As we have seen, the vast majority of indices in this range satisfy this bound. This means that when we draw a  length~$n+1$ polygon $\ga$ and mark with a black dot each vertex $\ga_j$, $n/4 \leq j \leq 3n/4$, with the property that 
  $\psaw^0_n \big( \Ga \closes \, \big\vert \, \vert \Ga^1 \vert = j , \Ga^1 = \ga_{[0,j]} \big) \geq n^{-1/2 + 2\eps}$,  
   most such~$\ga$ will appear with black dots in most of the available spots.  The left-hand sketch represents such a typical~$\ga$ and three of its many black dots. The pointed second part shows a sample of the law $\psaw^0_n \big( \cdot \, \big\vert \, \vert \Ga^1 \vert = j_k , \Ga^1 = \ga_{[0,j_k]} \big)$, with $k=3$, one that happens to close~$\ga_{[0,j_3]}$. The second part being sampled has length $n - j_3$. If we sample instead this law with $k=2$, then the second part has a greater length, $n - j_2$, which equals $n-j_3 + 2$ in this instance. However, to construct a charming snake, we want this length to stay the same as we move from one black dot to the next. Pattern exchange is the mechanism that achieves this. By turning one type~$I$ pattern into a pattern of type~$II$,  we push two units of length into the first part, so that, in the middle sketch, the random second part has the original length $n - j_3$. The process is iterated in the right-hand sketch. The first part is akin to a belayer who takes in rope, storing it in accumulating type~$II$ patterns, so that the second part climber maintains a constant length of rope. This process of pattern exchange can be maintained for an order of $n^{1/2}$ steps, because the Gaussian fluctuation between the two types of pattern means that the process of artificially altering pattern type does not push the system out of its rough equilibrium when the number of changes is of this order. In this way, black dots also mark charming snake terms for a snake of length of order~$n^{1/2}$.}\label{f.snakemethodpattern}
\end{figure}

The reader may now wish to view Figure~\ref{f.snakemethodpattern} and its caption for an expository overview of the snake method via Gaussian pattern fluctuation. We mention also that this Gaussian fluctuation has been utilized in~\cite{JROSTW} to prove a $n^{1/2 - o(1)}$ lower bound on the absolute value of the writhe of a typical length~$n$ polygon.

We will make some use of the notion of shell, but will predominantly consider 
a slightly different definition, the $(n+1)$-local shell, which we now develop.
This new notion concerns polygons rather than walks. 

Recall that the parameter $n \in 2\N +1$ has been fixed in Subsection~\ref{s.indexset}. (The upcoming  definitions do not require that $n \in 2\N + 1$ be fixed in this particular way in order to make sense, but, when we make use of the definitions, it will be for this choice of $n$.)
Define an equivalence relation $\sim$ on $\sap_{n+1}$ as follows. For any $\gamma \in \sap_{n+1}$, let $\gaemp$ denote the polygon in the shell of $\gamma$ that has no type $II$ patterns, (formed by switching every type $II$ pattern of $\ga$ into a type $I$ pattern). Thus, $\gaemp \in \sap_{n + 1 - 2T_{II}(\ga)}$, where $T_{II}(\ga)$ denotes the total number of type~$II$ patterns in $\ga$. A type $II$ pattern contains thirteen edges (in $d=2$; at least this number in higher dimensions), and these patterns are disjoint, so $T_{II}(\ga) \leq (n+1)/13$. Thus, the length of $\gaemp$ is at least $11(n+1)/13$. Let $S_1$ denote the set of slots in $\gamma$ that are slots in $\gaemp_{[0,(n+1)/10]}$, and, writing $\ellempty$ for the length of $\gaemp$, let $S_2$ denote the set of slots in $\gamma$ that are slots in $\gaemp_{[\ellempty - (n+1)/10,\ellempty]}$. Note that $S_1$ and $S_2$ are disjoint. We further write $N_I(\ga)$ and $N_{II}(\ga)$ for the number of patterns of the given type in the slots $S_1 \cup S_2$, and $N_{I}^1(\ga)$ and $N_{I}^2(\ga)$ for the number of type~$I$ patterns occupying slots in $S_1$ and in $S_2$; and similarly for $N_{II}^1(\ga)$ and $N_{II}^2(\ga)$. 

For $\ga,\ga' \in \sap_{n+1}$, we say that $\ga \sim \ga'$ 
if $\gamma'$ may be obtained from $\gamma$ by relocating the type $II$ patterns of $\gamma$ contained in the set of slots $S_1 \cup S_2$ for $\gamma$ to another set of locations among these slots. The relation $\sim$ is an equivalence relation, because the polygon $\gaemp$ formed by filling all the slots of $\ga$ with type $I$ patterns is shared by related polygons, so that the value of $S_1 \cup S_2$ is equal for such polygons. 
Elements of a given equivalence class have common values of length, $N_I$ and $N_{II}$, but not of $N_I^i$ or $N_{II}^i$ for $i \in \{1,2\}$. 
We call the equivalence classes $(n+1)$-local shells: the parameter $n+1$ appears to denote the common length of the member polygons, and the term {\em local} is included to indicate that members of a given class may differ only in locations that are close to the origin (in the chemical distance, along the polygon). Complementing the notation $\varsigma(\ga)$ for the shell of $\ga$,
 write $\vaso(\ga)$ for the $(n+1)$-local shell of $\ga \in \sap_{n+1}$.

For $\deltanew > 0$, write $\mathcal{G}_{n+1,\deltanew}$ for the set of $(n+1)$-local shells $\sigma \subseteq  \sap_{n+1}$ such that each of the quantities $\vert S_1(\sigma) \vert$, $\vert S_2(\sigma) \vert$, $N_I$ and $N_{II}$ is at least $\deltanew (n+1)$. Such ``good'' shells are highly typical if $\deltanew > 0$ is small, as we now see.

\begin{lemma}\label{l.kptconseq}
There exist constants $c > 0$ and $\deltanew > 0$ such that 
$$
\psap_{n+1} \Big( \vaso(\Gamma) \in \mathcal{G}_{n+1,\deltanew} \Big) \, \geq \, 1 - e^{-cn} \, .
$$ 
\end{lemma}
\noindent{\bf Proof.}
By Kesten's pattern theorem \cite[Theorem 1]{kestenone} (and Lemma~\ref{l.closeproblb}, in order that we may state the result for the polygon measure), there exist constants $c > 0$ and $\deltanew >0$ such that,
for any odd $n \geq d 3^d$,
\begin{equation}\label{e.patterns}
  \psap_{n+1} \Big( T_I(\Ga) \leq \deltanew m \Big) \le e^{-c n} 
  \quad \text{ and } \quad
  \psap_{n+1} \Big( T_{II}(\Ga) \leq \deltanew m \Big) \le e^{-c n} \, .
\end{equation}

Note that every slot in $S_1$ belongs to $\Gamma_{[0,(n+1)/4]}$, because such slots belong to  $\Gamma_{[0,(n+1)/10 + 2 T_{II}(\Gamma)]}$ and $T_{II}(\Gamma) \leq (n+1)/13$.
Thus, 
\begin{eqnarray*}
 & & \psap_{n+1} \Big( S_1(\Ga) \text{ contains fewer than $\deltanew (n+1)$ type I patterns} \Big) \\
  & \leq & \psap_{n+1} \Big( \Ga_{[0, (n+1)/4]} \text{ contains fewer than $\deltanew (n+1)$ type I patterns} \Big) \,  .
\end{eqnarray*}

  There exist constants $\deltanew ,c > 0$ such that 
  \begin{align}
    &\psap_{n + 1} \Big( \Ga_{[0,(n + 1)/4]} \text{ contains fewer than $\deltanew (n + 1)$ type I patterns} \Big) \nonumber \\
    &\le  \frac{c_{(n+1)/4} c_{3(n+1)/4}}{p_{n+1}} \cdot \psaw_{(n+1)/4} \Big( \Ga\text{ contains fewer than $\deltanew (n+1)$ type I patterns} \Big) \nonumber \\
    &\le e^{(2c_{{\rm HW}} + c) \sqrt n} \cdot \psaw_{(n+1)/4} \Big( \Ga\text{ contains fewer than $\deltanew (n+1)$ type I patterns} \Big) \nonumber
    < e^{-c n},
  \end{align}
where the second inequality comes from the Hammersley-Welsh bound (\ref{e.hw}) and Lemma~\ref{l.closeproblb}, and the third from (\ref{e.patterns}) (with a relabelling of $c > 0$). Thus, $\psap_{n+1} \big( N_{I}^1 \leq \deltanew (n+1) \big) \leq e^{-cn}$. The same holds for the quantity $N_{II}^1$. Considering 
$\Ga_{[3(n+1)/4, n+1]}$ in place of $\Ga_{[0,(n + 1)/4]}$, the same conclusion may be reached about $N_I^2$ and $N_{II}^2$. 
  It follows that 
  \begin{align*}
    \psap_{n+1} \Big( \min \big\{ \vert S_1 \vert , \vert S_2 \vert , N_I , N_{II} \big\} < \deltanew (n+1) \Big) < 4e^{-c n} \, .
  \end{align*}
This completes the proof. \qed  

\medskip

In the next lemma, we see how, for any $\ga \in \mathcal{G}_{n+1,\deltanew}$, the mixing of patterns that occurs when an element of $\vaso(\ga)$ is realized involves an asymptotically Gaussian fluctuation in the pattern number $N_I^1$. The statement and proof are minor variations of those of \cite[Lemma 3.5]{ontheprob}.

\begin{lemma}\label{l.balanced}

For any $\deltanew > 0$, there exists $c > 0$ and $N \in \N$ such that, for $n \geq N$ and $\ga \in \mathcal{G}_{n+1,\deltanew}$,

\medskip

\noindent{\bf (1)}  if $k \in \N$ satisfies 
  $\left| k - \frac{T_I |S_1|}{|S_1|+|S_2|}\right|  \leq  n^{1/2} (\log n)^{1/4}$,  then
$$ 
\psap_{n+1}\Big( N^1_I(\Ga) = k \, \Big\vert \, \Gamma \in \vaso(\ga) \Big)
    \geq n^{-1/2} \exp \big\{ - c (\log n)^{1/2} \big\} \, ;    
$$
\noindent{\bf (2)} and, for any  $g \in [1/8,1/2]$,
$$ 
\psap_{n+1}\bigg( \left| N^1_I(\Ga) -  \frac{N_I |S_1|}{|S_1|+|S_2|}\right|  \geq  n^{1/2} (\log n)^{g} \, \bigg\vert \, \Gamma \in \vaso(\ga) \bigg)
    \leq  \exp \big\{ - c (\log n)^{2g} \big\} \, .    
$$
\end{lemma}

\noindent{\bf Proof.}
Fix a choice of $\ga$ in $\mathcal{G}_{n+1,\deltanew}$.   If $\Ga$ is distributed according to $\psap_{n+1}$ conditionally on $\Gamma \in \vaso(\ga)$, then 
  $N_I$ type $I$ patterns and $N_{II}$ type $II$ patterns are distributed uniformly in the slots of $S_1 \cup S_2$. 
  Thus, for $k \in \big\{0, \cdots, \vert S_1 \vert \big\}$,
  \begin{align}\label{e.Tdistrib}
    \psap_{n+1}\Big( N^1_I(\Ga) = k  \, \Big\vert \, \Gamma \in \vaso(\ga) \Big) 
    = \frac{\binom{|S_1|}{k} \binom{|S_2|}{N_I - k}}{\binom{|S_1|+|S_2|}{N_I}}.
  \end{align}

  Write $m = |S_1| + |S_2|$, $|S_1| = \alpha m$ and   $N_I = \beta m$.
  By assumption $\alpha, \beta \in [\deltanew , 1-\deltanew] $ and $m \geq 2\deltanew (n+1)$. 
  Let $Z = \frac{N^1_I}{\alpha \beta m} -1$. 
  Under $\psap_{n+1} \left( \cdot \mcond \Gamma \in \vaso(\ga) \right)$, 
  $Z$ is a random variable of mean $0$,
  such that $\alpha \beta (1 + Z) m \in \ZZ \cap [0, \min \{ |S_1|,T_I \} ]$.
  
  First, we investigate the case where $Z$ is close to its mean.
  By means of a computation which uses Stirling's formula and \eqref{e.Tdistrib},  we find that
  \begin{align}\label{e.zdistrib}
 \psap_{n+1}\left(Z = z \mcond  \Gamma \in \vaso(\ga)  \right) 
    = \big( 1 + o(1) \big)
    \frac{ \exp \left(- \frac{\alpha \beta}{2(1-\alpha)(1-\beta)} m z^2 \right)}
    {\sqrt{2 \pi \alpha \beta (1 - \alpha)(1 - \beta) m}} \, ,     
 \end{align}
  where $o(1)$ designates a quantity tending to $0$ as $n$ tends to infinity, 
  uniformly in the acceptable choices of $\gamma$, $S_1$, $S_2$ and $z$, 
  with $|z| \leq \frac{2 n^{1/2}(\log n)^{1/2 +\eps}}{\alpha \beta m}$. We have obtained Lemma~\ref{l.balanced}(1).

  We now turn to the deviations of $Z$ from its mean. 
  From \eqref{e.Tdistrib}, one can easily derive that
  $\psap_{n+1}\left(Z = z \mcond \Gamma \in \vaso(\ga) \right)$ 
  is unimodal in $z$ with maximum at the value closest to $0$ that $Z$ may take.
  (We remind the reader that $Z$ takes values in $\frac{1}{\alpha \beta  m } \Z - 1$, 
  which contains $0$ only if ${\alpha \beta  m } \in \Z$.)  
The asymptotic equality
\eqref{e.zdistrib} thus implies the existence of constants $\mathrm{c}_0, \mathrm{c}_1 > 0$ depending only on $\deltanew$ such that, for $|z| \geq \frac{n^{1/2}(\log n)^{g}}{\alpha \beta m}$ and $n$ large enough, 
  \begin{align}\label{e.balanced0}
    \psap_{n+1} \Big( Z = z  \, \Big\vert \, \Gamma \in \vaso(\ga) \Big) \leq
  \mathrm{c}_1^{-1}  n^{-1/2}    \exp\left\{ - \mathrm{c}_0 ( \log n)^{2g}\right\}   \, ;
  \end{align}
  while for given $\eps > 0$, $|z| \geq \frac{n^{1/2}(\log n)^{1/2 + \eps}}{\alpha \beta m}$ and $n$ large enough,
  \begin{align}\label{e.balanced1}
    \psap_{n+1} \Big( Z = z \, \Big\vert \, \Gamma \in \vaso(\ga) \Big) \leq
   \mathrm{c}_1^{-1}  n^{-1/2}    \exp\left\{ - \mathrm{c}_0 ( \log n)^{1 + 2\eps}\right\}  \, \leq n^{- 2} \, .
  \end{align}

  Since $T^1_I$ takes no more than $n+1$ values, \eqref{e.balanced0} and~(\ref{e.balanced1}) imply Lemma~\ref{l.balanced}(2).\qed

\subsection{Mixing patterns by a random resampling}

Consider a random resampling experiment whose law we will denote by~$\pexper$. 
First an input polygon $\gamin$ is sampled according to the law $\psap_{n+1}$.
Then the contents of the slots in $S_1 \cup S_2$ are forgotten and independently resampled to form an output polygon $\gamout \in \sap_{n + 1}$. That is, given $\gamin$, $\gamout$ is chosen uniformly among the set of polygons $\ga \in \sap_{n + 1}$ for which $\ga \in \vaso(\gamin)$. Explicitly, if there are $j$ type $II$ patterns among $k$ slots in $S_1 \cup S_2$ in $\gamin$ (so that $k \geq j$), the polygon $\gamout$ is formed by choosing uniformly at random a subset of cardinality $j$ of these $k$ slots and inserting type $II$ patterns into the chosen slots.

Note the crucial property that $\gamout$ under $\pexper$ has the law $\psap_{n + 1}$: the resampling experiment holds the length~$n + 1$ random polygon at equilibrium. We mention in passing that a basic consequence of this resampling is a delocalization of the walk midpoint.
\begin{proposition}{\cite[Proposition 1.3]{ontheprob}}
Let $d \geq 2$. There exists $C > 0$ such that, for $m \in \N$,
$$
 \sup_{x \in \Z^d} \psaw_{m} \Big( \Ga_{\lfloor m/2 \rfloor} = x \Big) \leq C m^{-1/2} \, .
$$
\end{proposition}
It may be instructive to consider how to proof this result using the resampling experiment (or in fact a similar one involving walks rather than polygons) and Lemma~\ref{l.kptconseq}; a proof using such an approach is given in \cite[Section 3.2]{ontheprob}.

Any element $\phi \in \vaso(\ga)$ 
begins by tracing a journey over the region where slots in $S_1$ may appear, from the origin to $\gaemp_{(n+1)/10}$; it then follows its {\em middle section}, the trajectory of $\ga$ from  $\gaemp_{(n+1)/10}$ until  $\gaemp_{l' - (n+1)/10}$ (where $l' = \ellempty$); and it ends by moving from this vertex back to the origin, through the territory of slots in~$S_2$. 
Note that, in traversing the middle section, $\phi$ is exactly following a sub-walk of $\ga$, because no pattern changes have been made to this part of~$\ga$.  The timing of the middle section of this trajectory is advanced or retarded according to how many type~$II$ patterns are placed in the slots  in $S_1$.
Each extra such pattern retards the schedule by two units.
When $\phi$ has the {\em minimum} possible number $m : = \max \{ 0 , T_{II}(\ga)  - \vert S_2 \vert  \}$ of type $II$ patterns in the slots of $S_1$, the middle section is traversed by $\phi$ as early as possible,  the journey taking place during $\big[(n+1)/10 + 2m, n - 2 \big( T_{II}(\ga) - m \big) - (n+1)/10 \big]$. When $\phi$ has the maximum possible number $M : = \min \{ \vert S_1 \vert , T_{II}(\ga) \}$ of type~$II$ patterns in the slots of $S_1$, this traversal occurs as late as possible,  during $\big[ (n+1)/10 + 2 M , (n+1)  - (n+1)/10 - 2 \big( T_{II}(\ga) - M \big) \big]$. Since
$M \leq \vert S_1 \vert \leq (n+1)/13$ and $T_{II}(\ga) - m  \leq \vert S_2 \vert \leq (n+1)/13$, $\phi_j$ necessarily lies in the middle section whenever $j \in [(n+1)/10 + 2(n+1)/13,n+1 - (n+1)/10 -2(n+1)/13]$. Since the snake method parameter $\ell$ has been set to belong to the interval $[n/4,3n/4]$, we see that $\phi_j$ always lies in the middle section whenever $j \in \big[ \ell - n^{1/2}, \ell + n^{1/2} \big]$.

Taking $\ga \in \sap_{n + 1}$ and conditioning $\pexper$ on $\gamin = \ga$, note that the mean number of type $I$ patterns that end up in the slots in $S_1$ under $\gamout$ is given by $T_I(\ga) \cdot \tfrac{\vert S_1(\ga) \vert}{\vert S_1(\ga) \vert + \vert S_2(\ga) \vert}$, because this expression is the product of the number of type~$I$ patterns that are redistributed and the proportion of the available slots that lie in~$S_1$.

Consider now a polygon $\phi \in \vaso(\ga)$ that achieves as closely as possible the mean value for the number of type $I$ patterns among the slots in $S_1$: that is, $T_I^1(\phi)$ equals $\lfloor T_I(\ga) \cdot \tfrac{\vert S_1(\ga) \vert}{\vert S_1(\ga) \vert + \vert S_2(\ga) \vert} \rfloor$. As we have noted, $\phi_{\ell}$ is always reached during the middle section of $\phi$'s three-stage journey. Define  the {\em middle index}  $\lmid = \lmid(\ga)$ so that $\phi_{\ell} = \ga_{\lmid}$.   
Note that, given $\ga$ and this value of  $T_I^1(\phi)$, the value of this index is independent of the choice of $\phi$.

\subsection{Snakes of walks with high closing probability are typical}

Recall that the index parameter $\ell$ (and also $n$)
were fixed in Subsection~\ref{s.indexset}.
Moving towards the proof of Proposition~\ref{p.snakeg}, we take $\ga \in \sap_{n + 1}$ and define $\notick(\ga)$ to be the set  
\begin{eqnarray*}
     & & \Big\{ j \in \big(\ell - 2\N\big) \, \cap \, \big[ \lmid(\ga) - 2 n^{1/2} (\log n)^{1/4} \, , \,  \lmid(\ga) + 2 n^{1/2} (\log n)^{1/4} \big]: \\
   & & \qquad \qquad \qquad \qquad  \qquad \qquad  \qquad 
    \ga \, \, \textrm{is not $(1/2 - 2\eps,n,\ell)$-charming at $j$} \Big\} \, .
\end{eqnarray*}
Henceforth in this proof, charming will mean $(1/2 - 2\eps,n,\ell)$-charming. 
(The parity constraint that $\ell - j$ is even is applied above because the walks of length $j + n - \ell$ considered in the definition of charming at index~$j$ must be of odd length if they are to close.)

 We have seen that 
$\gamout_{\ell}$ is visited by $\gamout$ during its middle section, when it is traversing a subpath of $\gamin$ unchanged by pattern mixing. For this reason, we may define a random variable $L$ under $\pexper$ by setting  
$\gamout_{\ell} = \gamin_L$. 

We now state a key property of the resampling procedure. 
 \begin{lemma}\label{l.keyprop}
 The events $\big\{ \textrm{$\gamout$ is charming at $\ell$} \big\}$ and $\big\{ \textrm{$\gamin$ is charming at $L$} \big\}$ coincide.
\end{lemma}
\noindent{\bf Proof.} Note that the shells of $\vas\big( \gamin_{[0,L]} \big)$ and $\vas \big( \gamout_{[0,\ell]} \big)$ coincide, because   
$\gamout_{[0,\ell]}$ may be obtained from  $\gamin_{[0,L]}$ by modifying the $I/II$-status of some of its slots (these being certain slots in $S_1$). Thus,
Lemma~\ref{l.sharedavoidance} implies the statement. \qed

\begin{lemma}\label{l.jbound}
$$
\pexper \Big(   \big\vert \notick(\gamin) \big\vert \geq n^{1/2 - \eps/6}  \Big) \leq n^{-\eps/6} \, .
$$
\end{lemma}
\noindent{\bf Proof.}
Choosing $\deltanew > 0$ small enough and abbreviating $\mathcal{G} = \mathcal{G}_{n+1,\deltanew}$, Lemma~\ref{l.balanced}(1) implies that, 
for each $\ga \in \mathcal{G}$ and $k \in \big[  - n^{1/2} (\log n)^{1/4}, n^{1/2} (\log n)^{1/4} \big]$,
$$
 \pexper \Big(  L  = \lmid(\gamin) +  2k    \,  \Big\vert \,  \gamin =  \ga \Big) \geq  n^{-1/2} \exp \big\{ - c (\log n)^{1/2} \big\} \, .
$$ 
Thus, again taking any $\gamma \in \mathcal{G}$,
\begin{eqnarray*}
 & & \pexper \Big( \gamout \, \, \textrm{is not charming at $\ell$} \, \Big\vert \, \gamin  = \ga \Big) \\
 & \geq & \sum_{k =   - n^{1/2} (\log n)^{1/4}}^{n^{1/2} (\log n)^{1/4}}  \pexper \Big( \gamout \, \, \textrm{is not charming at $\ell$} \, , \, L =  \lmid(\gamin) + 2k \,  \Big\vert \, \gamin  = \ga \Big) \\
 & \geq & \sum_{k =   - n^{1/2} (\log n)^{1/4}}^{n^{1/2} (\log n)^{1/4}}  \pexper \Big(  \ga \, \, \textrm{is not charming at $\lmid(\gamin) + 2k$} \,   \, , \, L =  \lmid(\gamin) + 2k  \,  \Big\vert \, \gamin  = \ga \Big) \\
 & \geq &   n^{-1/2} \exp \big\{ - c (\log n)^{1/2} \big\} \sum_{k =   - n^{1/2} (\log n)^{1/4}}^{n^{1/2} (\log n)^{1/4}}  1\!\! 1_{    \ga \, \, \textrm{is not charming at $\lmid(\ga) + 2k$}}   \\
 & \geq &   n^{-1/2} \exp \big\{ - c (\log n)^{1/2} \big\} \cdot    \vert \notick(\ga) \vert  \, ,
\end{eqnarray*}
where the second inequality made use of Lemma~\ref{l.keyprop}.

Averaging over such $\gamma$, we find that
\begin{eqnarray*}
 & & \pexper \Big(  \gamout \, \, \textrm{is not charming at $\ell$} \, \Big\vert \, \gamin  \in \mathcal{G} \Big) \\
& \geq & c n^{-1/2}  \exp \big\{ - c (\log n)^{1/2} \big\}  \cdot \eexper \Big[ \,  \big\vert \notick(\gamin) \big\vert \, \Big\vert \, \gamin  \in \mathcal{G} \Big] \, ,
\end{eqnarray*}
where $\eexper$ denotes the expectation associated with the law $\pexper$.

Note that
\begin{eqnarray*}
 & & \pexper \Big( \gamout \, \, \textrm{is not charming at $\ell$} \, \Big\vert \, \gamin  \in \mathcal{G} \Big) \\
 & \leq &  2 \, \psap_{n+1} \Big(  \Gamma \, \, \textrm{is not charming at $\ell$} \Big) 
  =  2 \, \psap_{n+1} \Big( \Gamma_{[0,\ell]} \not\in \hPhi_{\ell,n}^{1/2 - 2\eps}  \Big) \leq 2 n^{-\eps} \, ,
\end{eqnarray*}
where in the first inequality we use that  $\gamout$ under $\pexper$ has the law $\psap_{n+1}$, and then apply Lemma~\ref{l.kptconseq}, to find that $\pexper \big( \gamin \in \mathcal{G} \big) = \psap_{n+1} \big( \Gamma \in \mathcal{G} \big) \geq 1 - e^{-cn} \geq 1/2$. The final inequality above used Lemma~\ref{l.polyclose}.

Thus,
$$
\eexper \Big[ \,  \big\vert \notick(\gamin) \big\vert \, \Big\vert \, \gamin  \in \mathcal{G} \Big] \leq  n^{1/2 - \eps/2} \, .
$$

We find that $\eexper \vert \notick \vert$ is at most
\begin{eqnarray*}
  &  & 
\eexper \Big[ \,  \big\vert \notick(\gamin) \big\vert \, \Big\vert \, \gamin  \in \mathcal{G} \Big] \, + \, \Big( 2 n^{1/2} \big( \log n \big)^{1/4} + 1 \Big)  \pexper \Big( \gamin \not\in \mathcal{G} \Big) \\
   & \leq & n^{1/2 - \eps/2} \, + \, \Big( 2 n^{1/2} \big( \log n \big)^{1/4} + 1 \Big) \, e^{-cn}  \leq n^{1/2 - \eps/3} \, . 
\end{eqnarray*}

Applying Markov's inequality yields Lemma~\ref{l.jbound}. \qed

\medskip

\noindent{\bf Proof of Proposition~\ref{p.snakeg}.}
Lemma~\ref{l.balanced}(2) implies that 
$$
\psap_{n+1} \Big( \ell \in \big[ \lmid(\Ga) -  n^{1/2} (\log n)^{1/4} ,  \lmid(\Ga) +  n^{1/2} (\log n)^{1/4} \big] \, \Big\vert \,  \Ga \in \mathcal{G}_{n,\deltanew} \Big) 
$$
is at least $1 - \exp \big\{ - c (\log n)^{1/2} \big\}$. Note that the interval centred on $\lmid(\Ga)$ considered here is shorter than its counterpart in the definition of $\notick(\ga)$ for $\ga \in \sap_{n+1}$.
Applying Lemmas~\ref{l.kptconseq} and~\ref{l.jbound}, we find that
\begin{eqnarray*}
 & &  \psap_{n+1} \Big( \# \Big\{  j \in \big[ \ell - n^{1/2} (\log n)^{1/4} ,  \ell + n^{1/2} (\log n)^{1/4} \big]: \\
   & & \qquad \qquad   \textrm{$\Gamma$ is not charming at $j$} \Big\} \geq n^{1/2 - \eps/6}   \Big) 
 \leq   n^{-\eps/6} +  e^{- c (\log n)^{1/2}}    + e^{-cn} \, .
\end{eqnarray*}

When the complementary event occurs, $\Gamma$ is charming for at least one-quarter of the indices in  $\big[ \ell - n^{1/2}, \ell  \big]$, so that $\Gamma_{[0,\ell]}$ is an element of $\mathsf{CS}_{1/2,0}^{1/2 - 2\eps,\ell,n}$. (We write one-quarter rather than one-half here, because one-half of such indices are inadmissible due to their having the wrong parity.) \qed

\part{The snake method in unison with polygon joining}\label{c.aboveonehalf}

By Theorem~\ref{t.polydev}, $\thet_n \geq 3/2 - o(1)$ for typical $n$, and thus by~(\ref{e.closepc}) and $c_n \geq \mu^n$,  the closing probability is at most order~$n^{-1/2}$. This order is also the limit of the snake method via Gaussian pattern fluctuation. This part is devoted to the proof of the two results that push the closing probability upper bound below the level~$n^{-1/2}$: Theorem~\ref{t.thetexist}(1) and~(2).

In each case, our method will combine the two main approaches used thus far: it is the snake method implemented via a certain form of the polygon joining technique. It is the use of this technique which means that we work in dimension $d=2$ throughout Part~\ref{c.aboveonehalf}.

Since the length parameter $\beta$ in the snake method must be at most one, the snake method may at best prove closing probability upper bounds of the order $n^{-1}$.
The hypotheses of  Theorem~\ref{t.thetexist}(2) certainly entail the existence of the closing exponent~$-\lim_{n \in 2\N + 1} \big( \log n \big)^{-1} \log \psaw_n \big( \Gamma  \closes \big)$, since under these hypotheses, this exponent must equal $\thet + \xi -1$ due to~(\ref{e.closepc}); 
and Theorem~\ref{t.thetexist}(2) covers one-third of the interval from the solved~$1/2$ to the possible-in-principle~$1$, to reach a bound of~$2/3$ on the value of this exponent. Despite its conditional nature, the result is attractive for expository reasons: the proof of Theorem~\ref{t.thetexist}(1) shares the same framework, but  a number of technicalities that arise in its proof are absent in the conditional result's, permitting a more attentive focus on the main concepts. For this reason, we present the latter proof first.

\section{Proving Theorem~\ref{t.thetexist}(2)}

It is enough to suppose the existence of  $\theta$ and $\xi$
and to argue that $\thet + \xi \geq 5/3$, since the latter assertion of  Theorem~\ref{t.thetexist}(2) was justified as it was stated.  In order to prove this inequality by finding a contradiction,  we suppose that $\thet + \xi < 5/3$.
The sought contradiction is exhibited in the next result.

\begin{proposition}\label{p.thetexist}
Let $d = 2$. Assume that the limits~$\thet : = \lim_{n \in 2\N} \thet_n$ and $\xi =  \lim_{n \in \N} \xi_n$ exist. Suppose further that $\thet + \xi < 5/3$.  Then, for some $c > 1$ and $\delta >0$, the set of $n \in 2\N + 1$ for which 
$$
\psaw_n \big( \Ga \closes \big) \leq c^{-n^{\delta}}
$$
intersects the shifted dyadic scale $\big[ 2^i - 1,2^{i+1} -1 \big]$ for all but finitely many $i \in \N$.
\end{proposition}
\noindent{\bf Proof of Theorem~\ref{t.thetexist}(2).}
The non-negative limits $\thet$ and $\xi$ are hypothesised to exist. If their sum is finite, then the formula  (\ref{e.closepc.formula}) exhibits a polynomial decay. 
When the sum is less than $5/3$, this 
  contradicts Proposition~\ref{p.thetexist}. \qed

It remains of course to derive the proposition. 
Assume throughout the derivation that the proposition's hypotheses are satisfied. Define an exponent $\chi$ so that  $\theta + \xi = 3/2 + \chi$. 
Since $\xi \geq 0$ and $\theta \geq 3/2$ (classically and by Theorem~\ref{t.polydev}), $\chi$ is non-negative.
Note that $\psaw_n \big( \Ga \closes \big)$ equals $n^{-1/2 - \chi + o(1)}$, and also that  we are supposing  that $\chi < 1/6$. 

The proof of Proposition~\ref{p.thetexist} is an application of the snake method.
We now set the snake method exponent parameters for this application, taking 
\begin{itemize}
\item the snake length $\beta$ equal to one; 
\item
the inverse charm $\alpha$ equal to $1/2 + 2\chi + 4\eps$;
\item and the deficit $\esm$ equal to $\chi + 8 \eps$.
\end{itemize}
 where henceforth in proving Proposition~\ref{p.thetexist}, $\eps > 0$ denotes a given but arbitrarily small quantity. The quantity $\delta = \beta - \esm - \alpha = 1/2 - 3\chi - 12\eps$ must be positive if the snake method is to work, and thus we choose $\eps < (1/6 - \chi)/4$.

Recall the heuristic derivation via polygon joining that $\thet \geq 5/2$ in~$d=2$, and the associated notion of a $(\eone,\etwo,\ethree)$-argument from Section~\ref{s.sumthree}.
The snake method derivation of  
Proposition~\ref{p.thetexist} may also be viewed as
an $(1/2,1,\chi-\varphi)$-argument for $\varphi > 0$ arbitrarily small. Some care is needed in interpreting this point of view, however: given the $(1/2,1,0)$-argument that yields Theorem~\ref{t.polydev}, it would seem that, in proposing an $(1/2,1,\chi-\varphi)$-argument
to reach  
Proposition~\ref{p.thetexist}, we would aim to prove that the $\psap_n$-probability of a macroscopic join plaquette is at most $n^{-\chi + \varphi}$. In fact, we will not succeed in making this geometric inference. Rather, our assumption that $\chi < 1/6$ will force a scenario in which this probability exceeds the value  $n^{-1/6}$, so that macroscopically joined polygons are fairly commonplace. The snake method will be implemented by exploiting the prevalence of such polygons. Ultimately the assumption that $\chi < 1/6$ will be contradicted. Though we will have learnt that $\chi \geq 1/6$, we will not have learnt that the probability of macroscopically joined polygons is low; indeed, in terms of the polygon joining  $(\eone,\etwo,\ethree)$-argument, it could be that $\chi \geq 1/6$ is due to say $\eone$ being at least $1/2 + 1/6$ (i.e., the Flory exponent governing diameter being at least $2/3$), rather than because $\ethree$ is positive. Despite this, it is meaningful to say that we are undertaking a $(1/2,1,\chi -\varphi)$-argument, rather than say a  $(1/2 + \chi -\varphi,1,0)$-argument, because, during this proof by contradiction, it is the assumption that the $\ethree$-value is bounded above by $\chi - \varphi$ that forces the prevalence of macroscopically joined polygons that is central to the argument.

\section{Three steps to the proof of Proposition~\ref{p.thetexist}}

The argument leading to Proposition~\ref{p.thetexist} has three principal steps. We explain these now in outline, stating the conclusion of each one. The proofs follow, in the three sections that follow this one.

\subsection{Step one: finding many first parts with high closing probability}

The outcome of the first step may be expressed as follows, using the high conditional closing probability set notation~$\hPhi$ introduced in Subsection~\ref{s.possparts}.
\begin{proposition}\label{p.notquite}
For $i \in \N$ sufficiently high, there exists $m' \in 2\N \cap \big[ 2^{i+3}, 2^{i+4} \big]$ such that, writing $\lset$ for the set of values $k \in \N$, $2^i \leq k \leq 2^i + 2^{i-2}$, that satisfy
$$
\psap_{m'+k} \Big( \Ga_{[0,k]} \in \hPhi_{k,k+m' - 1}^\alpha \Big)
 \geq 1 - 2^{-i\chi} \, ,
$$
we have that $\vert \lset \vert \geq 2^{-8}m'$.
\end{proposition}

This first step is a soft, Fubini, argument: such index pairs $(k,m'+k)$ are characteristic given our assumption on closing probability decay, and the proposition will follow from Lemma~\ref{l.closecard}.

\subsection{Step two: individual snake terms are often charming}
To apply the snake method, we must verify its fundamental hypothesis~(\ref{e.afewcs}) that charming snakes are not rare under the uniform polygon law. Proposition~\ref{p.notquite}
is not adequate for verifying this hypothesis, because the polygon law index $m' + k$ is a variable that changes with $k$.
We want a version of the proposition where this index is fixed. 
The next assertion, which is the conclusion of our second step, is such a result. Note from the differing forms of the right-hand sides in the two results 
that a behaviour determined to be highly typical in Proposition~\ref{p.notquite} has a counterpart in Proposition~\ref{p.almost} which is merely shown to be not unusual.
\begin{proposition}\label{p.almost}
For each $i \in \N$ sufficiently high, there exist $m \in 2\N \cap \big[ 2^{i+4}, 2^{i+5} \big]$ and $m' \in \big[ 2^{i+3},2^{i+4} \big]$ such that, writing $K$ for the set of values $k \in \N$, $1 \leq k \leq m-m'$, that satisfy
$$
\psap_m \Big( \Ga_{[0,k]} \in \hPhi_{k,k+m' - 1}^\alpha \Big)
 \geq   m^{-\esm + \eps} \, ,
$$
we have that $\vert K \vert \geq 2^{-9} m$.
\end{proposition}

\begin{figure}
  \begin{center}
    \includegraphics[width=1.0\textwidth]{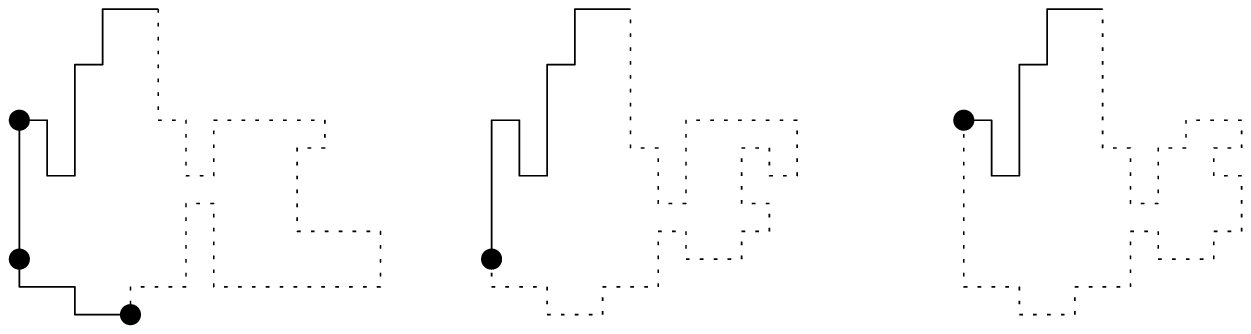}
  \end{center}
  \caption{The snake method's new guise: compare with Figure~\ref{f.snakemethodpattern}. Again black dots mark vertices of  polygons $\ga$: fixing a typical $\tilden \in 2\N$, assign a black dot at each vertex $\ga_j$, $\tilden/4 \leq j \leq 3\tilden/4$, of a polygon $\ga$ of length of order~$\tilden$, that has the property that $\saw^0_{\tilden+j} \big( \Ga \closes \, \big\vert \, \vert \Ga^1 \vert = j , \Ga^1 = \ga_{[0,j]} \big) \geq (\tilden)^{-1/2 - \chi - \eps}$.
There is a difference to Figure~\ref{f.snakemethodpattern}, because polygons $\ga$ may differ in length. 
Our assumptions readily imply that most polygons $\ga$ with length $n'$ of order say~$2\tilden$ will have a black dot at $\ga_{n' - \tilden}$. What about at other locations $\ga_j$ where $j$ is also of order $\tilden$? A black dot is likely to appear at $\Ga_j$ when $\Ga$ is $\psap_{j+\tilden}$-distributed. If we can make a {\em comparison of measure}, showing that the laws $\psap_{j+\tilden}$ and $\psap_{{n'}}$ are to some degree similar, then the black dot will be known to also appear at~$\ga_j$ in a typical sample of $\psap_{{n'}}$. Applying this for many such~$j$, black dots will appear along the course of a typical length~$n'$ polygon. These black dots index charming snake terms and permit the use of the snake method.
Comparison of measure  will be undertaken in Proposition~\ref{p.similarity} by the use of polygon joining. Lemma~\ref{l.polyepsregthet} shows that a non-negligible proportion of polygons are regulation global join polygons; thus, such polygons themselves typically have black dots~$\tilden$ steps from their end.
The law of the left polygon under the uniform law on regulation polygons of a given length is largely unchanged as that total length is varied on a given dyadic scale,
because the length discrepancy can be absorbed by altering the length of the right polygon. Indeed, in the above sketches, we
deplete the length of the right polygon in a possible extension of the depicted first part, as the length of this first part shortens; in this way, we show that first parts typically arising in regulation polygons at one length are also characteristic in such polygons with lengths on the same scale.}\label{f.snakemethodpolygon}
\end{figure}

The second step, leading to the proof of this proposition, is quite subtle. Expressed in its simplest terms, the idea of the proof is  that the subscript index in $\psap$ can be changed from the $k$-dependent $m' + k$ in Proposition~\ref{p.notquite}  to the $k$-independent~$m$ in the new result by using the polygon joining technique.
This index change involves proving a {\em similarity of measure} between polygon laws with distinct indices on a given dyadic scale. The set of regulation global join polygons introduced in Section~\ref{s.pjprep} provides a convenient reservoir of macroscopically joined polygons whose left polygons are shared between polygon laws of differing index, so that this set provides a means by which this similarity may be proved. 
Figure~\ref{f.snakemethodpolygon} offers some further outline of how we will exploit the ample supply of regulation polygons to prove similarity of measure.

\subsection{Step three: applying the snake method Theorem~\ref{t.snakesecondstep}}
In this step, it is our goal to use Proposition~\ref{p.almost}
in order to verify the snake method hypothesis~(\ref{e.afewcs}). 
In the step, we will  fix the method's two index parameters: for $i \in \N$ given, we will use Proposition~\ref{p.almost} to take $n$ equal to~$m-1$, and $\ell$ equal to $m - m'$. 

We may now emphasise why the proposition is useful for this verification goal.
Note that $\max K \leq \ell$.
Since $n - \ell = m' - 1$, we have that, for $k \in [1,\ell]$, the walk $\gamma \in \fpart_{\ell,n}$ is $(\alpha,n,\ell)$-charming at index $k$ if and only if 
$$
 \psaw^0_{k+m'-1} \Big( \Ga \closes \, \Big\vert \, \big\vert \Gamma^1 \big\vert = k \, , \,  \Ga^1_{[0,k]} = \ga_{[0,k]} \Big) > n^{-\alpha} \, ;
$$ 
note the displayed condition is almost the same as $\ga_{[0,k]} \in \hPhi_{k,m' +k-1}^{\alpha}$, the latter occurring when $n^{-\alpha}$ is replaced by~$(m'+k-1)^{-\alpha}$, a change which involves only a bounded factor. 

Thus, Proposition~\ref{p.almost} shows that individual snake terms are charming with probability at least~$m^{-\chi - o(1)}$.
It will be a simple Fubini argument that will allow us to verify~(\ref{e.afewcs}) by confirming that it is not rare that such terms gather together to form charming snakes. Theorem~\ref{t.snakesecondstep} may be invoked to prove Proposition~\ref{p.thetexist} immediately. The third step will provide this Fubini argument.

The next three sections are devoted to implementing these three steps. In referring to the dyadic scale index $i \in \N$ in results stated in these sections, we will sometimes neglect to record that this index is assumed to be sufficiently high.

\section{Step one: the proof of Proposition~\ref{p.notquite}}\label{s.notquite}

Before proving this result, we review its meaning in light of  Figure~\ref{f.snakemethodpolygon}. 
Proposition~\ref{p.notquite} is a precisely stated counterpart to the assertion in the figure's caption that a black dot typically appears $\tilden$ steps from the end of a polygon drawn from $\psap_{n'}$ where $n'$ has order $2\tilden$ (or say $2^{i+3}$). Note the discrepancy that we have chosen $\alpha$ to be marginally above $1/2 + 2\chi$ while in the caption a choice just above $1/2 + \chi$ was made. The extra margin permits $1 - 2^{-i\chi}$ rather than $1 - o(1)$ in the conclusion of the proposition. We certainly need the extra margin; we will see how  at the end of step two.

Consider two parameters $\eclo > \eroom > 0$. (During the proof of Proposition~\ref{p.almost}, we will ultimately set  $\eclo = 4\eroom$ and $\eroom = \eps$, where $\eps > 0$ is the parameter that has been used to specify the three snake method exponents.)

\begin{definition}\label{d.cpt}
Let $i \in \N$. An index pair $(k,j) \in \N \times 2\N$ will be called {\em closing probability typical} on dyadic scale $i$ if 
\begin{itemize}
\item $j \in \big[ 2^{i+4},  2^{i+5} \big]$ and $k \in \big[ 2^i , 2^i + 2^{i-2} \big]$;
\item and, for any $a \geq 0$,
\begin{equation}\label{e.acondition}
 \psap_{j} \Big( \Gamma_{[0,k]} \notin \phicl{k}{j-1}{1/2 + (a+1)\chi + \eclo}  \Big) \leq  j^{-(a\chi + \eclo - \eroom)} \, .
\end{equation}
\end{itemize}
\end{definition}

\begin{proposition}\label{p.manyethet}
There exists $k \in 2\N \cap \big[ 2^{i+4}, 2^{i+5} - 2^{i-2} \big]$ for which there are at least $2^{i-4}$ values of $j \in 2\N \cap [0,2^{i-2}]$ with the index pair
$\big( 2^i + j, k + j \big)$ closing probability typical on dyadic scale $i$.
\end{proposition}

When we prove Proposition~\ref{p.almost} in step two, we will in fact do so using a direct consequence of Proposition~\ref{p.manyethet}, rather than Proposition~\ref{p.notquite}. 
The latter result, a byproduct of the proof of  Proposition~\ref{p.manyethet}, has been stated merely because it permitted us to make a direct comparison between steps one and two in the preceding outline.
As such, our goal in step one is now to prove Proposition~\ref{p.manyethet}.
\begin{definition}\label{d.e}
Define $E$ to be the set of pairs $(k,j) \in \N \times 2\N$,  
\begin{itemize}
\item where $j \in \big[ 2^{i+4},  2^{i+5} \big]$ and $k \in \big[ 2^i , 2^i + 2^{i-2} \big]$;
\item and the pair $(k,j)$ is such that
\begin{eqnarray}
 & & \# \Big\{ \gamma \in \saw^0_{j-1}: \vert \gamma^1 \vert = k \Big\} \label{e.rchithet} \\
 & < & (j-1)^{1/2 + \chi + \eroom} \cdot  \# \Big\{ \gamma \in \saw^0_{j-1}: \vert \gamma^1 \vert = k , \gamma \closes \Big\} \, . \nonumber
\end{eqnarray}
\end{itemize}
\end{definition}

We will now apply Lemma~\ref{l.closecard} to show that membership of $E$ is typical. The application is possible because by hypothesis there is a lower bound on the  decay of the closing probability. Indeed, we know that
$\psaw_n \big( \Gamma \closes \big) \geq n^{-1/2 - \chi - o(1)}$. Thus, the lemma with $\alphamac = 1/2 + \chi + \eroom/2$ and $\deltamac = \eroom/2$ proves the following.
\begin{lemma}\label{l.etypical}
Provided the index $i \in \N$ is sufficiently high,
for each $j \in  2 \N \cap \big[ 2^{i+4},  2^{i+5} \big]$, the set of $k \in \big[ 2^i , 2^i + 2^{i-2} \big]$ such that $(k,j) \not\in E$ has cardinality at most $2 j^{1 - \eroom/2} \leq 2 \big(  2^{i+5} \big)^{1-\eroom/2}$.
\end{lemma}

\begin{lemma}\label{l.rjsjthet}
Any pair $(k,j) \in E$
is closing probability typical on dyadic scale $i$.
\end{lemma}
\noindent{\bf Proof.} We must verify that $(k,j)$
verifies~(\ref{e.acondition}) for any $a \geq 0$. Recalling the notation $\psaw_{j-1}^0$ from Subsection~\ref{s.sawfree}, note that the assertion
\begin{equation}\label{e.lststhet}
 \psaw^0_{j-1} \Big( \Gamma^1 \notin \phicl{k}{j-1}{1/2 + (a+1)\chi + \eclo} \, \, \Big\vert \, \, \Gamma \closes \, ,  \vert \Gamma^1 \vert  = k \Big) \leq j^{-(a\chi + \eclo - \eroom)} 
\end{equation}
is a reformulation of this condition.
Note that
\begin{equation}\label{e.cardgamma}
  \# \Big\{ \gamma \in \saw^0_{j-1}: \vert \gamma^1 \vert =  k \, , \, \gamma^1 \notin \phicl{k}{j-1}{1/2 + (a+1)\chi + \eclo} \Big\} 
\end{equation}
is at least the product of $(j-1)^{1/2 + (a+1)\chi + \eclo}$
and
$$
 \# \Big\{ \gamma \in \saw^0_{j-1}: \vert \gamma^1 \vert = k \, , \, \gamma^1 \notin \phicl{k}{j-1}{1/2 + (a+1)\chi + \eclo} \, , \, \gamma \closes \Big\}  \, .
$$
The quantity~(\ref{e.cardgamma}) is bounded above by~(\ref{e.rchithet}), 
and~(\ref{e.lststhet}) is obtained. \qed

\medskip

\noindent{\bf Proof of Proposition~\ref{p.manyethet}.}
By Lemma~\ref{l.rjsjthet}, it is enough to argue that there exists $k \in 2\N \cap \big[ 2^{i+4}, 2^{i+5} - 2^{i-2} \big]$ for which there are at least $2^{i-4}$ values of $j \in 2\N \cap [0,2^{i-2}]$ with 
$\big( 2^i + j, k + j \big) \in E$.
 
By Lemma~\ref{l.etypical},
\begin{eqnarray}
 & & \sum_{j \in 2\N \cap [ 2^{i+4} + 2^{i-2} , 2^{i+4} + 2^{i} ] } \sum_{\ell = 2^i}^{2^i + 2^{i-2}} 
 1\!\!1_{(\ell,j) \in E} \label{e.ineqe} \\
 & \geq & \Big( 2^{i-2} + 1 - 2 \big( 2^{i+5} \big)^{1 - \eroom/2} \Big) \cdot \tfrac{1}{2} \big( 2^{i} - 2^{i-2} \big) \, . \nonumber
\end{eqnarray}
    
\begin{figure}
    \begin{center}
      \includegraphics[width=0.5\textwidth]{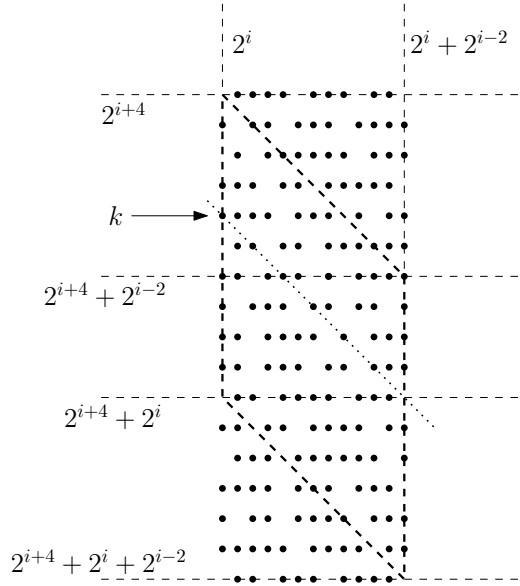}
    \end{center}
    \caption{Illustrating the proof of Proposition~\ref{p.manyethet}.    The dots depict elements of $E$. 
 The quantity~(\ref{e.ineqe})  is the number of dots in the middle rectangle; (in these terms, $\ell$ is the $x$-coordinate in the sketch and $j$ is the negative $y$-coordinate). Similarly,     
    the quantity~(\ref{e.doublesum}) is the number of dots in the parallelogram bounded by dashed lines.
    In~(\ref{e.kterm}), we find a diagonal, such as the one depicted with a pointed line, on which many dots lie. }\label{f.sumexplain}
  \end{figure}

As Figure~\ref{f.sumexplain} illustrates, the left-hand side here is bounded above by 
\begin{equation}\label{e.doublesum}
 \sum_{k \in 2\N \cap [ 2^{i+4},2^{i+4} + 2^{i} ]} \sum_{\ell = 2^i}^{2^i + 2^{i-2}}
 1\!\!1_{(\ell,k + \ell - 2^i) \in E} \, .
\end{equation}
There being $2^{i - 1} + 1 \leq 2^{i}$ indices $k \in 2\N \cap \big[ 2^{i+4},2^{i+4} + 2^{i} \big]$, one such~$k$ satisfies
\begin{equation}\label{e.kterm}
\sum_{\ell = 2^i}^{2^i + 2^{i-2}} 1\!\!1_{(\ell,k + \ell - 2^i) \in E} \geq \tfrac{3}{8} \, \Big( 2^{i-2} + 1 - 2 \big( 2^{i+5} \big)^{1 - \eroom/2} \Big) \, .
\end{equation}
Noting that $\tfrac{3}{4} \big( 2^{i+5} \big)^{1 - \eroom/2} \leq \tfrac{1}{8} 2^{i-2}$, we obtain the sought statement and so conclude the proof. \qed

\medskip

\noindent{\bf Proof of Proposition~\ref{p.notquite}.} 
Proposition~\ref{p.manyethet} furnishes the existence of a minimal $Q \in 2\N \cap \big[ 2^{i+4}, 2^{i+4} + 2^{i} \big]$ 
such that there are at least $2^{i-4}$ values of $k \in 2\N \cap [0,2^{i-2}]$ with  $\big( 2^i + k, Q + k \big)$ closing probability typical on dyadic scale~$i$. 

Let $\big( k_1,\ldots,k_{2^{i-4}} \big)$ be an increasing sequence of such values of $k$ associated to the value~$Q$. 
Set $m' = Q - 2^i$ and note that $m' \in \big[2^{i+3},2^{i+4} \big]$. 
Take $a=1$ in~(\ref{e.acondition}) and specify that $\eclo  = 4\eps$ alongside $\eclo > \eroom$  to find that, for all such~$j$,  $2^i + k_j \in \lset$. Thus, $\vert \lset \vert \geq 2^{i-4} \geq 2^{-8} m'$.

This completes the proof of the proposition. \qed

\medskip

We end this section by recording the explicit artefact that will be used in step two.
\begin{definition}\label{d.rjsj}
Set $\indexc = 2^{i-4}$. We specify a {\em constant difference} sequence of  index pairs $(r_j,s_j)$, $1 \leq j \leq \indexc$. Its elements are closing probability typical for dyadic scale $i$, the difference $s_j - r_j$ is independent of $j$ and takes a value in $[2^{i+4}- 2^i,2^{i+4}]$, and $\big\{ r_j: 1 \leq j \leq \indexc \big\}$ is increasing.
\end{definition}
Indeed, we may work with the sequence constructed in the preceding proof, 
setting $r_j = 2^i + k_j$ and $s_j = Q + k_j$  for $j \in \big[ 1, \indexc  \big]$. 

We are {\em en route} to Proposition~\ref{p.almost},
and take the opportunity to mention how its two parameters $m$ and $m'$ will be specified. The latter will be chosen, as it was in the proof of Proposition~\ref{p.notquite}, to be the value of the constant difference in the sequence in Definition~\ref{d.rjsj}, while the former will be taken equal to $s_{\indexc}$ in this definition.

\section{Step two: deriving Proposition~\ref{p.almost}}\label{s.almost}

The key tool for this proof is similarity of measure
between polygon laws with distinct indices, with the reservoir of regulation global join polygons providing the mechanism for establishing this similarity. Recall that these polygons were introduced and discussed in Section~\ref{s.pjprep}; the reader may wish to review this material~now. 

This step two section has three subsections. In the first, we provide a piece of apparatus, which augments the regulation polygon properties discussed in Section~\ref{s.pjprep}, that we will need in order to state our similarity of measure result. This result, Proposition~\ref{p.similarity}, is stated and proved in the next subsection. The final subsection gives the proof of Proposition~\ref{p.almost}.

\subsection{Relating polygon laws $\psap_n$ 
for distinct $n$ via global join polygons}\label{s.regnotrare}

 The next lemma shows that there is an ample supply of regulation polygons in the sense that they are not rare among all polygons.
In contrast to the regulation polygon tools in Section~\ref{s.pjprep}, Lemma~\ref{l.polyepsregthet}
is contingent on our standing assumption that the hypotheses of Proposition~\ref{p.thetexist} are satisfied.
(In fact, it is the existence and finiteness of $\theta$ and $\xi$ that are needed, rather than the condition $\theta + \xi < 5/3$.)

Define $\reg_{i} : = 2\N \cap \big[ 2^{i}, 2^{i+1} \big]$. 

\begin{lemma}\label{l.polyepsregthet}
Let $\varphi > 0$. 
For any $i \in \N$ sufficiently high, and $n \in 2\N \cap \big[ 2^{i+4},2^{i+5} \big]$,
$$
\psap_n \bigg( \Gamma \in \bigcup_{k \in \reg_{i+2}} \dubbub_{k,n-16-k}  \bigg) \geq n^{-\chi - \varphi} \, .
$$
\end{lemma}

 The set  $\reg_{i}$ may be thought as a set of {\em regular} indices $j \in 2\N \cap \big[ 2^{i}, 2^{i+1} \big]$, with~$\thet_j$ for example being neither atypically high nor low. All indices are regular as we prove Theorem~\ref{t.thetexist}(2); the set will be respecified non-trivially in the proof of Theorem~\ref{t.thetexist}(1).


(As we explained in the three-step outline, the snake method parameters, including $n$, will be set for the proof of Proposition~\ref{p.thetexist} at the start of step three. Until then, we are treating $n$ as a free variable.)

\medskip

\noindent{\bf Proof of Lemma~\ref{l.polyepsregthet}.} Note that the probability in question is equal to the ratio of the cardinality of the union $\bigcup_{j \in \reg_{i+2}} \dubbub_{j,n-16-j}$ and the polygon number $p_n$.

In order to bound below the numerator in this ratio, we aim to apply Proposition~\ref{p.polyjoin.aboveonehalf} with the role of $i \in \N$ played by $i+2$ and with $\mathsf{R} = \reg_{i+2}$.
In order to do so, set $\overline{\theta} = \sup_{n \in 2\N} \theta_n$, a quantity that is finite by hypothesis. 
Set the proposition's $\Theta$ equal to  $\overline{\theta}$. 
The index $\kmac$ in the proposition may then be any element of $\reg_{i+2}$.
Applying the proposition, we find that there exists $c = c(\Theta) > 0$ such that
$$
 \Big\vert \, \bigcup_{j \in \reg_{i+2}} \dubbub_{j,n-16-j} \, \Big\vert \, \geq \, c \,
      \frac{n^{1/2}}{\log n} \,
\sum_{j \in \reg_{i+2}} p_j p_{n-16-j} \, .
$$

The existence of~$\thet$
implies that $p_m \mu^{-m} \in \big[ m^{-\thet - \varphi} , m^{-\thet + \varphi} \big]$ whenever $m \geq m_0(\varphi)$.
Since $i \in \N$ is supposed sufficiently high, we thus find that
\begin{eqnarray*}
& & \psap_n \Big(  \Gamma \in \bigcup_{j \in \reg_{i+2}} \dubbub_{j,n-16-j}   \Big) \\
 & \geq & c \, \frac{n^{1/2}}{\log n}  \sum_{j \in \reg_{i+2}} j^{-(\thet + \varphi)} (n-16-j)^{-(\thet + \varphi)} \cdot n^{\theta - \varphi} \geq 2^{-4} c \, \frac{n^{3/2 - \thet - 3\varphi}}{\log n} \, ,
\end{eqnarray*}
where we used $\# \reg_{i+2} \geq 2^{i+1} \geq 2^{-4}n$. Noting that $\thet \leq 3/2 + \chi$ and relabelling~$\varphi > 0$ completes the proof. \qed

\subsection{Similarity of measure}


We may now state our similarity of measure Proposition~\ref{p.similarity}. Recall from Definition~\ref{d.rjsj} the constant difference sequence $(r_j,s_j)$, $1 \leq j \leq \indexc$, and  also the regular index set $\reg_{i+2}$ in the preceding subsection.
\begin{proposition}\label{p.similarity}
For each $j \in [1,L]$, there is a subset $D_j \subseteq \fpart_{r_j}$ 
satisfying
\begin{equation}\label{e.similarity}
 \psap_{s_j} \Big( \Gamma_{[0,r_j]} \not\in D_j \, , \, \Gamma \in \dubbub_{k,s_j - 16 - k} \, \, \, \textrm{for some $k \in \reg_{i+2}$} \Big) \leq 64 \mu^{16}  s_j^{-10} 
\end{equation}
such that, for each $\phi \in D_j$,
\begin{eqnarray*}
 & & \psap_{s_{\indexc}} \Big( \Gamma_{[0,r_j]} = \phi  \Big) \\
 & \geq & \tfrac{1}{40} \,  \conindex^{-1} \big( \log s_\indexc \big)^{-1} \big( \tfrac{3}{16} \big)^{3/2 + \chi} \, s_{\indexc}^{-4\eps}  \\
 & & \qquad \qquad \qquad \times \, \, \, 
  \psap_{s_j} \Big(  \Gamma_{[0,r_j]} = \phi \, , \,  \Gamma \in \dubbub_{k,s_j - 16 - k} \, \, \textrm{for some $k \in \reg_{i+2}$} \Big) \, .
\end{eqnarray*}
\end{proposition}

Preparing for the proof, we define, for each $\gamma \in \sap_{s_j}$, $1 \leq j \leq s_\indexc$, the set of $\ga$'s regulation global join indices 
$$
\tau_\gamma = \Big\{ k \in \reg_{i+2} : \gamma \in \dubbub_{k,s_{\indexcj} - 16 - k} \Big\} \, .
$$

 For $\gamma \in \sap_{s_{\indexcj}}$, the map sending $k \in \tau_\gamma$ to the junction plaquette associated to the join of elements of $\sapell_k$ and $\sapr_{s_\indexcj - 16 - k}$
is an injective map into $\globaljoin_\ga$: see the first assertion of Proposition~\ref{p.rgjpnumber}. Thus,  
$\vert \tau_\gamma \vert \leq \big\vert \globaljoin_\gamma \big\vert$.  Corollary~\ref{c.globaljoin} implies that   there exists a constant $\conindex > 0$ such that
\begin{equation}\label{e.constl}
 \psap_{s_{\indexc}} \Big( \big\vert \tau_\Gamma \big\vert \geq \conindex \log s_{\indexc} \Big) \leq \conindex \, s_{\indexc}^{-4\overline\theta - 22} \, ,
\end{equation}
where here we have reintroduced $\overline\theta = \sup_{n \in 2\N} \theta_n$, a quantity which we recall is finite by hypothesis; (this finiteness permits the application of the corollary).

We make use of this constant to specify for each $j \in [1,\indexc]$ a
set $D_j \subseteq \fpart_{r_j}$ of length~$r_j$ first parts that, when extended to form a length~$s_\indexc$ polygon, do not typically produce far-above-average numbers of regulation joinings. We define
$$
D_j = \Big\{ \phi \in \fpart_{r_j}: \psap_{s_{\indexc}} \Big( \big\vert \tau_\Gamma \big\vert \leq \conindex \log s_{\indexc} \, \Big\vert \, \Gamma_{[0,r_j]} = \phi    \Big) \geq 1 - \conindex \, s_{\indexc}^{-2 \overline\theta - 11} \Big\} \, ,
$$

The use of the next lemma lies in its final part, which establishes the $\psap_{s_j}$-rarity
assertion~(\ref{e.similarity}).
The first two parts of the lemma are steps towards the third.
The rarity under  $\psap_{s_\indexc}$ of the length $r_j$ initial path failing to be in $D_j$  is a direct consequence of definitions, as the first part of the lemma shows, but a little work is needed to make the change of index $\indexc \to j$.

\begin{lemma}\label{l.dcub}
Let $j \in [1,\indexc]$. 
\begin{enumerate}
\item We have that 
$$
\psap_{s_{\indexc}} \Big( \Gamma_{[0,r_j]} \in D_j^c \Big) \leq  s_{\indexc}^{-2 \overline\theta - 11} \, .
$$
\item
Let $k \in \reg_{i+2}$. Then
$$
\psapleft k \Big( \Gamma_{[0,r_j]} \not\in D_j \Big)
\leq   16 \mu^{16}  s_\indexc^{-11} \, .
$$
\item We have that
$$
 \psap_{s_j} \Big( \Gamma_{[0,r_j]} \not\in D_j \, , \, \Gamma \in \dubbub_{k,s_j - 16 - k} \, \, \, \textrm{for some $k \in \reg_{i+2}$} \Big) \leq 64 \mu^{16}  s_j^{- 10} \, .
$$
\end{enumerate}
\end{lemma}
\noindent{\bf Proof (1).}
Note that
\begin{eqnarray*}
 & & \psap_{s_{\indexc}} \Big( \big\vert \tau_\Gamma \big\vert \geq \conindex \log s_{\indexc} \Big) \geq \psap_{s_{\indexc}} \Big( \big\vert \tau_\Gamma \big\vert \geq \conindex \log s_{\indexc} \, , \, \Gamma_{[0,r_j]} \in D_j^c \Big) \\
  & =  &   \psap_{s_{\indexc}} \Big(\Gamma_{[0,r_j]} \in D_j^c \Big)  \psap_{s_{\indexc}} \Big( \big\vert \tau_\Gamma \big\vert \geq \conindex \log s_{\indexc} \, \Big\vert \, \Gamma_{[0,r_j]} \in D_j^c \Big) \\
 &  \geq & \conindex s_{\indexc}^{-2 \overline\theta - 11} \,  \psap_{s_{\indexc}} \Big(\Gamma_{[0,r_j]} \in D_j^c \Big)  \, . 
\end{eqnarray*}
Lemma~\ref{l.dcub}(1) thus follows from (\ref{e.constl}).  

\medskip

\noindent{\bf (2).}
In light of the preceding part, it suffices to prove 
\begin{equation}\label{e.sjsl}
\psapleft{k} \Big( \Gamma_{[0,r_j]} \not\in D_j \Big)
\leq  16 \mu^{16} s_\indexc^{2\overline{\theta}} \,
\psap_{s_\indexc} \Big( \Gamma_{[0,r_j]} \not\in D_j \Big) \, .
\end{equation}
We will derive this bound with a simple joining argument. 
Note that $s_\indexc \geq k + 16$ since $s_\indexc \geq 2^{i+4}$, $k \leq 2^{i+3}$ and $i \geq 1$.
Consider the map 
$$
\Psi:  \big\{ \phi \in \sapell_k : \phi_{[0,r_j]} \in D_j^c \big\} \times  \sapr_{s_\indexc - k - 16} \to  \big\{ \phi \in \sap_{s_\indexc} : \phi_{[0,r_j]} \in D_j^c \big\}
$$
defined by setting $\Psi(\phi_1,\phi_2)$ equal to the join polygon $J(\phi_1,\phi_2 + \vec{u})$, where the second argument is chosen so that the pair is Madras joinable, with the second polygon lying strictly below the $x$-axis.  
Lemma~\ref{l.leftlongbasic} and $r_j \leq s_j/2$ imply that
   $\Psi ( \phi , \phi' )_{[0,r_j]} =  \phi_{[0,r_j]}$, so that $\Psi$'s image lies in the stated set. The decomposition uniqueness in Proposition~\ref{p.rgjpnumber} shows that $\Psi$ is injective.
   Thus, by considering $\Psi$, we find that
$$
\# \big\{ \phi \in \sap_{s_\indexc} : \phi_{[0,r_j]} \in D_j^c \big\} \geq 
\# \big\{ \phi \in \sapell_k : \phi_{[0,r_j]} \in D_j^c \big\}  \# \big\{ \phi \in \sapr_{s_\indexc - k - 16} \big\} \, .  
$$   
By Lemma~\ref{l.polysetbound},
$$
 \psap_{s_\indexc} \big( \Gamma_{[0,r_j]} \in D_j^c \big)
  \geq \tfrac{1}{16} \psapleft{k} \big( \Gamma_{[0,r_j]} \in D_j^c \big) p_{s_\indexc}^{-1} p_k p_{s_\indexc - k - 16} \, .
$$  
Note that $p_{s_\indexc}^{-1} p_k p_{s_\indexc - k - 16} \geq \mu^{-16} s_\indexc^{-2\overline{\thet}}$, where recall that $\overline{\thet} = \sup_{n \in 2\N} \thet_n$. In this way, we obtain~(\ref{e.sjsl}). 

\medskip

\noindent{\bf (3).} 
 We first make the key observation that,
for $k \in \reg_{i+2}$ and $j' \in [1,\indexc]$,
\begin{equation}\label{e.keyob}
 \psap_{s_{j'}} \Big( \Gamma_{[0,r_j]} = \phi \, \Big\vert \, \Gamma \in \dubbub_{k,s_{j'} - 16 - k} \Big)  =  \psapleft{k} \Big( \Gamma_{[0,r_j]} = \phi \Big) \, .
\end{equation}
Proposition~\ref{p.leftlong} implies this, provided that its hypotheses that $r_j \leq k/2$ and $k/2 \leq s_{j'} - 16 - k \leq 35 k$ are valid. These bounds follow from $2^{i+3} \geq k \geq 2^{i+2}$, $2^{i+5} \geq s_{j'} \geq 2^{i+4}$ and $i \geq 2$.

The quantity in Lemma~\ref{l.dcub}(3) is at most
\begin{eqnarray*}
 & & \sum_{k \in \reg_{i+2}} \psap_{s_j} \Big( \Gamma_{[0,r_j]} \not\in D_j \, , \, \Gamma \in \dubbub_{k,s_j - 16 - k}  \Big)
\end{eqnarray*}
and thus by~(\ref{e.keyob}) with $j' = j$, at most
$\vert \reg_{i+2} \vert \max_{k \in \reg_{i+2}} 
\psapleft{k} \big( \Gamma_{[0,r_j]} \not\in D_j \big)$. 
We may now apply the first part of the lemma as well as the bounds
$\vert \reg_{i+2} \vert \leq 4 s_\indexc$ and $s_\indexc \geq s_j$ to obtain Lemma~\ref{l.dcub}(3). 
\qed

\begin{lemma}\label{l.compare}
Recall that $\indexc = 2^{i-4}$. 
Let $j \in [1,\indexc]$. 
\begin{enumerate}
\item 
For $\phi \in D_j$,
$$   
   \psap_{s_{\indexc}} \Big( \Gamma_{[0,r_j]} = \phi  \Big) \\
  \geq 
      \tfrac{1}{400} \, \conindex^{-1} \big( \log s_\indexc \big)^{-1} s_{\indexc}^{-2\eps} \mu^{-16} \sum_{k \in \reg_{i+2}} \psapleft{k} \Big( \Gamma_{[0,r_j]} = \phi  \Big)
   k^{1/2} p_k \mu^{-k} \, .
   $$ 
\item For $\phi \in \fpart_{r_j}$,
\begin{eqnarray*}
   & & \psap_{s_j} \Big(  \Gamma_{[0,r_j]} = \phi \, , \,  \Gamma \in \dubbub_{k,s_j - 16 - k} \, \, \textrm{for some $k \in \reg_{i+2}$} \Big) \\
   & \leq &  \tfrac{1}{10} \big( \tfrac{16}{3} \big)^{3/2 + \chi} \mu^{-16} \cdot s_j^{2\eps}  \sum_{k \in \reg_{i+2}} 
   \psapleft{k}  \Big(  \Gamma_{[0,r_j]} = \phi \Big)  k^{1/2} p_k \mu^{-k} \, .  
\end{eqnarray*}
\end{enumerate}
\end{lemma}
\noindent{\bf Proof:} {\em (1).}
Note that $\psap_{s_{\indexc}} \big( \Gamma_{[0,r_j]} = \phi  \big)$ is at least
\begin{eqnarray*}
 & &  \psap_{s_{\indexc}} \Big( \Gamma_{[0,r_j]} = \phi \, , \, \Gamma \in \bigcup_{k \in  \reg_{i+2}} \dubbub_{k,s_{\indexc} - 16 - k}  \, , \, \big\vert \tau_{\Gamma} \big\vert \leq \conindex \log s_\indexc \Big) \\
 & \geq & \big( \conindex \log s_\indexc \big)^{-1}
 \sum_{k \in  \reg_{i+2}}
 \psap_{s_{\indexc}} \Big( \Gamma_{[0,r_j]} = \phi \, , \, \Gamma \in \dubbub_{k,s_{\indexc} - 16 - k}  \, , \, \big\vert \tau_{\Gamma} \big\vert \leq \conindex \log s_\indexc \Big)  
  \, . 
\end{eqnarray*}
The summand in the last line may be written
\begin{eqnarray*}
 & & \psap_{s_{\indexc}} \Big( \Gamma_{[0,r_j]} = \phi \, , \, \Gamma \in \dubbub_{k,s_{\indexc} - 16 - k}  \Big)  \\
  & - &
 \psap_{s_{\indexc}} \Big( \Gamma_{[0,r_j]} = \phi \, , \, \Gamma \in \dubbub_{k,s_{\indexc} - 16 - k}  \, , \, \big\vert \tau_{\Gamma} \big\vert > \conindex \log s_\indexc \Big)  \, ,
\end{eqnarray*}
with the subtracted term here being at most
$$
  \psap_{s_{\indexc}} \Big( \Gamma_{[0,r_j]} = \phi \, , \, \big\vert \tau_{\Gamma} \big\vert > \conindex \log s_\indexc \Big) \leq   \conindex s_{\indexc}^{-2 \overline\theta - 11}  \psap_{s_{\indexc}} \Big( \Gamma_{[0,r_j]} = \phi \Big) 
$$
since $\phi \in D_j$.   
Note that $\vert \reg_{i+2} \vert \leq 4 s_\indexc$ since $\vert \reg_{i+2} \vert = 2^{i+1} + 1$ and $s_\indexc \geq 2^i$.
Thus,  $\psap_{s_{\indexc}} \big( \Gamma_{[0,r_j]} = \phi  \big)$ is at least
\begin{eqnarray}
 & & \big( \conindex \log s_\indexc \big)^{-1}
 \sum_{k \in  \reg_{i+2}} \psap_{s_{\indexc}} \Big( \Gamma_{[0,r_j]} = \phi \, , \, \Gamma \in \dubbub_{k,s_{\indexc} - 16 - k}  \Big) \label{e.conindexbound}  \\
 & & \qquad  \, - \, \, \,  
\big( \conindex \log s_\indexc \big)^{-1} 4 s_\indexc \cdot  \conindex s_{\indexc}^{-2 \overline\theta - 11} \psap_{s_{\indexc}} \Big( \Gamma_{[0,r_j]} = \phi \Big)  \, . \nonumber 
\end{eqnarray}

Taking $j' = \indexc$ in~(\ref{e.keyob}) and using Proposition~\ref{p.rgjpnumber}, we learn that the summand in the first line satisfies 
\begin{eqnarray}
 & = &  \psap_{s_{\indexc}} \Big( \Gamma_{[0,r_j]} = \phi \, , \, \Gamma \in \dubbub_{k,s_{\indexc} - 16 - k} \Big) \label{e.summandeq} \\
 & = &  \psapleft{k} \Big( \Gamma_{[0,r_j]} = \phi \Big) \psap_{s_{\indexc}} \Big( \Gamma \in \dubbub_{k,s_{\indexc} - 16 - k} \Big) \nonumber \\
 & = & 
\psapleft{k} \Big( \Gamma_{[0,r_j]} = \phi \Big) \cdot p^{-1}_{s_{\indexc}} \cdot \lfloor \tfrac{1}{10} \, k^{1/2}  \rfloor \cdot \big\vert \sapell_k \big\vert \, \big\vert \sapr_{s_{\indexc} - 16 - k} \big\vert \, , \nonumber
\end{eqnarray} 
and thus by Lemma~\ref{l.polysetbound} (as well as $k \geq 2^{i+2}$ and $i \geq 7$) is at least
$$ 
\psapleft{k} \Big( \Gamma_{[0,r_j]} = \phi \Big) \cdot p^{-1}_{s_{\indexc}} \cdot \tfrac{1}{20} \, k^{1/2}  \cdot
 \tfrac{1}{8} p_k \cdot \tfrac{1}{2} p_{s_{\indexc} - 16 - k} \, . 
$$

Returning to~(\ref{e.conindexbound}), we find that
  $\psap_{s_{\indexc}} \big( \Gamma_{[0,r_j]} = \phi  \big)$ is at least
\begin{eqnarray}
 & & \Big( 1 +  
4  \big( \log s_\indexc \big)^{-1}    s_{\indexc}^{-2 \overline\theta - 10} \Big)^{-1}  \big( \conindex \log s_\indexc \big)^{-1} \tfrac{1}{320} \, p^{-1}_{s_{\indexc}} \nonumber \\
& & \qquad \qquad \qquad \qquad \qquad \times \, \, \, \sum_{k \in \reg_{i+2}} \psapleft{k} \Big( \Gamma_{[0,r_j]} = \phi  \Big)
   k^{1/2} p_k p_{s_{\indexc} - 16 - k} \label{e.willbesame} \\
  & \geq &   \big( \conindex \log s_\indexc \big)^{-1} \tfrac{1}{400} \, p^{-1}_{s_{\indexc}} \mu^{s_\indexc - 16} \nonumber \\
& & \qquad \qquad \qquad \times \, \, \, \sum_{k \in \reg_{i+2}} \psapleft{k} \Big( \Gamma_{[0,r_j]} = \phi  \Big)
   k^{1/2} p_k \mu^{-k} (s_{\indexc} - 16 - k)^{-3/2 - \chi - \eps} \, . \nonumber
\end{eqnarray}
Using $p_{s_{\indexc}} \leq \mu^{s_{\indexc}} s_{\indexc}^{-3/2 - \chi + \eps}$, we obtain Lemma~\ref{l.compare}(1). 
   
\medskip

\noindent{{\em (2).}} For $\phi \in \fpart_{r_j}$,
\begin{eqnarray*}
 & &   \psap_{s_j} \Big(  \Gamma_{[0,r_j]} = \phi \, , \,  \Gamma \in \dubbub_{k,s_j - 16 - k} \, \, \textrm{for some $k \in \reg_{i+2}$} \Big) \\
 & \leq & \sum_{k \in \reg_{i+2}} \psap_{s_j}  \Big(  \Gamma_{[0,r_j]} = \phi \, , \,   \Gamma \in \dubbub_{k,s_j - 16 - k} \Big) \, .
 \end{eqnarray*}
By taking $j' = j$ in (\ref{e.keyob}), we see that the two equalities in the display~(\ref{e.summandeq}) hold with $s_j$ in place of $s_\indexc$. Since  $\big\vert \sapell_k \big\vert \leq p_k$ and $\big\vert \sapr_{s_j - 16 - k} \big\vert \leq p_{s_j - 16 - k}$, the last sum is seen to be at most
 \begin{eqnarray}
&  &  \tfrac{1}{10} \, p_{s_j}^{-1} \sum_{k \in \reg_{i+2}} 
   \psapleft{k}  \Big(  \Gamma_{[0,r_j]} = \phi \Big)  k^{1/2} p_k p_{s_j - 16 - k} \label{e.willbesametwo} \\
 & \leq &  \tfrac{1}{10} s_j^{3/2 + \chi + \eps} \mu^{-16} \sum_{k \in \reg_{i+2}} 
   \psapleft{k}  \Big(  \Gamma_{[0,r_j]} = \phi \Big)  k^{1/2} p_k \mu^{-k} (s_j - 16 - k)^{-3/2 - \chi + \eps} \nonumber \\
 & \leq &  \tfrac{1}{10}  \big( \tfrac{16}{3} \big)^{3/2 + \chi} \mu^{-16}  \cdot s_j^{2\eps}  \sum_{k \in \reg_{i+2}} 
   \psapleft{k}  \Big(  \Gamma_{[0,r_j]} = \phi \Big)  k^{1/2} p_k  \mu^{-k}  \, , \nonumber   
 \end{eqnarray}
where we used $2^{i+5} \geq s_j \geq 2^{i+4}$, $k \leq 2^{i+3}$ and $i \geq 3$
in the form $s_j - 16 - k \geq 3s_j/16$.    \qed
  
\medskip

\noindent{\bf Proof of Proposition~\ref{p.similarity}.}
A consequence of Lemmas~\ref{l.dcub}(3) and~\ref{l.compare}. \qed

\subsection{Obtaining Proposition~\ref{p.almost} via similarity of measure}
This subsection is devoted to the next proof.
\medskip

\noindent{\bf Proof of Proposition~\ref{p.almost}.}
Note that the right-hand side of the bound
$$
  \psap_{s_{\indexc}} \Big( \Gamma_{[0,r_j]} \in \phicl{r_j}{s_j -1}{1/2 + 2\chi + \eclo} \Big) 
  \geq  
  \psap_{s_{\indexc}} \Big( \Gamma_{[0,r_j]} \in \phicl{r_j}{s_j -1}{1/2 + 2\chi + \eclo} 
  \, , \,  
\Gamma_{[0,r_j]} \in D_j  \Big)  
$$
is, by the similarity of measure Proposition~\ref{p.similarity}, at least the product of the quantity $c  \big( \log s_\indexc \big)^{-1} s_{\indexc}^{-4\eps}$ and the probability
$$
 \psap_{s_j} \Big( \Gamma_{[0,r_j]} \in \phicl{r_j}{s_j -1}{1/2 + 2\chi + \eclo} \, , \,  \Gamma \in \dubbub_{k,s_j - 16 - k} \, \, \textrm{for some $k \in \reg_{i+2}$} \, , \,  
\Gamma_{[0,r_j]} \in D_j   \Big) \, ,
$$
where the constant $c$ equals $\tfrac{1}{40}  \conindex^{-1} \big( \tfrac{3}{16} \big)^{3/2 + \chi}$.
Note further that
\begin{itemize}
\item
$\psap_{s_j} \Big( \Gamma_{[0,r_j]} \in \phicl{r_j}{s_j -1}{1/2 + 2\chi + \eclo} \Big) \geq  1 - s_j^{-(\chi + \eclo - \eroom)}$, 
by (\ref{e.acondition}) with $a=1$,
\item
and $\psap_{s_j} \Big(   \Gamma \notin \dubbub_{k,s_j - 16 - k} \, \, \textrm{for any $k \in \reg_{i+2}$}   \Big) \leq 1 - s_j^{-(\chi + \eclo/2 )}$, by Lemma~\ref{l.polyepsregthet} with $\varphi = \eclo/2$;
\item and recall~(\ref{e.similarity}).
\end{itemize}
Thus,
\begin{eqnarray*}
 & & \psap_{s_{\indexc}} \Big( \Gamma_{[0,r_j]} \in \phicl{r_j}{s_j -1}{1/2 + 2\chi + \eclo} \Big) \\
 & \geq & c \, \big( \log s_\indexc \big)^{-1} s_{\indexc}^{-4\eps}  \, \bigg(  s_j^{-(\chi + \eclo/2 )}  - s_j^{-(\chi + \eclo - \eroom)} \, - \, 64 \mu^{16} s_{\indexcj}^{-10}  \bigg) \, . 
\end{eqnarray*}

Provided that the index $i \in \N$ is supposed to be sufficiently high, the right-hand side is at least $\tfrac{1}{2} c  \big( \log s_\indexc \big)^{-1} s_{\indexc}^{-\chi - \eclo/2 - 4\eps}$ if we insist that $\eclo > 2\eroom$ as well as $\chi + \eclo < 10$ (the latter following from $\chi < 1/6$ and the harmlessly assumed $\eclo < 1$).

We are now ready to obtain Proposition~\ref{p.almost}.
We must set the values of the quantities $m$ and $m'$ in the proposition's statement. We take $m = s_{\indexc}$
and $m'$ equal to the constant difference $s_j - r_j$
between the terms in any pair in the sequence constructed in Definition~\ref{d.rjsj}. Recalling the bound on this difference stated in the definition, we see that 
 $2^{i+3} \leq m' \leq 2^{i+4} \leq m \leq 2^{i+5}$. 
 Set $\eclo = 4\eroom$ and $\eroom = \eps$, so that  $\chi + \eclo/2 + 4\eps = \chi + 6\eps < \esm - \eps$ and $1/2 + 2\chi + \eclo = \alpha$.  
 Note that $r_j$ belongs to the set $K$ specified in Proposition~\ref{p.almost} whenever $j \in [1,\indexc]$, because $r_j = s_j - m' \leq m - m'$. Thus, $\# K \geq 2^{i-4} \geq 2^{-9} m$, and 
Proposition~\ref{p.almost} is proved. \qed

\medskip

\noindent{\em A concluding remark for step two.} We review our approach in the context of Figure~\ref{f.snakemethodpolygon}. Note that, in the three bullet points above, the second line asserts that regulation global join polygons, while not necessarily typical, are not so rare under~$\psap_{s_j}$. The first line asserts that the initial subpath $\Ga_{[0,r_j]}$ is highly likely to have typical conditional closing probability when a length $s_j - r_j$ extension is considered. Since $s_j - r_j$ is independent of $j$, we have labelled it, calling it $m'$, thereby precisely specifying this quantity from the explanation that accompanies Figure~\ref{f.snakemethodpolygon}. The first bullet point, similarly to Proposition~\ref{p.notquite}, corresponds to the informal claim in the figure's caption that black dots are very typically present at $\tilden$ steps from the end of a polygon. The extra margin in the choice of $\alpha = 1/2 + 2\chi + o(1)$ that we alluded to the start of Section~\ref{s.notquite}  ensures that such black dots are so highly typical that they remain typical even among regulation polygons; i.e., the error in the first bullet point is smaller than in the second. For these regulation polygons,
 we have implemented comparison of measure in Proposition~\ref{p.similarity}. Pursuing the caption's story: black dots make the jump as length changes from $s_j$ to $s_\indexc$ with probability~$s_{\indexc}^{-4\eps}$ and are thus proved in Proposition~\ref{p.almost} to have a non-negligible probability of appearing in generic locations in a uniform polygon of length $m = s_\indexc$.

\section{Step three: completing the proof of Proposition~\ref{p.thetexist}}
As we explained in this step's outline, we set the snake method index parameters
at the start of step three.  For a given sufficiently high  choice of $i \in \N$,  Proposition~\ref{p.almost} specifies  $m$ and $m'$. We then set the two index parameters, with $n$ equal to  $m-1$ and $\ell$ equal to $m - m'$.

For $\gamma \in \sap_{n+1}$, set 
$$
X_\gamma =  \sum_{k=0}^\ell  1\!\!1_{\gamma_{[0,k]} \in \phicl{k}{m'+k-1}{\alpha}} \, .
$$

Recall that $\gamma_{[0,\ell]} \in \fpart_{\ell,n}$ is $(\alpha,n,\ell)$-charming at an index $k \in [1,\ell]$ if
$$
 \saw^0_{k + n - \ell} \Big( \Ga \closes \, \Big\vert \, \big\vert \Ga^1 \big\vert = k \, , \, \Ga^1 = \ga_{[0,k]} \Big) > n^{-\alpha} \, ;
$$
on the other hand, for given $k \in [1,\ell]$, the property that such a $\gamma$ satisfies $\ga_{[0,k]} \in \hPhi_{k,n-\ell+k}^\alpha$ takes the same form with $n^{-\alpha}$ replaced by the slightly larger quantity $(n - \ell + k)^{-\alpha}$. Since $n - \ell = m' - 1$, we find that, for any $\gamma \in \sap_{n+1}$,  
\begin{equation}\label{e.relate}
 X_\gamma \geq n^{\beta - \esm}/4 \, \, \, \textrm{implies that} \, \, \gamma_{[0,\ell]} \in \mathsf{CS}_{\beta,\eps}^{\alpha,\ell,n} \, .
\end{equation}

We now show that
\begin{equation}\label{e.claimx}
 \psap_{n+1} \Big( X_\Gamma \geq   n^{1- \esm + \eps/2}   \Big) \geq    n^{-\esm + \eps/2} \, .
\end{equation}

To derive this,  consider the expression
$$
S = \sum_{\gamma \in \sap_{n+1}}  \sum_{k=0}^{\ell}
1\!\!1_{\gamma_{[0,k]} \in \phicl{k}{m'+k}{\alpha}} \, .
$$

Recall that $p_{n+1}$ denotes $\# \sap_{n+1}$; using Proposition~\ref{p.almost} in light of~$n + 1 = m$,
\begin{eqnarray*}
  S  & = & 
 p_{n+1} \, \sum_{k=0}^{\ell}  \psap_{n+1} \Big(  \Gamma_{[0,k]} \in \phicl{k}{m'+k}{\alpha}   \Big) \\
 & \geq &   p_{n+1} \cdot \# K \cdot (n+1)^{-\esm + \eps} \geq  p_{n+1} \cdot 2^{-9} n \cdot  2^{-1} n^{-\esm + \eps}  \, .  
\end{eqnarray*}

Let $q$ denote the left-hand side of~(\ref{e.claimx}).
 Note that
$$
  S  \leq  p_n   \cdot \Big( q \big( \ell + 1 \big) + (1-q) n^{1-\esm + \eps/2}   \Big) \, .
$$
From the lower bound on $S$, and $n \geq \ell$,
$$
 q  \big( n + 1 \big) +  n^{1-\esm + \eps/2}
 \geq   2^{-10}  \, n^{1- \esm + \eps} \, ,
$$
which implies that $q \geq  n^{-\esm + \eps/2}$; in this way, we obtain~(\ref{e.claimx}).

It is unsurprising that $n^{1-\esm + \eps/2} \geq n^{1 - \esm}/4$ (for all $n \in \N$, including our choice of $n$). Since the snake length exponent $\beta$ is set to one, we learn from~(\ref{e.relate}) and~(\ref{e.claimx}) that
$$
  \psap_{n+1} \Big( \Gamma_{[0,\ell]} \in \mathsf{CS}_{\beta,\eps}^{\alpha,\ell,n}   \Big) \geq \psap_{n+1} \Big( X_\Gamma \geq   n^{1-\esm + \eps/2}   \Big) \geq    n^{-\esm + \eps/2} \, .
 $$
Thus, if the dyadic scale parameter $i \in \N$ is chosen so that $n \geq 2^{i+4}$ is sufficiently high,  the charming snake presence hypothesis~(\ref{e.afewcs}) is satisfied. By Theorem~\ref{t.snakesecondstep}, $\psaw_n \big( \Gamma \closes \big) \leq 2(n+1)c^{-n^{\delta}/2}$. This deduction has been made for some value of $n \in (2\N + 1) \cap \big[2^{i+4}-1,2^{i+5}-1]$, where here $i \in \N$ is arbitrary after a finite interval. Relabelling $c > 1$ to be any value in $(1,c^{1/2})$
completes the proof of Proposition~\ref{p.thetexist}. \qed

\section{Preparing for the proof of Theorem~\ref{t.thetexist}(1)}

This argument certainly fits the template offered by the proof of Theorem~\ref{t.thetexist}(2). We will present the proof by explaining how to modify the statements and arguments used in the preceding derivation.

Recall first of all that the hypotheses of Theorem~\ref{t.thetexist}(2) implied the existence of the closing exponent, and that we chose to label $\chi \in \R$ such that $\psap_n\big( \Gamma \closes \big) = n^{-1/2 - \chi + o(1)}$ for odd~$n$;
 moreover, $\chi$ could be supposed to be non-negative. 
 Specifically, it was the lower bound  $\psap_n\big( \Gamma \closes \big) \geq n^{-1/2 - \chi - o(1)}$ that was invoked, in part one of the preceding proof.  

In the present proof, the closing exponent is not hypothesised to exist, and so we must abandon this usage of $\chi$. However, this parameter will be used again, with the basic role of~$1/2 + \chi$ as the exponent  in such a closing probability lower bound  being maintained.

\begin{definition}\label{d.highclose}
For $\maceta > 0$, define the set of indices $\highclose_\maceta \subseteq 2\N$ of $\maceta$-{\em high closing probability},
$$
 \highclose_\maceta = \Big\{ n \in 2\N : \psaw_{n-1} \big( \Gamma \closes \big) \geq n^{- \maceta} \Big\} \, .
$$
\end{definition}

For $\chi > 0$ arbitrary, we introduce the {\em closing probability}

\medskip

\noindent{\bf Hypothesis  $\mathsf{CP}_{\chi}$.} The set $2\N \setminus \hcb$ has limit supremum density in $2\N$ less than $1/1250$. 

\medskip

We begin the proof by the same type of reduction as was used in the earlier derivation. Here is the counterpart to  Proposition~\ref{p.thetexist}.
\begin{proposition}\label{p.closingprob}
Let $d = 2$. Let $\chi \in (0,1/{14})$, and assume Hypothesis~$\mathsf{CP}_\chi$. For some $c > 1$ and $\delta >0$, the set of odd integers $n$ satisfying $n + 1 \in \highclose_{1/2 + \chi}$ and 
$$
\psaw_n \big( \Ga \closes \big) \leq c^{-n^{\delta}}
$$
intersects the dyadic scale $\big[ 2^i,2^{i+1} \big]$ for all but finitely many $i \in \N$.
\end{proposition}
\noindent{\bf Proof of Theorem~\ref{t.thetexist}(1).}
The two properties of the index $n$ asserted by the proposition's conclusion are evidently in  contradiction for $n$ sufficiently high. Thus  Hypothesis~$\mathsf{CP}_\chi$ is false whenever  $\chi \in [0,1/{14})$. This implies the result. \qed

\medskip

The rest of the article is dedicated to proving Proposition~\ref{p.closingprob}. The three step plan of attack for deriving Proposition~\ref{p.thetexist} will also be adopted for the new proof.
Of course some changes are needed. 
We must cope with a deterioration in the known regularity of the $\theta$-sequence.
We begin by overviewing 
what we know in the present case and how the new information will cause changes in the three step plan.

In the proof of Theorem~\ref{t.thetexist}(2) via Proposition~\ref{p.thetexist}, we had the luxury of 
assuming the existence of the closing exponent. Throughout the proof of Proposition~\ref{p.closingprob} (and thus henceforth), we instead {\em fix $\chi \in [0,1/{14} )$, and suppose that Hypothesis~$\mathsf{CP}_\chi$ holds}. Thus, our given information is that the closing probability is at least $n^{-1/2 - \chi}$ for  a uniformly positive proportion of indices $n$ in any sufficiently long initial interval of positive integers.

We may state and prove right away the $\thet$-sequence regularity offered by this information.

\begin{proposition}\label{p.thetint}
Let $\eps \in (0,1)$. For all but finitely many values of $i \in \N$,
$$
 \# \Big\{ j \in 2\N \, \cap \, \big[2^{i},2^{i+1} \big]:  3/2  - \eps \leq \thet_ j \leq 3/2 + \chi + \eps \Big\}  \geq 2^{i-1} \big( 1 - \tfrac{1}{300} \big) \, .
$$
\end{proposition}
\noindent{\bf Proof.}
Theorem~\ref{t.polydev} implies that,  for all but finitely many values of $i \in \N$,
\begin{equation}\label{e.thirtyone}
 \# \Big\{ j \in 2\N \, \cap \, \big[2^{i},2^{i+1} \big]: \thet_ j \leq  3/2   - \eps \Big\} \leq   \tfrac{1}{600} \,  2^{i-1}  \, .
\end{equation}
Thus, it is enough to show that, also for all but finitely many $i$,
\begin{equation}\label{e.thirtytwo}
 \# \Big\{ j \in  2\N \, \cap \, \big[2^{i},2^{i+1} \big]: \thet_j > 3/2 + \chi + \eps \Big\} \leq   \tfrac{1}{600} \,  2^{i-1}  \, . 
\end{equation}

We will invoke Hypothesis $\mathsf{CP}_{\chi}$ as we verify~(\ref{e.thirtytwo}). 
First, we make a 

\medskip

\noindent{\bf Claim.}
For any
 $\varphi > 0$, the set  $\hcb  \cap   \big\{ n \in 2 \N: \thet_ n > 3/2 + \chi + \varphi \big\}$ is finite.

To see this, note that $2n p_n/c_{n-1} \geq n^{-1/2 - \chi}$ if $n \in \hcb$. From $c_{n-1} \geq \mu^{n-1}$ and $p_n = \mu^n n^{-\thet_n}$ follows $n^{-\thet_n} \geq \tfrac{1}{2\mu} n^{-3/2 - \chi}$ for such~$n$, and thus the claim.

The claim implies that, 
under Hypothesis $\mathsf{CP}_{\chi}$, the set
$\big\{ n \in 2\N : \thet_ n > 3/2 + \chi + \eps \big\}$
has limit supremum density  at most $1/1250$. From this,~(\ref{e.thirtytwo}) follows directly. \qed


\medskip

That is, the weaker regularity information forces a positive proportion of the polygon number deficit exponents $\thet_n$ to lie in any given open set containing the interval~$[3/2,3/2 + \chi]$. 
In the preceding proof, we knew this for the one-point 
set $\{ 3/2 + \chi \}$, for all high~$n$. 

In the three step plan, we will adjust step one so that the constant difference sequence constructed there incorporates the $\theta$-regularity that Proposition~\ref{p.thetint} shows to be typical, as well as a closing probability lower bound permitted by Hypothesis~$\mathsf{CP}_{\chi}$.

In step two, the weaker regularity will lead to a counterpart Lemma~\ref{l.polyepsreg} to the regulation polygon ample supply Lemma~\ref{l.polyepsregthet}. Where before we found a lower bound on the probability that a polygon is regulation global join of the form $n^{-\chi - o(1)}$, now we will find only a bound $n^{-2\chi - o(1)}$. 

Continuing this step, the mechanism of measure comparison for initial subpaths of polygons drawn from the laws $\psap_n$ for differing lengths~$n$ is no longer made via all regulation global join polygons but rather via such polygons whose length index lies in a certain regular set. The new similarity of measure result counterpart to Proposition~\ref{p.similarity} will be Proposition~\ref{p.similaritynew}. It will be applied at the end of step two to prove a slight variant of Proposition~\ref{p.almost}, namely Proposition~\ref{p.almostnew}.

This last result will yield Proposition~\ref{p.closingprob}
in step three by a verbatim argument to that by which Proposition~\ref{p.thetexist} followed from Proposition~\ref{p.almost}.

 \section{The three steps for Proposition~\ref{p.closingprob}'s proof in detail}
 
In a first subsection of this section, we reset the snake method's parameters to handle the weaker information available.
In the following three subsections, we state and prove the assertions associated to each of the three steps.

\subsection{The snake method exponent parameters}
Since $\chi < 1/{14}$, we may fix a parameter  $\eps \in \big(0,(1/2 - 7\chi)/14 \big)$, and do so henceforth. The three exponent parameters are then set so that
\begin{itemize}
\item $\beta = 1$;
\item $\alpha = 1/2 + 3\chi + 5\eps$; 
\item and $\esm = 4\chi + 9\eps$.
\end{itemize}
Note that the quantity $\delta = \beta - \esm - \alpha$, which must be positive if the method to work, is equal to $1/2 - 7\chi - 14\eps$.
The constraint imposed on $\eps$ ensures this positivity.

\subsection{Step one}

The essential conclusion of the first step one was the construction of the constant difference sequence in Definition~\ref{d.rjsj}. Now our aim is similar.

Fixing as previously two parameters $\eclo > \eroom > 0$ (which in the present case we will ultimately set $\eclo = 5\eroom$ and $\eroom = \eps$ in terms of our fixed parameter $\eps$), we adopt the Definition~\ref{d.cpt} of a closing probability typical index pair on dyadic scale~$i$, where of course the parameter $\chi$ is fixed in the way we have explained. An index pair $(k,j)$ that satisfies this definition is called {\em regularity typical} for dyadic scale~$i$ if moreover  $j \in \highclose_{1/2 + \chi}$ and  
\begin{equation}\label{e.tineqs}
 3/2  - \epn  \leq \thet_j \leq 3/2 + \chi + \epn  \, .
\end{equation}

We then replace Definition~\ref{d.rjsj} as follows.
\begin{definition}\label{d.rjsj.new}
Set $\indexc = 2^{i-4}$. We specify a {\em constant difference} sequence of  index pairs $(r_j,s_j)$, $1 \leq j \leq \indexc$, whose elements are regularity typical for dyadic scale $i$, with $s_j - r_j$ independent of $j$ and valued in $[2^{i+4}- 2^i,2^{i+4}]$,  and $\big\{ r_j: 1 \leq j \leq \indexc \big\}$ increasing.
\end{definition}
Our job in step one is to construct such a sequence.
First, recall the set $E$ of index pairs $(k,j)$ specified in Definition~\ref{d.e}. Let $E'$ be the subset of $E$ consisting of such pairs for which~$j \in \highclose_{1/2 + \chi}$ and~(\ref{e.tineqs}) holds.
Counterpart to Lemma~\ref{l.etypical}, we have:
\begin{lemma}\label{l.eprimetypical}
Provided the index $i \in \N$ is sufficiently high,
for at most $\tfrac{1}{100} 2^{i+3}$ values of $j \in  2 \N \cap \big[ 2^{i+4},  2^{i+5} \big]$, the set of $k \in \big[ 2^i , 2^i + 2^{i-2} \big]$ such that $(k,j) \not\in E'$ has cardinality exceeding $2 j^{1 - \eroom/2} \leq 2 \big(  2^{i+5} \big)^{1-\eroom/2}$.
\end{lemma}
\noindent{\bf Proof.}
Proposition~\ref{p.thetint} implies that the set of $j \in 2\N \cap \big[ 2^{i+4},2^{i+5} \big]$ satisfying~(\ref{e.tineqs}) -- a set containing all second coordinates of pairs in $E'$ -- has cardinality at least $\big( 1 -\tfrac{1}{300} \big) 2^{i+3}$. Invoking Hypothesis~$\mathsf{CP}_\chi$, we see that, of these values of~$j$, at most~$\tfrac{1}{9600} 2^{i+3}$ fail the test of membership of~$\highclose_{1/2 + \chi}$.
We claim that the set of remaining values of $j$ 
satisfies the statement in the lemma. Indeed, all these values
belong to $\highclose_{1/2 + \chi}$. As such, we are able to invoke Lemma~\ref{l.closecard}, as we did in proving Lemma~\ref{l.etypical}, but this time with $\alphamac = 1/2 + \chi$ and $\deltamac = \eroom$, to find that at most  $2 \big(  2^{i+5} \big)^{1-\eroom}$ values of  $k \in \big[ 2^i , 2^i + 2^{i-2} \big]$ violate the condition~(\ref{e.rchithet}). \qed

Analogously to Proposition~\ref{p.manyethet},
\begin{proposition}\label{p.manyethet.new}
There exists $k \in 2\N \cap \big[ 2^{i+4}, 2^{i+5} - 2^{i-2} \big]$ for which there are at least $2^{i-4}$ values of $j \in 2\N \cap [0,2^{i-2}]$ with the index pair
$\big( 2^i + j, k + j \big)$ regularity typical on dyadic scale $i$.
\end{proposition}
\noindent{\bf Proof.} Using Lemma~\ref{l.eprimetypical},
we find, analogously to~(\ref{e.ineqe}),
\begin{eqnarray*}
 & & \sum_{j \in 2\N \cap [ 2^{i+4} + 2^{i-2} , 2^{i+4} + 2^{i} ] } \sum_{\ell = 2^i}^{2^i + 2^{i-2}} 
 1\!\!1_{(\ell,j) \in E'} \\ 
 & \geq &  \Big( 2^{i-2} + 1 - 2 \big( 2^{i+5} \big)^{1 - \eroom} \Big) \cdot  \Big( \tfrac{1}{2} \big( 2^{i} - 2^{i-2} \big) - \tfrac{1}{100} 2^{i+3} \Big) \, . 
\end{eqnarray*}
The argument for Proposition~\ref{p.manyethet} now yields the result. \qed

\medskip

We may now conclude step one, because the proof of Proposition~\ref{p.notquite} is valid with the notion of regularity typical replacing closing probability typical; thus, we have specified a constant difference sequence in the sense of Definition~\ref{d.rjsj.new}. (Incidentally, the unused Proposition~\ref{p.notquite} is seen to be valid with right-hand side strengthened to $1 - 2^{-2i\chi}$ if we choose $\eclo = 5\eps$ and $\eclo > \eroom$ by considering $a=2$ in~(\ref{e.acondition}).)

\subsection{Step two}

The new version of Proposition~\ref{p.almost} differs from the original only in asserting that the constructed $m$ belongs to $\highclose_{1/2 + \chi}$.
\begin{proposition}\label{p.almostnew}
For each $i \in \N$ sufficiently high, there exist $m \in 2\N \cap \big[ 2^{i+4}, 2^{i+5} \big] \cap \highclose_{1/2 + \chi}$ and $m' \in \big[ 2^{i+3},2^{i+4} \big]$ such that, writing $K$ for the set of values $k \in \N$, $1 \leq k \leq m-m'$, that satisfy
$$
\psap_m \Big( \Ga_{[0,k]} \in \hPhi_{k,k+m'-1}^\alpha \Big)
 \geq  m^{-\esm + \eps} \, ,
$$
we have that $\vert K \vert \geq 2^{-9} m$.
\end{proposition}

In this step, we provide a suitable similarity of measure result, Proposition~\ref{p.similaritynew}, and use it to derive Proposition~\ref{p.almostnew}.

Before doing so, we must adapt our regulation polygon tools to cope with the weaker regularity that we are hypothesising. The ample regulation polygon supply Lemma~\ref{l.polyepsregthet}
will be replaced by Lemma~\ref{l.polyepsreg}. 

In preparing to state the new lemma, recall that Lemma~\ref{l.polyepsregthet} 
was stated in terms of a regular index set $\reg_{i+2}$ that was specified in a trivial way.  
Recalling  that $\indexc = 2^{i-4}$, we now specify a non-trivial counterpart set $\reg_{n,i+2}^{\chi + \eps}$ when $n$ is an element in the sequence $\big\{ s_j: 1 \leq j \leq \indexc \big\}$ from step one. 
\begin{definition}
For $j \in \big[1, \indexc \big]$,   define
$\regnew{s_j}{i+2}{\chi + \epn}$ to be the set of $k \in 2\N \cap \big[ 2^{i+2} , 2^{i+3} \big]$ such that 
$$
  \max \big\{ \thet_k , \thet_{s_j-16-k} ,  
\thet_{s_\indexc - 16 - k}  \big\}  \leq  3/2 +  \chi + \eps
\, \, \, \textrm{and} \, \, \,
\thet_{s_j - 16 - k} \geq 3/2 - \epn  \, \, .   
$$
\end{definition}

Here is the new ample supply lemma.
\begin{lemma}\label{l.polyepsreg}
Let $j \in [1,\indexc]$. 
Then
$$
\psap_{s_j} \Big( \Gamma \in \dubbub_{k,s_j-16-k} \, \, \textrm{for some $k \in  \regnew{s_j}{i+2}{\chi + \epn}$} \Big) \geq \frac{c}{\log s_j} \, {s_j}^{-(2\chi + 3\epn)} \, ,
$$
where $c > 0$ is a universal constant.
\end{lemma}

The similarity of measure Proposition~\ref{p.similarity} is replaced by the next result. 
\begin{proposition}\label{p.similaritynew}
For each $j \in [1,L]$, there is a subset $D_j \subseteq \fpart_{r_j}$ 
satisfying
\begin{equation}\label{e.similaritynew}
 \psap_{s_j} \Big( \Gamma_{[0,r_j]} \not\in D_j \, , \, \Gamma \in \dubbub_{k,s_j - 16 - k} \, \, \, \textrm{for some $k \in \regnew{s_j}{i+2}{\chi + \epnmac}$} \Big) \leq 64 \mu^{16}  s_j^{-10} 
\end{equation}
such that, for each $\phi \in D_j$,
\begin{eqnarray*}
  \psap_{s_{\indexc}} \Big( \Gamma_{[0,r_j]} = \phi  \Big) 
  & \geq  & \tfrac{1}{40} \,  \conindex^{-1} 
\big( \tfrac{3}{16} \big)^{3/2}  \big( \log s_\indexc \big)^{-1}  s_\indexc^{-2\chi -  4\epn} \\ 
 & &    \times \, \, \, \psap_{s_j} \Big(  \Gamma_{[0,r_j]} = \phi \, , \,  \Gamma \in \dubbub_{k,s_j - 16 - k} \, \, \textrm{for some $k \in \regnew{s_j}{i+2}{\chi + \epnmac}$} \Big) \, .
\end{eqnarray*}
\end{proposition}

In the rest of this subsection of step two proofs, we prove in turn Lemma~\ref{l.polyepsreg} and Propositions~\ref{p.similaritynew} and~\ref{p.almostnew}.

\subsubsection{Deriving Lemma~\ref{l.polyepsreg}}

We will use the next result.

\begin{lemma}\label{l.rsize}
For  $i \in \N$ with $i \geq 6$, and $n \in 2\N \cap \big[ 2^{i+4},2^{i+5} \big]$,
$$
  \big\vert \regnew{s_j}{i+2}{\chi + \epn}  \big\vert \geq  \big( \tfrac{9}{10} - \tfrac{1}{32} \big) 2^{i+1}    \, .
$$

\end{lemma}
\noindent{\bf Proof.} It is enough to prove two claims.

For $\varphi > 0$, set
$\regul_{s_j,i+2}^\varphi =     \big\{ k \in 2\N \, \cap \, \big[2^{i+2},2^{i+3} \big]:  \max \big\{ \thet_k , \thet_{s_j-16-k}  \big\}  \leq  3/2 +  \varphi \,  \big\}$. 

\noindent{\em Claim $1$.}
The set  $\regul_{s_j,i+2}^{\chi + \epn} \setminus \regnew{s_j}{i+2}{\chi + \epn}$ has cardinality at most $2^{i-4}$.

\noindent{\em Claim $2$.}
$\big\vert \regul_{s_j,i}^{\chi + \epn} \big\vert \geq \tfrac{9}{10} \, 2^{i-1}$.

\noindent{\em Proof of Claim~$1$.} Note that $k \in \regul_{s_j,i+2}^{\chi + \epn} \setminus \regnew{s_j}{i+2}{\chi + \epn}$ implies that either $s_j - 16 - k$ belongs to the union up to index~$i+4$ (in the role of~$i$) of the sets in~(\ref{e.thirtyone}), or $s_\indexc - 16 - k$ belongs to the comparable union of the sets in $(\ref{e.thirtytwo})$. Thus, the cardinality in question is at most $\tfrac{4}{600} 2^{i+3}$, and 
we obtain Claim~$1$.

\noindent{\em Proof of Claim~$2$.}
For $\varphi > 0$, define the set of $\varphi$-high $\thet$ values on dyadic scale $i$,
$$
\macess_{i}^\varphi = \Big\{ k \in 2\N \, \cap \, \big[2^{i},2^{i+1} \big]: \thet_k \geq  3/2  + \varphi \Big\}  \, ,
$$
and also set
$\irr_{s_j,i+2}^\varphi = 2\N \cap \big[ 2^{i+2} , 2^{i+3} \big] \setminus \regul_{s_j,i+2}^\varphi$.
We have then that,
 for any $\varphi > 0$,
$$
 \irr_{s_j,i+2}^\varphi \subseteq  \macess_{i+2}^\varphi \, \cup \, \Big\{ k \in 2\N \cap \big[ 2^{i+2} , 2^{i+3} \big]: \thet_ {s_j-16-k} \geq   3/2 + \varphi \Big\} \, ;
 $$
 the latter event in the right-hand union is a subset of
 $$ 
 2\N \cap \big[ 2^{i+3} - 2^4,2^{i+4} \big] \, \, \cup \, \, \Big\{ k \in 2\N \cap \big[ 2^{i+2},2^{i+3} - 2^4 \big]: s_j-16-k \in \macess_{i+3}^\varphi \cup \, \macess_{i+4}^\varphi \Big\} \, ,
$$
because $2^{i+3} \leq s_j-16-k \leq 2^{i+5} - 2^{i+2} \leq 2^{i+5}$ for $k \in \big[ 2^{i+2},2^{i+3} - 2^4 \big]$. Thus, 
$$
\# \irr_{s_j,i+2}^{\chi + \epn} \leq  \# \macess_{i+2}^{\chi + \epn} +  \# \macess_{i+3}^{\chi + \epn} + \# \macess_{i+4}^{\chi + \epn} + 9 \, ,
$$
 so that~(\ref{e.thirtytwo}) implies $\big\vert \irr_{s_j,i+2}^{\chi + \epn} \big\vert \leq \tfrac{7}{600} \, 2^{i+1} + 9 \leq \tfrac{1}{10}  \, 2^{i+1}$ (for $i \geq 6$). Thus, we obtain Claim~$2$. \qed

\medskip

\noindent{\bf Proof of Lemma~\ref{l.polyepsreg}.}
 The proof is similar to Lemma~\ref{l.polyepsregthet}'s. Now, the probability in question equals the ratio of $\big\vert \bigcup_{l \in  \regnew{s_j}{i+2}{\chi + \epn}} \dubbub_{l,s_j-16-l} \big\vert$ and $p_{s_j}$.

To bound below the numerator, we again seek to  apply Proposition~\ref{p.polyjoin.aboveonehalf} with the role of $i \in \N$ played by $i+2$, but now with $n = s_j$ and $\mathsf{R} =  \regnew{s_j}{i+2}{\chi + \epn}$. Set $\Theta$ in the proposition to be $3/2 + \chi + \eps$. Any element of $\regnew{s_j}{i+2}{\chi + \epn}$ may act as $\kmac \in \mathsf{R}$
in the way that the proposition demands, by the definition of $\regnew{s_j}{i+2}{\chi + \epn}$; (and Lemma~\ref{l.rsize} shows that $\regnew{s_j}{i+2}{\chi + \epn}$ is non-empty if $i \geq 6$).
Applying the proposition and then using Lemma~\ref{l.rsize}, $2^{i+1} \geq 2^{-4} s_j$,
and the definition of  $\regnew{s_j}{i+2}{\chi + \epn}$, we find that 
$$
 \bigg\vert \, \bigcup_{l \in  \regnew{s_j}{i+2}{\chi + \epn}} \dubbub_{l,s_j-16-l} \, \bigg\vert \geq c \, \frac{s_j^{1/2}}{\log s_j} \cdot  \big( \tfrac{9}{10} - \tfrac{1}{32} \big) 2^{-4} s_j  \cdot \mu^{s_j - 16} s_j^{-(3 + 2\chi + 2\eps)}  \, ,   
$$
which is at least a small constant multiple of  $\mu^{s_j} \big(\log s_j\big)^{-1}s_j^{-3/2 - 2\chi - 2\eps}$.

The ratio's denominator is bounded above with the bound  $p_{s_j} \leq \mu^{s_j} \, {s_j}^{-3/2  + \epn}$ which is due to the lower bound in~(\ref{e.tineqs}) with $s_j$ playing the role of $j$. The lemma follows by relabelling~$c > 0$. \qed

\medskip

\subsubsection{Deriving Proposition~\ref{p.similaritynew}.} We begin by respecifying the set $\tau_\gamma$.
For each $\gamma \in \sap_{s_j}$, $1 \leq j \leq s_\indexc$, we now set
$$
\tau_\gamma = \Big\{ k \in \regnew{s_j}{i+2}{\chi + \epnmac} : \gamma \in \dubbub_{k,s_{\indexcj} - 16 - k} \Big\} \, .
$$

Note that $\regnew{s_j}{i+2}{\chi + \epnmac}$ has replaced $\reg_{i+2}$ here. As we comment on the needed changes, it is understood that this replacement is always made.

We reinterpret the quantity $\overline{\theta}$ to be $3/2 + \chi + \epn$. The constant $\conindex > 0$ may again be chosen in order that~(\ref{e.constl}) be satisfied.
Indeed, since $(r_\indexc,s_\indexc)$ is by construction regularity typical on dyadic scale~$i$, we have that $s_\indexc$  is an element of $\mathsf{HCP}_{1/2 + \chi}$ and thus also of $\mathsf{HPN}_{3/2 + \chi + o(1)}$.
 Corollary~\ref{c.globaljoin} is thus applicable with $n = s_\indexc$, so that~(\ref{e.constl}) results.

The set $D_j$ for $j \in [1,\indexc]$ is specified as before, and Lemma~\ref{l.dcub}(1) is also obtained, by the verbatim argument. 

The definition of $\overline\theta$ has been altered in order to maintain the validity of Lemma~\ref{l.dcub}(3). The moment in the proof of this result where the new definition is needed is in deriving the bound 
\begin{equation}\label{e.boundadapt}
p_{s_\indexc}^{-1} p_k p_{s_\indexc - 16 - k} \geq \mu^{-16} s_\indexc^{-2\overline{\theta}}
\end{equation}
for $k \in \regnew{s_j}{i+2}{\chi + \epnmac}$
 in the proof of  Lemma~\ref{l.dcub}(2) (which is needed to prove the third part). The membership of this set by the index $k$ ensures that both $\thet_k$ and $\thet_{s_\indexc - 16 - k}$ are at most $3/2 + \chi + \epnmac = \overline{\theta}$, from which follows~(\ref{e.boundadapt}), and thus Lemma~\ref{l.dcub}(2) and~(3).
 
The weaker form of the inference that we may make in view of the weaker available regularity information is evident in the counterpart to Lemma~\ref{l.compare}. Note the presence of the new factors of $s_\indexc^{-\chi}$ and $s_j^{\chi}$ in the two right-hand sides.
\begin{lemma}\label{l.comparenew}
Recall that $\indexc = 2^{i-4}$.  Let $j \in [1,\indexc]$.
\begin{enumerate}
\item 
For $\phi \in D_j$,
\begin{eqnarray*}
   \psap_{s_{\indexc}} \Big( \Gamma_{[0,r_j]} = \phi  \Big) 
 & \geq &
      \tfrac{1}{400} \, \conindex^{-1} \big( \log s_\indexc \big)^{-1}  s_{\indexc}^{-\chi - 2\epn} \mu^{-16} \\
      & & \qquad \qquad \times \, \, \, \sum_{k \in  \regnew{s_j}{i+2}{\chi + \epn}} \psapleft{k} \Big( \Gamma_{[0,r_j]} = \phi  \Big)
   k^{1/2} p_k \mu^{-k} \, .
\end{eqnarray*}
\item  For $\phi \in \fpart_{r_j}$,
\begin{eqnarray*}
   & & \psap_{s_j} \Big(  \Gamma_{[0,r_j]} = \phi \, , \,  \Gamma \in \dubbub_{k,s_j - 16 - k} \, \, \textrm{for some $k \in  \regnew{s_j}{i+2}{\chi + \epn}$} \Big) \\
   & \leq &  \tfrac{1}{10}  \big( \tfrac{16}{3} \big)^{3/2} \mu^{- 16} \, \cdot s_j^{\chi + 2\epn}  \sum_{k \in \regnew{s_j}{i+2}{\chi + \epnmac}} 
   \psapleft{k}  \Big(  \Gamma_{[0,r_j]} = \phi \Big)  k^{1/2} p_k \mu^{-k} \, .  
\end{eqnarray*}
\end{enumerate}
\end{lemma}
\noindent{\bf Proof (1).} 
The earlier proof runs its course undisturbed until 
$\psap_{s_{\indexc}} \big( \Gamma_{[0,r_j]} = \phi  \big)$
is found to be at least the expression in the double line ending at~(\ref{e.willbesame}). The proof of this part is then completed by noting the bound
$$
 p_{s_\indexc}^{-1} p_{s_\indexc - 16 - k} \geq \mu^{-16 - k} s_\indexc^{-\chi - 2\epn} 
$$
for $k \in  \regnew{s_j}{i+2}{\chi + \epn}$.
This bound is due to $\thet_{s_\indexc - 16 -k} \leq 3/2 + \chi + \epn$ and $\thet_{s_\indexc} \geq 3/2  - \epn$, the latter being the lower bound in~(\ref{e.tineqs}) with $j = s_\indexc$.

\medskip

\noindent{{\bf (2).}} 
The quantity in question is bounded above by~(\ref{e.willbesametwo}). The proof is then completed by invoking the bound
\begin{equation}\label{e.neededbound}
 p_{s_j}^{-1} p_{s_j - 16 - k} \leq \mu^{-16 - k} \big( \tfrac{16}{3} \big)^{3/2} s_j^{\chi + 2\epn}
\end{equation}
for $k \in  \regnew{s_j}{i+2}{\chi + \epn}$.
To derive this bound, note that
$$
 p_{s_j - 16 -k} \leq \mu^{s_j - 16 - k} \big( s_j - 16 - k \big)^{-3/2  + \epn} \leq \mu^{s_j - 16 -k} \big( \tfrac{16}{3} \big)^{3/2} s_j^{-3/2  + \epn} \, , 
$$
where here we used $\thet_{s_j - 16 -k} \geq 3/2  - \epn$, which is due to $k \in   \regnew{s_j}{i+2}{\chi + \epn}$, as well as the harmlessly supposed $\chi + \epn < 3/2$
and the inequality $s_j - 16 - k \geq 3s_j/16$, which we saw at the corresponding moment in the derivation that is being adapted.

We also have  
$$
 p_{s_j} \geq \mu^{s_j} s_j^{-3/2 - \chi - \epn}
$$
due to the upper bound in~(\ref{e.tineqs}) with the role of $j$ played by $s_j$.

The last two displayed inequalities combine to yield~(\ref{e.neededbound}). \qed

\medskip

\noindent{\bf Proof of Proposition~\ref{p.similaritynew}.}
A consequence of Lemmas~\ref{l.dcub}(3) and~\ref{l.comparenew}. \qed

\subsubsection{Proof of Proposition~\ref{p.almostnew}.}
Note that the right-hand side of the bound
$$
  \psap_{s_{\indexc}} \Big( \Gamma_{[0,r_j]} \in \phicl{r_j}{s_j -1}{1/2 + 3\chi + \eclo} \Big) 
  \geq  
  \psap_{s_{\indexc}} \Big( \Gamma_{[0,r_j]} \in \phicl{r_j}{s_j -1}{1/2 + 3\chi + \eclo} 
  \, , \,  
\Gamma_{[0,r_j]} \in D_j  \Big)  
$$
is, by Proposition~\ref{p.similaritynew}, at least the product of $c_1  \big( \log s_\indexc \big)^{-1} s_\indexc^{-2\chi -  4\epn }$ and
$$
 \psap_{s_j} \Big( \Gamma_{[0,r_j]} \in \phicl{r_j}{s_j -1}{1/2 + 3\chi + \eclo} \, , 
    \Gamma \in \dubbub_{k,s_j - 16 - k} \, \, \textrm{for some $k \in \regnew{s_j}{i+2}{\chi + \epn}$} \, , \,  
\Gamma_{[0,r_j]} \in D_j   \Big)
$$
where $c_1 = \tfrac{1}{40} \,  \conindex^{-1}  \big( \tfrac{3}{16} \big)^{3/2}$.

Recall that the index pair $(r_j,s_j)$ is regularity, and thus closing probability, typical on dyadic scale~$i$. We may thus apply 
(\ref{e.acondition}) with $s_j$ playing the role of $j$ and with $a=2$ to find that 
\begin{itemize}
\item
$\psap_{s_j} \Big( \Gamma_{[0,r_j]} \in \phicl{r_j}{s_j - 1}{1/2 + 3\chi + \eclo} \Big) \geq  1 - s_j^{-(2\chi + \eclo - \eroom)}$.
\end{itemize}

Note also that,
\begin{itemize}
\item 
by  Lemma~\ref{l.polyepsreg}, the quantity 
$\psap_{s_j} \Big(   \Gamma \notin \dubbub_{k,s_j - 16 - k} \, \, \textrm{for any $k \in \regnew{s_j}{i+2}{\chi + \epn}$}   \Big)$ is at most $1 \, - \, c \, \big( \log s_j \big)^{-1}  s_j^{-(2\chi + 3\epn)}$,
\item and recall (\ref{e.similaritynew}).
\end{itemize}
Using the latter two inequalities to bound subtractions from the right-hand side of the first bullet point bound, we see that
\begin{eqnarray*}
 & & \psap_{s_{\indexc}} \Big( \Gamma_{[0,r_j]} \in \phicl{r_j}{s_j - 1}{1/2 + 3\chi + \eclo} \Big) \\
 & \geq & c_1  \big( \log s_\indexc \big)^{-1} s_{\indexc}^{-2\chi -  4\epn}  \,  
   \, \bigg(  \, \frac{c}{\log s_j} \, s_j^{-(2\chi + 3\epn)} \, - \, s_j^{-(2\chi + \eclo - \eroom)} \, - \,  64 \mu^{16} s_{\indexcj}^{-10} \, \bigg) \, . 
\end{eqnarray*}

Provided that the index $i \in \N$ is supposed to be sufficiently high, the right-hand side is at least $\tfrac{1}{2} c_1 c \big( \log s_\indexc \big)^{-2} s_{\indexc}^{-4\chi - 7\eps} \geq s_{\indexc}^{-4\chi - 8\eps} = s_{\indexc}^{-\eta + \eps}$ if we insist that $\eclo - \eroom > 3 \epn$ as well as $2\chi + \eclo - \eroom < 10$.
We may set $\eclo = 5\eroom$ and $\eroom = \eps$, because we imposed at the outset that $\chi \in [0,1/{14})$ and a condition on $\epn > 0$ entailing that $\epn \in (0,1/{28})$, and these choices do ensure the stated bounds. Note that $\alpha = 1/2 + 3\chi + \eclo$ results from this choice.

We may now obtain  
Proposition~\ref{p.almostnew}. We set the values of $m$ and $m'$ as we did in proving Proposition~\ref{p.almost}, taking  $m = s_{\indexc}$ and 
$m'$ equal to the constant difference $s_j - r_j$ associated to any element in the sequence $(r_j,s_j)$ from Definition~\ref{d.rjsj.new} (in place of course of Definition~\ref{d.rjsj}).
The bound $\# K \geq 2^{i-4} \geq 2^{-9} m$ holds as in the earlier proof. Note however that we are also claiming that  $m \in \highclose_{1/2 + \chi}$.  
This holds because $m = s_\indexc$ is the second coordinate of a regularity typical index pair for dyadic scale~$i$.  
Proposition~\ref{p.almostnew} is proved. \qed

\subsection{Step three}

At the start of the earlier step three, the index parameters $n$ and $\ell$ were respectively set equal to $m-1$ and $m -m'$, where Proposition~\ref{p.almost} provided $m$ and $m'$. We now make the same choice, with Proposition~\ref{p.almostnew} providing the latter two quantities.

\medskip

\noindent{\bf Proof of Proposition~\ref{p.closingprob}.} The proof of Proposition~\ref{p.thetexist} applies verbatim after recalling that the snake method's parameter $\delta = \beta - \esm - \alpha$ is positive, and noting that the method's index parameter $n$ satisfies $n+1 \in \highclose_{1/2 + \chi}$
because it has been set equal to~$m - 1$ where $m$ is specified in Proposition~\ref{p.almostnew}. \qed

\bibliographystyle{plain}

\bibliography{saw}

\end{document}